\documentclass[11pt]{article}%
 \usepackage[dvips]{color}
\usepackage{amsfonts}
\usepackage{amssymb}
\usepackage{amsmath}
\usepackage{graphicx,pictex}
\usepackage{theorem}
\usepackage[a4paper]{geometry}
\usepackage[toc]{appendix}
\usepackage{hyperref}
\usepackage{makeidx}%
\theoremstyle{break}
\newtheorem{theorem}{Theorem}[section]

\newtheorem{condition}[theorem]{Condition}

\newtheorem{corollary}[theorem]{Corollary}

\newtheorem{definition}[theorem]{Definition}

\newtheorem{lemma}[theorem]{Lemma}

\newtheorem{proposition}[theorem]{Proposition}

\theorembodyfont{\upshape}

\newenvironment{proof}[1][Proof]{\textbf{#1.} }{\ \rule{0.5em}{0.5em}}
\numberwithin{equation}{section}
\makeindex
\begin{document}

\title{Multiscale analysis of exit distributions for random walks in random environments}
\author{Erwin Bolthausen
\and Ofer Zeitouni}
\date{March 19, 2006}
\maketitle

\abstract{We present a multiscale analysis for the exit measures from
large balls in $\mathbb{Z}^d$, $d\geq 3$,
of random walks in certain i.i.d.
random environments which are small perturbations of the fixed
environment corresponding to
simple random walk. Our main assumption is an isotropy assumption
on the law of the environment, introduced by Bricmont and Kupianien.
The analysis is based on propagating estimates on the variational
distance between the exit measure and that of simple random walk,
in addition to estimates on the variational distance between
smoothed versions of these quantities.}

\section{Introduction}

\label{sec-intro} We consider random walks in random environments on
$\mathbb{Z}^{d}$, $d\geq3$, when the environment is a small perturbation of
the fixed environment corresponding to simple random walk. More precisely, let
$\mathcal{P}$ be the set of probability distributions on $\mathbb{Z}^{d},$
charging only neighbors of $0.$ If $\varepsilon\in(0,1/2d),$ we set, with
$\{e_{i}\}_{i=1}^{d}$ denoting the standard basis of $\mathbb{R}^{d}$,%
\begin{equation}
\mathcal{P}_{\varepsilon}\overset{\mathrm{def}}{=}\left\{  q\in\mathcal{P}%
:\left\vert q\left(  \pm e_{i}\right)  -\frac{1}{2d}\right\vert \leq
\varepsilon,~\forall i\right\}  . \label{eq-Pepsdef}%
\end{equation}
$\Omega\overset{\mathrm{def}}{=}\mathcal{P}^{\mathbb{Z}^{d}}$ is equipped with
the natural product $\sigma$-field $\mathcal{F}.$ We call an element
$\omega\in\Omega$ a \textit{random environment}.
%For $A\subset\mathbb{Z}%
%^{d},$ we write $\mathcal{F}_{A}$ for the
%$\sigma$-field generated by the
%projections $\omega\rightarrow\omega_{x},\ x\in A.$
For $\omega\in\Omega,$ and $x\in\mathbb{Z}^{d},$ we consider the transition
probabilities $p_{\omega}\left(  x,y\right)  \overset{\mathrm{def}}{=}%
\omega_{x}\left(  y-x\right)  ,$ if $\left\vert x-y\right\vert =1,$ and
$p_{\omega}\left(  x,y\right)  =0$ otherwise, and construct the random walk
in random environment (RWRE)
$\{X_{n}\}_{n\geq0}$ with initial position $x\in\mathbb{Z}^{d}$ which is,
given the environment $\omega$, the Markov chain with $X_{0}=x$ and transition
probabilities
\[
P_{\omega,x}(X_{n+1}=y|X_{n}=z)=\omega_{z}(y-z)\,.
\]
(By a slight abuse of notation, for consistency with the sequel we also write
$P_{\omega,x}=P_{p_{\omega},x}.$)

We are mainly interested in the case of a \textit{random }$\omega.$ Given a
probability measure $\mu$ on $\mathcal{P},$ we consider the product measure
$\mathbb{P}_{\mu}\overset{\mathrm{def}}{=}\mu^{\otimes\mathbb{Z}^{d}}$ on
$\left(  \Omega,\mathcal{F}\right)  .$ We usually drop the index $\mu$ in
$\mathbb{P}_{\mu}.$ In all that follows we make the following basic assumption

\begin{condition}
\label{Cond_Mu}$\mu$ is invariant under lattice isometries, i.e. $\mu
f^{-1}=\mu$ for any orthogonal mapping $f$ which leaves $\mathbb{Z}^{d}$
invariant, and $\mu\left(  \mathcal{P}_{\varepsilon}\right)  =1$ for some
$\varepsilon\in(0,1/2d)$ which will be specified later\texttt{.}
\end{condition}

The model of RWRE has been studied extensively. We refer to \cite{sznitmanLN}
and \cite{zeitouniLN} for recent surveys. A major open problem is the
determination, for $d>1$, of laws of large numbers and central limit theorems
in full generality (the latter, both under the \textit{quenched} measure, i.e.
for $\mathbb{P}_{\mu}$-almost every $\omega$, and under the \textit{annealed}
measure $\mathbb{P}_{\mu}\otimes P_{x,\omega}$). Although much progress has
been reported in recent years (\cite{BSZ,sznitman1,sznitman2}), a full
understanding of the model has not yet been achieved.

In view of the above state of affairs, attempts have been made to understand
the perturbative behavior of the RWRE, that is the behavior of the RWRE when
$\mu$ is supported on $\mathcal{P}_{\varepsilon}$ and $\varepsilon$ is small.
The first to consider such a perturbative regime were \cite{BK}, who
introduced Condition \ref{Cond_Mu} and showed that in dimension $d\geq3$, for
small enough $\varepsilon$ a quenched CLT holds\footnote{As the examples in
\cite{BSZ} demonstrate, for every $d\geq 7$ and
$\varepsilon>0$ there are measures $\mu$
supported on $\mathcal{P}_{\varepsilon}$, with $\mathbb{E}_{\mu}\left[
\sum_{i=1}^{d} e_i
(q(e_{i})-q(-e_{i}))\right]  =0$, such that $X_{n}/n\to
_{n\to\infty}v\neq0$, $\mathbb{P}_{\mu}$-a.s. One of the goals of Condition
\ref{Cond_Mu} is to prevent such situations from occurring.}. Unfortunately,
the multiscale proof in \cite{BK} is rather difficult, and challenging to
follow. This in turns prompted the derivation, in \cite{SZ}, of an alternative
multiscale approach, in the context of diffusions in random environments. One
expects that the approach of \cite{SZ} could apply to the discrete setup, as
well.

Our goal in this paper is somewhat different: we focus on the exit law of the
RWRE from large balls, and develop a multiscale analysis that allows us to
conclude that the exit law approaches, in a suitable sense, the uniform
measure. Like in \cite{SZ}, the hypothesis propagated involves smoothing. In
\cite{SZ}, this was done using certain H\"{o}lder norms of (rescaled)
transition probabilities. Here, we focus on two ingredients. The
first is  a propagation of
the variational distance between the exit laws of the RWRE from 
balls and those of simple
random walk (which distance
remains small but does not decrease as the scale
increases). The second is the propagation of
the variation distance between
the convolution of the exit law of the RWRE with the exit law
of a simple random walk from a ball of (random) radius, 
and the corresponding convolution of the exit law of simple random walk
with the same smoothing, which 
distance decreases to zero
as scale increases (a precise statement can be 
found in Theorems \ref{Th_Main} and \ref{Th_Main1}; the latter, which is our
main result, provides a local limit law for the exit measures). 
This approach is of a different nature than the one in \cite{SZ} and, we
believe, simpler. In future work we hope to combine our exit law approach with
suitable exit time estimates in order to deduce a (quenched) CLT for the RWRE.

The structure of the article is the following. In the next section, we
introduce our basic notation and state our induction step and our main
results. In Section \ref{Sect_Preliminaries}, we present our basic
perturbation expansion, coarsening scheme for random walks, and auxiliary
estimates for simple random walk. The (rather standard)
proofs of the latter estimates are
presented in the appendices. Section \ref{Sect_Smooth} is devoted to the
propagation of the smoothed estimates, whereas
 Section \ref{Sect_NonSmooth} is
devoted to the propagation of the variation distance 
estimate (the non-smooth estimate). Section \ref{sec-finalpush}
completes the proof of our main result by
using the estimates  of Sections \ref{Sect_Smooth} and
\ref{Sect_NonSmooth}. 
\section{Basic notation and main result\label{Sect_Basic}}

\medskip\noindent\textbf{Sets: }For $x\in\mathbb{R}^{d},$ $\left\vert
x\right\vert $ is the Euclidean norm. If $A,B\subset\mathbb{Z}^{d},$
%$x\in\mathbb{Z}^{d},$ 
we set 
%$d\left(  x,A\right)  \overset{\mathrm{def}}%
%{=}\inf\left\{  \left\vert x-y\right\vert :y\in A\right\}  ,\ 
$d\left(
A,B\right)  \overset{\mathrm{def}}{=}$ $
\inf\left\{  |x-y|:\right.$ $\left.  x\in
A, y\in B\right\}  .$ If $L>0,$ 
we write $V_{L}\overset{\mathrm{def}}{=}%
\{x\in\mathbb{Z}^{d}:\left\vert x\right\vert \leq L\},$ and for $x\in
\mathbb{Z}^{d},$ $V_{L}\left(  x\right)  \overset{\mathrm{def}}{=}x+V_{L}.$ If
$V\subset\mathbb{Z}^{d},$ $\partial V=\{x\in V^c: d(x,V)=1\}$ 
is the outer boundary.
%, i.e. the set of
%points outside $V$ which have a neighbor point in $V$. 
If $x\in V,$ we set
$d_{V}\left(  x\right)  \overset{\mathrm{def}}{=}d\left(  x,\partial V\right)
.$ We also set $d_{L}(x)=L-|x|$ (note that $d_{L}(x)\neq d_{V_{L}}(x)$ with
this convention).
%In case $V=V_{L},$ we simply write $d_{L}\left(  x\right)  .$
For $0\leq a<b\leq L,$ we define 
%the \textquotedblleft shells\textquotedblright%
\begin{equation}
\operatorname*{Shell}\nolimits_{L}\left(  a,b\right)  \overset{\mathrm{def}%
}{=}\left\{  x\in V_{L}:a\leq d_{L}\left(  x\right)  <b\right\}
,\ \operatorname*{Shell}\nolimits_{L}\left(  b\right)  \overset{\mathrm{def}%
}{=}{\operatorname*{Shell}\nolimits}_{L} \left(  0,b\right)  .
\label{Def_Shell}%
\end{equation}

\medskip\noindent\textbf{Functions: }If $F,G$ are functions $\mathbb{Z}%
^{d}\times\mathbb{Z}^{d}\rightarrow\mathbb{R}$ we write $FG$ for the (matrix)
product: $FG\left(  x,y\right)  \overset{\mathrm{def}}{=}\sum_{u}
F\left(
x,u\right)  G\left(  u,y\right)  ,$ provided the right hand side is absolutely
summable. $F^{k}$ is the $k$-th power defined in this way, and $F^{0}\left(
x,y\right)  \overset{\mathrm{def}}{=}\delta_{x,y}.$ We interpret $F$ also as a
kernel, operating from the left on functions $f:\mathbb{Z}^{d}\rightarrow
\mathbb{R},$ by $Ff\left(  x\right)  \overset{\mathrm{def}}{=}\sum_y F\left(
x,y\right)  f\left(  y\right)  $. If $W\subset\mathbb{Z}^{d},$ we use $1_{W}$
not only as the indicator function but, by slight abuse of notation, also to
denote the kernel $\left(  x,y\right)  \rightarrow1_{W}\left(  x\right)
\delta_{x,y}.$

For a function $f:\mathbb{Z}^{d}\rightarrow\mathbb{R},$ $\left\Vert
f\right\Vert _{1}\overset{\mathrm{def}}{=}\sum_{x}\left\vert f\left(
x\right)  \right\vert ,$ and $\left\Vert f\right\Vert _{\infty}\overset
{\mathrm{def}}{=}\sup_{x}\left\vert f\left(  x\right)  \right\vert ,$ as
usual. If \thinspace$F$ is a kernel then, by an abuse of notation, we write
$\left\Vert F\right\Vert _{1}$ for its norm as operator on $L_{\infty},$ i.e.%
\begin{equation}
\left\Vert F\right\Vert _{1}\overset{\mathrm{def}}{=}\sup_{x}\left\Vert
F\left(  x,\cdot\right)  \right\Vert _{1}. \label{Def_OperatorNorm}%
\end{equation}
%We set $\operatorname*{supp}f\overset{\mathrm{def}}{=}\left\{  x:f\left(
%x\right)  \neq0\right\}  .$ If $f,g:\mathbb{Z}^{d}\rightarrow\mathbb{R},$ we
%write $f\ast g$ for the usual convolution:
%$f\ast g(x)=\sum_{y\in\mathbb{Z}^d} f(x-y)g(y)$.

\medskip\noindent\textbf{Transition probabilities: }For transition
probabilities $p=\left(  p\left(  x,y\right)  \right)  _{x,y\in\mathbb{Z}^{d}%
},$ not necessarily nearest neighbor, we write $P_{p,x}$ for the law of a
Markov chain $X_{0}=x,X_{1},\ldots$ 
%on $\mathbb{Z}^{d}$
having $p$ as transition
probabilities.  If $V\subset
\mathbb{Z}^{d},$ $\tau_{V}\overset{\mathrm{def}}{=}\inf\left\{  n\geq
0:X_{n}\notin V\right\}  $ is the first exit time from $V$, and $T_{V}%
\overset{\mathrm{def}}{=}\tau_{V^{c}}$ the first entrance time. We set%
\[
\operatorname*{ex}\nolimits_{\scriptscriptstyle V}\left(  x,z;p\right)
\overset{\mathrm{def}}{=}P_{p,x}\left(  X_{\tau_{V}}=z\right)  .
\]
For $x\in V^{c},$ one has $\operatorname*{ex}\nolimits_{\scriptscriptstyle V}%
\left(  x,z;p\right)  =\delta_{x,z}.$ A special case is the standard simple
random walk $p\left(  x,\pm e_{i}\right)  =1/2d,$ where $e_{1},\ldots,e_{d}%
\in\mathbb{Z}^{d}$ is the standard base. We abbreviate this as $p^{\mathrm{RW}%
},$ and set $P_{x}^{\mathrm{RW}}\overset{\mathrm{def}}{=}P_{x,p^{\mathrm{RW}}%
}.$ Also, exit distributions for the simple random walk are written as
$\pi_{\scriptscriptstyle V}\left(  x,z\right)  \overset{\mathrm{def}}%
{=}\operatorname*{ex}\nolimits_{\scriptscriptstyle V}\left(
x,z;p^{\mathrm{RW}}\right)  .$

We will coarse-grain \textit{nearest-neighbor} transition probabilities $p$ in
the following way. Given $W\subset\mathbb{Z}^{d},$ we choose for any $x\in W$
either a fixed finite 
subset $U_{x}\subset W,$ $x\in U_{x},$ or a probability
distribution $s_{x}$ on such sets. Of course, a fixed choice $U_{x}$ is just a
special choice for the distribution $s_{x},$ namely the one point distribution
on $U_{x}.$

\begin{definition}
\label{Def_CoarseGrainingScheme}A collection $\mathcal{S}=\left(
s_{x}\right)  _{x\in W}$ is called a \textbf{coarse graining scheme}%
\textit{\ }on $W.$ Given such a scheme, and nearest neighbor transition
probabilities $p$, we define the coarse grained transitions by%
\begin{equation}
{p}_{\scriptscriptstyle\mathcal{S},W}^{\mathrm{CG}}\left(  x,\cdot\right)
\overset{\mathrm{def}}{=}\sum_{U:x\in U\subset W}s_{x}\left(  U\right)
\operatorname*{ex}\nolimits_{\scriptscriptstyle U}\left(  x,\cdot;p\right)  .
\label{Smoothing1}%
\end{equation}
\end{definition}
In the case of the standard nearest neighbor random walk, we use the notation
$\pi_{\mathcal{S},W}$ instead of $\left(  {p}^{\mathrm{RW}}\right)
_{\mathcal{S},W}^{\mathrm{CG}}.$ 

Using the Markov property, we have, whenever $W$ is finite,%
\begin{equation}
\operatorname*{ex}\nolimits_{\scriptscriptstyle W}\left(  x,\cdot;p\right)
=\operatorname*{ex}\nolimits_{\scriptscriptstyle W}\left(  x,\cdot
;p_{\scriptscriptstyle\mathcal{S},W}^{\mathrm{CG}}\right)  .
\label{EqualExits}%
\end{equation}

We will choose the coarse-graining scheme in special ways. Fix once for all a
probability density%
\begin{equation}
\varphi:\mathbb{R}^{+}\rightarrow\mathbb{R}^{+},\ \varphi\in C^{\infty
},\ \operatorname*{support}\left(  \varphi\right)  =\left[  1,2\right]  .
\label{SmootingFunction}%
\end{equation}
If $m\in\mathbb{R}^{+},$ the rescaled density is defined by $\varphi
_{m}\left(  t\right)  \overset{\mathrm{def}}{=}\left(  1/m\right)
\varphi\left(  t/m\right)  .$ The image measure of $\varphi_{m}\left(
t\right)  dt$ under the mapping $t\rightarrow V_{t}\left(  x\right)  \cap W$
defines a probability distribution on subsets of $W$ containing $x$.
We may also choose $m$ to depend on $x,$ i.e. consider a field ${\Psi}=\left(
m_{x}\right)  _{x\in W}$ of positive real numbers on $W.$ Such a field then
defines via the above scheme coarse grained transition probabilities, which by
a slight abuse of notation we denote as $p_{\scriptscriptstyle{\Psi}%
,W}^{\mathrm{CG}}.$ In case $W=\mathbb{Z}^{d},$ we simply drop $W$ in the
notation. In case $p$ is the standard nearest neighbor random walk, we write
$\hat{\pi}_{{\Psi}}$ instead of $p_{\scriptscriptstyle{\Psi}}^{\mathrm{CG}}.$

\medskip\noindent\textbf{The random environment: } We recall from the
introduction the notation $\mathcal{P}_{\varepsilon}$, $\Omega$, $p_{\omega
}\left(  x,y\right)  $, and the natural product $\sigma$-field $\mathcal{F}.$
For $A\subset\mathbb{Z}^{d},$ we write $\mathcal{F}_{A}=\sigma(\omega_x: x\in A)$.
%for the $\sigma
%$-field generated by the projections $\omega\rightarrow\omega_{x},\ x\in A.$
We also recall the probability measure $\mu$ on $\mathcal{P},$ the product
measure $\mathbb{P}_{\mu}$, and Condition \ref{Cond_Mu}, which is assumed throughout.

For a random environment $\omega\in\Omega$,
%a \textit{random environment,} and in this random case,
we typically write $\Pi_{\scriptscriptstyle V,\omega}\overset{\mathrm{def}}%
{=}\operatorname*{ex}\nolimits_{\scriptscriptstyle
V}\left(  \cdot,\cdot;p_{\omega}\right)  $ and occasionally drop $\omega$ in
the notation. So $\Pi_{V}$ should always be understood as a \textit{random}
exit distribution. We will also use $\hat{\Pi}_{\mathcal{S},W}$ for $\left(
p_{\omega}\right)  _{\mathcal{S},W}^{\mathrm{CG}}.$
%but we typically drop $\omega.$

For $x\in\mathbb{Z}^{d},$ $L>0,$ and ${\Psi}:\partial V_{L}\left(  x\right)
\rightarrow\mathbb{R}^{+}$, we define the random variables%
\begin{equation}
D_{L,{\Psi}}\left(  x\right)  \overset{\mathrm{def}}{=}\left\Vert \left(
\left[  \Pi_{\scriptscriptstyle V_{L}\left(  x\right)  }-\pi
_{\scriptscriptstyle V_{L}\left(  x\right)  }\right]  \hat{\pi}_{{\Psi}%
}\right)  \left(  x,\cdot\right)  \right\Vert _{1}, \label{Def_DL}%
\end{equation}%
\begin{equation}
\label{eq-280905a}D_{L,0}\left(  x\right)  \overset{\mathrm{def}}%
{=}\left\Vert \Pi_{V_{L}\left(  x\right)  }\left(  x,\cdot\right)  -\pi
_{V_{L}\left(  x\right)  }\left(  x,\cdot\right)  \right\Vert _{1},
\end{equation}
and with $\delta>0,$ we set%
\begin{align*}
&  b_{i}\left(  L,{\Psi},\delta\right)  \overset{\mathrm{def}}{=}%
\mathbb{P}\left(  \left(  \log L\right)  ^{-9+\frac{9(i-1)}{4}}<
D_{L,\Psi}\left(
0\right)
\leq\left(  \log L\right)  ^{-9+\frac{9i}{4}},\ D_{L,0}\left(  0\right)  \leq
\delta\right)\,,\; i=1,2,3, \\
%&  b_{2}\left(  L,{\Psi},\delta\right)  \overset{\mathrm{def}}{=}%
%\mathbb{P}\left(  \left(  \log L\right)  ^{-6.75}< D_{L,\Psi}\left(  0\right)
%\leq\left(  \log L\right)  ^{-4.5},\ D_{L}^{0}\left(  0\right)  \leq
%\delta\right) \\
%&  b_{3}\left(  L,{\Psi},\delta\right)  \overset{\mathrm{def}}{=}%
%\mathbb{P}\left(  \left(  \log L\right)  ^{-4.5}< D_{L,\Psi}\left(  0\right)
%\leq\left(  \log L\right)  ^{-2.25},\ D_{L}^{0}\left(  0\right)  \leq
%\delta\right) \\
&  b_{4}\left(  L,{\Psi},\delta\right)  \overset{\mathrm{def}}{=}%
\mathbb{P}\left(  \left\{  \left(  \log L\right)  ^{-2.25}< D_{L,\Psi}\left(
0\right)  \right\}  \cup\left\{  D_{L,0}\left(  0\right)  >\delta\right\}
\right)\,, \\
&  b\left(  L,{\Psi},\delta\right)  \overset{\mathrm{def}}{=}\sum_{i=1}^4 
b_{i}\left(
L,{\Psi},\delta\right) \,.
\end{align*}

We write $\mathcal{M}_{L}$ for the set of functions ${\Psi}:\partial
V_{L}\rightarrow\left[  L/2,2L\right]  $ which are restrictions of functions
defined on $\left\{  x\in\mathbb{R}^{d}:L/2\leq\left\vert x\right\vert
\leq2L\right\}  $ that have smooth third derivatives bounded by $10L^{-2}$ and
fourth derivatives bounded by $10L^{-3}$. We write $\Psi_t=(m_x=t)_{x\in
\mathbb{Z}^d}$ for the
coarse-graining scheme that consists  of constant coarse-graining
at scale $t$. Of course,
 $\Psi_t\in \mathcal{M}_L$ for all $t,L$.
\begin{condition}
\label{Cond_Main}Let $L_{1}\in\mathbb{N},$ and $\delta>0.$ We say that
condition $\operatorname*{Cond}\left(  \delta,L_{1}\right)  $ holds provided
that for all $L\leq L_{1},$ and for all ${\Psi}\in\mathcal{M}_{L}$,%
\begin{equation}
b_{i}\left(  L,{\Psi},\delta\right)  \leq\frac{1}{4}\exp\left[  -\left(
1-\left(  4-i\right)  /13\right)  \left(  \log L\right)  ^{2}\right]
,\ i=1,2,3,4. \label{BoundBad}%
\end{equation}
\end{condition}
In particular, if $\operatorname*{Cond}\left(  \delta,L_{1}\right)  $ is
satisfied, then for any $L\leq L_{1},$ and any ${\Psi}\in\mathcal{M}_{L}$,%
\begin{equation}
\mathbb{P}\left(  \{D_{L,0}(0)>\delta\}\cup\{D_{L,\Psi}(0)>(\log L)^{-9}\}
\right)  \leq\exp\left[  - \frac{10}{13} \left(  \log L\right)  ^{2}\right]
\label{eq-220805a}%
\end{equation}

Our main technical inductive result is

\begin{proposition}
\label{Prop_Main} There exist $\delta_{0}>0$
%,\ L_{0}\in\mathbb{N},$ such that
such that for all $\delta\in(0,\delta_{0}]$ there exists $\varepsilon
_{0}\left(  \delta\right)  $ and $L_{0}\in\mathbb{N}$ such that if
$\varepsilon\leq\varepsilon_{0},$ $L_{1}\geq L_{0},$ and $\mu$ is such that
Condition \ref{Cond_Mu} holds for $\varepsilon,$ then%
\[
\operatorname*{Cond}\left(  \delta,L_{1}\right)  \Longrightarrow
\operatorname*{Cond}\left(  \delta,L_{1}\left(  \log L_{1}\right)
^{2}\right)  .
\]

\end{proposition}

%\begin{remark}
Given $L_{0},\delta_{0},$ we can always choose $\varepsilon_{0}$ so small that
if Condition \ref{Cond_Mu} is satisfied with $\varepsilon\leq
\varepsilon_{0},$ then
$\operatorname*{Cond}\left(  \delta_{0},L_{0}\right)  $ holds trivially. 
Proposition \ref{Prop_Main} then implies that for any
$\delta<\delta_{0}$, there exists 
$\varepsilon_{0}=\varepsilon_{0}(\delta)$ small enough such that
if Condition \ref{Cond_Mu} is satisfied with $\varepsilon\leq 
\varepsilon_{0},$ then
$\operatorname*{Cond}\left(  \delta,L\right)  $ 
holds for all $L$.
%$\varepsilon
%\leq\varepsilon_{0},$ and all $L$. 
In particular, one obtains immediately from
Proposition \ref{Prop_Main} the following theorem (recall that
$\Psi_t$ denotes constant coarse-graining at scale $t$).
%, which is the main result of
%this paper.

\begin{theorem}
\label{Th_Main} For each $\delta>0$
there exists an $\varepsilon_{0}=\varepsilon_0(\delta)>0$ such that if 
Condition \ref{Cond_Mu} is satisfied with $\varepsilon\leq 
\varepsilon_{0},$ then
for any integer $r\geq 0$,
\[
\limsup_{L\rightarrow\infty}L^r 
b\left(  L,\Psi_L,\delta\right)  =0\,.
\]
%where $m_{L}$ denotes the element of $\mathcal{M}_{L}$ that consists of
%constant smoothing at scale $L$, that is $m_x=L$ for all $x$. 
\end{theorem}
%The Borel-Cantelli lemma then implies that under the
%conditions of Theorem \ref{Th_Main}, $D_{L,m_L}(0)\to_{L\to\infty} 0$,
%$\mathbb{P}_{\mu}$-a.s.

%\end{remark}

Our induction will also provide 
the following theorem which is the
main result of our paper. It provides a local limit theorem for the exit law.
\begin{theorem}
\label{Th_Main1} There exists $\varepsilon_{0}>0,$ such that
if 
Condition \ref{Cond_Mu} is satisfied with $\varepsilon\leq 
\varepsilon_{0},$ then
for any $\delta>0$,
and for any integer $r\geq 0$,
\[
\lim_{t\rightarrow\infty}\limsup_{L\rightarrow\infty}L^r b\left(  L,\Psi_t,\delta
\right)  =0\,.
\]
\end{theorem}
The Borel-Cantelli lemma then implies that under the
conditions of Theorem \ref{Th_Main}, 
$$\limsup_{L\to\infty}
D_{L,\Psi_t}(0)\leq  c_t\,,\quad
\mathbb{P}_{\mu}-a.s.,$$
where $c_t$ is a constant such that 
$c_t\to_{t\to\infty} 0$.
%\begin{remark}
%\label{rem-281205}
%In fact, the same proof as in the induction step in 
%Proposition \ref{Prop_Main}, with 
%$\Psi$ in  $\hat{\pi}_{{\Psi}}$ replaced by $m_{L/(\log L)^8}$,
%yields the following variant of Theorem \ref{Th_Main}:
%For each $\delta<\delta_{0}$ 
%there exists an $\varepsilon_{0}>0$ such that if
%Condition \ref{Cond_Mu} is satisfied with $\varepsilon_{0},$ then
%\[
%\limsup_{L\rightarrow\infty}b\left(  L,m_{L/(\log L)^8},\delta\right) 
% =0\,.
%\]
%In particular, a mesoscopic smoothing is sufficient in order
%to get the variation distance between the 
%exit measure of simple random walk and the RWRE to converge to $0$.
%\end{remark}

A remark about the wording which we use.
When we say that something
holds for \textquotedblleft large enough $L$\textquotedblright, we mean that
there exists $L_{0},$ \textit{depending only on the dimension}, such that the
statement holds for $L\geq L_{0}.$ We emphasize that $L_{0}$ then \textit{does
not depend on }$\varepsilon$.
%or $\delta$.

We write $C$ for a generic positive constant, not necessarily the same at
different occurrences. $C$ may depend on the dimension $d$ of the lattice, but
on nothing else, except when indicated explicitly. Other constants, such as
$c_{0},c_{1},\bar c, k_{0},K,C_{1}$ etc., follow the same convention
concerning what they depend on ($d$ only, unless explicitly stated
otherwise!), but their value is fixed throughout the paper and does not change
from line to line.

\section{Preliminaries \label{Sect_Preliminaries}}

\subsection{The perturbation expansion
\label{Subsect_Perturbation}}

Let $p=\left(  p\left(  x,y\right)  \right)  _{x,y\in\mathbb{Z}^{d}}$ be a
Markovian transition kernel on $\mathbb{Z}^{d},$ not necessarily nearest
neighbor, but of finite range, and let $V\subset\subset\mathbb{Z}^{d}.$ The
Green kernel on $V$ with respect to $p$ is defined by%
\[
g_{\scriptscriptstyle V}\left(  p\right)  \left(  x,y\right)  \overset
{\mathrm{def}}{=}\sum_{k\geq0}\left(  1_{\scriptscriptstyle V}p\right)
^{k}\left(  x,y\right)  .
\]
Evidently, if $z\notin V,$ then
\begin{equation}
g_{\scriptscriptstyle V}\left(  p\right)  \left(  \cdot,z\right)
=\operatorname*{ex}\nolimits_{\scriptscriptstyle V}\left(  \cdot,z;p\right)  .
\label{Green&Exit}%
\end{equation}
If $p,q$ are two transition kernels, write $\Delta_{p,q}=
1_V(p-q)$. The resolvent equation gives for every
$n\in\mathbb{N}$,%
\begin{align}
\label{Pert1}
&g_{\scriptscriptstyle V}\left(  p\right)   -g_{\scriptscriptstyle V}\left(
q\right) 
=g_{\scriptscriptstyle V}\left(  q\right)  \Delta_{p,q}
%1_{\scriptscriptstyle V}%
%\left(  p-q\right)  
g_{\scriptscriptstyle V}\left(  p\right) 
\nonumber\\
&  
=\sum_{k=1}^{n-1}\left[
g_{\scriptscriptstyle V}\left(  q\right)  \Delta_{p,q}
%1_{\scriptscriptstyle V}\left(
%p-q\right)  
\right]  ^{k}g_{\scriptscriptstyle V}\left(  q\right)
+\left[  g_{\scriptscriptstyle V}\left(  q\right)  \Delta_{p,q}
\right]  ^{n}g_{\scriptscriptstyle V}\left(  p\right)
=\sum_{k=1}^{\infty}\left[
g_{\scriptscriptstyle V}\left(  q\right)  \Delta_{p,q}
\right]  ^{k}g_{\scriptscriptstyle V}\left(  q\right)  ,
\end{align}
assuming convergence of the infinite series, which will always be trivial in
cases of interest to us, due to ellipticity and $V$ being finite.
We will occasionally slightly modify the above expansion, but the basis is
always the first equality in (\ref{Pert1}).

\subsection{The coarse graining 
schemes on $V_{L}\label{SubSect_Smoothing}$}

Our proof of Theorems \ref{Th_Main} and \ref{Th_Main1} is based
on a couple of explicit coarse graining schemes, whose definitions we now present.
Set
\begin{equation}
r\left(  L\right)  \overset{\mathrm{def}}{=}L/\left(  \log L\right)
^{10},\ s\left(  L\right)  \overset{\mathrm{def}}{=}L/\left(  \log L\right)
^{3} \,, \ 
 \operatorname*{Sh}\nolimits_{L}
\overset{\mathrm{def}}{=}
 \operatorname*{Shell}\nolimits_{L}\left(  r(L)\right)%
\label{Def_sL&rL}\,, \ 
\end{equation}
and
%$K$ was chosen independently of $\gamma,$ and thus we can set%
\begin{equation}
\gamma\overset{\mathrm{def}}{=}\min\left(  \frac{1}{10}, \frac{1}{2}\left(
1-\left(  \frac{2}{3}\right)  ^{1/\left(  d-1\right)  }\right)  \right)  .
\label{Def_Gamma}%
\end{equation}
%(The reason for the middle term in the right hand side of definition
%(\ref{Def_Gamma}) will only become clear in Section \ref{Sect_NonSmooth}, see
%(\ref{eq-260805b}).)

We fix a $C^{\infty}$-function $h:\mathbb{R}^{+}\rightarrow\mathbb{R}^{+},$
which satisfies $h\left(  u\right)  =u$ for $u\leq1/2,$ $h\left(  u\right)
=1$ for $u\geq2,$ and is strictly monotone and concave on $\left(
1/2,2\right)  .$ For $x\in V_{\scriptscriptstyle L},$ we set%
\begin{equation}
h_{\scriptscriptstyle L}\left(  x\right)  \overset{\mathrm{def}}{=}{\gamma
s\left(  L\right)  }h\left(  \frac{d_{L}\left(  x\right)  }{s\left(  L\right)
}\right)  . \label{Def_hL}%
\end{equation}
%The constant $\gamma>0$ should be thought
%to be small, and will be specified
%later, c.f. (\ref{Def_Gamma})
%($\gamma$ will depend only on $d$).
Remark that for $d_{L}\left(  x\right)  \geq2s\left(  L\right)  ,$ we have
$h_{\scriptscriptstyle L}\left(  x\right)  =\gamma s\left(  L\right)  .$

\begin{lemma}
\label{lem-130705} Fix $\delta_{1}>0$. Then, there is a constant $\bar
k_{0}=\bar k_{0}(\delta_{1})$ such that if $k\geq\bar k_{0}(\delta_{1})$,
and $\Delta(x,y)=\Pi_{V_{kr(L)}(x)}(x,y)-\pi_{V_{kr(L)}(x)}(x,y)$ then
for all $L$ large, if for some $\delta>0$, $d_{L}\left(  x\right)  \leq
r\left(  L\right)  $ and
$D_{kr\left(  L\right),0  }\left(  x\right)
\leq\delta,$ then
\begin{equation}
\sum_{y\in V_{L}\cap 
 \operatorname*{Sh}\nolimits_{L}
}\left\vert \Delta(x,y)\right\vert
\leq\delta+\delta_{1}\,. \label{eq-130705a}%
\end{equation}
\end{lemma}
\begin{proof}
Fix $k$. We have
\begin{align*}
  \sum_{y\in V_{L}\cap
 \operatorname*{Sh}\nolimits_{L}
}|\Delta(x,y)|
&  \leq\Pi_{V_{kr\left(  L\right)  }\left(  x\right)  }\left(  x,V_{L}%
\cap\operatorname*{Sh}\nolimits_{L}
\right)  +\pi_{V_{kr\left(  L\right)  }\left(  x\right)  }\left(  x,V_{L}%
\cap\operatorname*{Sh}\nolimits_{L}
\right) \\
&  \leq\delta+2\pi_{V_{kr\left(  L\right)  }\left(  x\right)  }\left(
x,V_{L}\cap\operatorname*{Sh}\nolimits_{L}
  \right)  .
\end{align*}
Choosing $k$ large enough completes the proof.
\end{proof}

We can now define our coarse graining schemes on $V_L$. 
The first will depend on a
constant $k_{0}>1$ that will be chosen below, based on some a-priori
estimates concerning simple random walk, see (\ref{k_0_large}).
\begin{definition}
\label{Def_SmoothingScheme}
\begin{enumerate}
\item[a)] The coarse graining scheme $\mathcal{S}_{1} =\mathcal{S}_{1,L,k_{0}}
=\left(  s_{x}\right)  _{x\in V_{L}}$ is defined
for $d_{L}\left(
x\right)  \leq r\left(  L\right)  $ by $s_{x}=\delta_{V_{k_{0}r\left(
L\right)}\cap V_L  }$, 
i.e. for such an $x,$ the coarse graining is done by choosing
the exit distribution from $V_{k_{0}r\left(  L\right)  }\left(  x\right)  \cap
V_{L}.$ For $d_{L}\left(  x\right)  >r\left(  L\right)  $, we take
$m_x=h_L(x)$ and define $s_{x}$
according to the description following 
(\ref{SmootingFunction}).
%\overset{\mathrm{def}}{=}\varphi_{h_{L} \left(
%x\right)  }\left(  t\right)  dt.$
\item[b)] The coarse graining scheme $\mathcal{S}_{2}= \mathcal{S}_{2,L}=\left(
s_{x}\right)  _{x\in V_{L}}$ is defined for all $x$ 
%with $d_L(x)\geq 
%1/2\gamma$ 
by 
$m_x=h_L(x)$.
%, and for all $x$ with $d_L(x)<  1/2\gamma$ by 
%$s_x=\delta_{x}$ (that is, no coarse graining at all).
%by $s_{x}\left(  dt\right)
%\overset{\mathrm{def}}{=}\varphi_{h_{L}\left(  x\right)  }\left(  t\right)
%dt$ for all $x.$
\end{enumerate}
\end{definition}
We will need the second scheme only in Section \ref{Sect_NonSmooth}, when
propagating  the part of the
estimate $b_{4}(L,\Psi,\delta)$ involving the expression $D_{L,0}(x)$ of
(\ref{eq-280905a}).  Note that under $\mathcal{S}_2$, if $d_L(x)<1/2\gamma$
then there is no coarse graining at all, i.e. $s_x=\delta_x$.

We write $\rho_{i,L}\left(  x\right)  $  for the range of the coarse
graining scheme at $x$ in scheme $i,$ $i=1,2$, i.e.
\begin{equation}
\rho_{1,L}\left(  x\right)  \overset{\mathrm{def}}{=}\left\{
\begin{array}
[c]{ll}%
k_{0}r\left(  L\right)  & \mathrm{for\ }d_{L}\left(  x\right)  \leq r\left(
L\right) \\
2h_{L}\left(  x\right)  & \mathrm{for\ }r\left(  L\right)  <d_{L}\left(
x\right)
\end{array}
\right.\,, \
\rho_{2,L}=
2h_{L}\left(  x\right) \,. 
\label{Def_Rho1}%
\end{equation}
%Up to Section \ref{Sect_NonSmooth}, we therefore only work
%with $\mathcal{S}_{1}.$
\begin{figure}[t]
\begin{picture}(10,200)(-80,0)\input{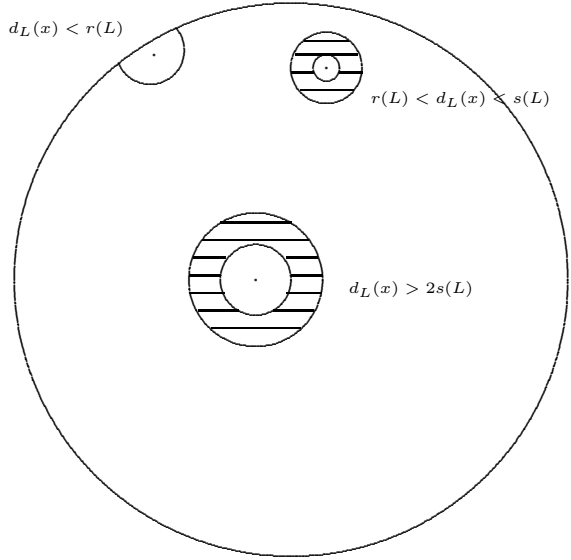}
\end{picture}
\caption{The coarse graining scheme $\mathcal{S}_{1}$}
\end{figure}

\subsection{Estimates on exit distributions and the Green's
function\label{Subsect_Exit&Green}}

For notational convenience, we write $\pi_{L}$ instead of $\pi_{V_{L}},$ and
similarly in other expressions. For instance, we write $\tau_{L}$ instead of
$\tau_{V_{L}}.$

\begin{lemma}
\label{Le_Lawler_Exit}
\begin{enumerate}
\item[a)] For $x\in \partial V_L$,
\[
\frac{1}{C}L^{-d+1}\leq\pi_{L}\left(  x\right)  \leq CL^{-d+1}.
\]
\item[b)] Let $x$ be a vector of unit length in $\mathbb{R}^{d},$ let
$0<\theta<1,$ and define the cone $C_{\theta}\left(  x\right)  \overset
{\mathrm{def}}{=}\left\{  y\in\mathbb{Z}^{d}:\left\langle y,x\right\rangle
\geq\left(  1-\theta\right)  \left\vert y\right\vert \right\}  .$ For any
$\theta,$ there exists $\eta\left(  \theta\right)  >0,$ such that for all
$L$ large enough, and all $x$%
\begin{equation}
\pi_{L}\left(  0,C_{\theta}\left(  x\right)  \right)  \geq\eta\left(
\theta\right)  . \label{ConeEst}%
\end{equation}
\item[c)] Let $0<l<L,$ and $x\in\mathbb{Z}^{d}$ satisfy $l<\left\vert
x\right\vert <L.$ Then%
\[
P_{x}^{\mathrm{RW}}\left(  \tau_{L}<T_{V_{l}}\right)  =\frac{ l^{-d+2}%
-\left\vert x\right\vert ^{-d+2}+ O\left(  l^{-d+1}\right)  }{l^{-d+2}%
-L^{-d+2}}%
\]
\end{enumerate}
\end{lemma}
\begin{proof}
a) is Lemma 1.7.4 of \cite{Lawler}. b) is immediate from a). c) is Proposition
1.5.10 of \cite{Lawler}.
\end{proof}

We will repeatedly make use of the following lemma.
\begin{lemma}
\label{Le_MainExit}Assume $x,y\in V_{L},$ $1\leq a\leq5d_{L}\left(  y\right)
,$ $x\notin V_{2a}\left(  y\right)  .$ Then%
\begin{equation}
P_{x}\left(  T_{V_{a}\left(  y\right)  }<\tau_{V_{L}} \right)  \leq
C\frac{a^{d-2}d_{L}\left(  y\right)  d_{L}\left(  x\right)  }{\left\vert
x-y\right\vert ^{d}} \label{eq-080605gg}%
\end{equation}
\end{lemma}
The proof will be given in Appendix \ref{App_A}.

We will need a corresponding result for the Brownian motion. We write $\pi
_{L}^{\mathrm{BM}}(y,dy^{\prime})$ for the exit distribution of the Brownian
motion from the ball $C_{L}$ of radius $L$ in $\mathbb{R}^{d}.$ The following
lemma is an easy consequence of the Poisson formula, see \cite[(1.43)]{Lawler}.

\begin{lemma}
\label{Le_ExitsBM}For any $y\in C_{L}$, it holds that
\begin{equation}
\frac{C^{-1}d(y,\partial C_{L})}{|y-y^{\prime}|^{d}}\leq\frac{\pi
_{L}^{\mathrm{BM}}(y,dy^{\prime})}{dy^{\prime}}\leq\frac{Cd(y,\partial C_{L}%
)}{|y-y^{\prime}|^{d}}, \label{eq-200305ff}%
\end{equation}
where $dy^{\prime}$ is the surface measure on $\partial C_{L}.$
\end{lemma}

We will also
need a comparison between smoothed exit distribution of the random
walk, and that of Brownian motion. Given $L>0,$ and ${\Psi}\in\mathcal{M}%
_{L},$ let
\begin{equation}
\phi_{L,{\Psi}}\overset{\mathrm{def}}{=}\pi_{L}\hat{\pi}_{{\Psi}}.
\label{Def_PhiLM}%
\end{equation}
We consider also the corresponding Brownian kernel on $\mathbb{R}^{d}$,%
\begin{equation}
	\label{170306bb}
\phi_{L,{\Psi}}^{\mathrm{BM}}\left(  y,dz\right)  \overset{\mathrm{def}}%
{=}\int_{\partial C_{L}\left(  0\right)  }\pi_{C_{L}\left(  0\right)
}^{\mathrm{BM}}\left(  y,dw\right)  \int\pi_{C_{t}\left(  w\right)
}^{\mathrm{BM}}\left(  w,dz\right)  \varphi_{m_{w}}\left(  t\right)  dt,
\end{equation}
where $\Psi=\left(  m_{w}\right)  ,$ and where we write $\phi_{L,{\Psi}%
}^{\mathrm{BM}}\left(  y,z\right)  $ for the density of
$\phi_{L,{\Psi}}^{\mathrm{BM}}\left(  y,dz\right)$
with respect to
$d$-dimensional Lebesgue measure.
\begin{lemma}
\label{Le_Approx_Phi_by_BM} There exists a constant $C$ such that for $L>0,$
and ${\Psi}\in\mathcal{M}_{L},$ we have%
\[
\sup_{y\in V_{L}}\sup_{z\in\mathbb{Z}^{d}}\left\vert \phi_{L,{\Psi}}\left(
y,z\right)  -\phi_{L,{\Psi}}^{\mathrm{BM}}\left(  y,z\right)  \right\vert \leq
CL^{-d-1/5}\,.%
\]
\end{lemma}
\begin{lemma}
\label{Le_ThirdDerivative} There exists a constant $C$ such that for $L>0$ and
${\Psi}\in\mathcal{M}_{L},$ we have%
\[
\sup_{y,z}\left\Vert \partial_{y}^{i}\phi_{L,{\Psi}}^{\mathrm{BM}}\left(
y,z\right)  \right\Vert \leq CL^{-d-i}\,, i=1,2,3\,.%
\]
\end{lemma}
The proofs of these two lemmas are again in Appendix \ref{App_A}.

We can draw two immediate conclusions from these results:
\begin{proposition}
\label{Prop_LipshitzPhi}
\begin{itemize}
\item[a)] Let $y,y^{\prime}$ be in $V_{L}$, and $\Psi\in\mathcal{M}_{L}.$ Then%
\begin{equation}
\left\vert \phi_{L,{\Psi}}\left(  y,z\right)  -\phi_{L,{\Psi}}\left(
y^{\prime},z\right)  \right\vert \leq C\left(  L^{-d-1/5}+\left\vert
y-y^{\prime}\right\vert L^{-d-1}\right)  . \label{OneDerivative}%
\end{equation}
\item[b)] Let $x\in V_{L},$ and $l$ be such that $V_{l}\left(  x\right)
\subset V_{L}.$ Consider a signed measure $\mu$ on $V_{l}$ with total mass $0$
and total variation norm $|\mu|$,
which is invariant under lattice isometries. Then%
\begin{equation}
\left\vert \sum\nolimits_{y}\mu\left(  y-x\right)  \phi_{L,{\Psi}}\left(
y,z\right)  \right\vert \leq C\left\vert \mu\right\vert \left(  L^{-d-1/5}%
+\left(  \frac{l}{L}\right)  ^{3}L^{-d}\right)  \,. \label{ThreeDerivatives}%
\end{equation}
\end{itemize}
\end{proposition}
\begin{proof}
[Proof of Proposition \ref{Prop_LipshitzPhi}]a) is immediate from Lemmas
\ref{Le_Approx_Phi_by_BM} and \ref{Le_ThirdDerivative}.
As for b), we get from Lemma \ref{Le_Approx_Phi_by_BM} that
%the same lemmas
\[
\left\vert \sum\nolimits_{y}\mu\left(  y-x\right)  \phi_{L,{\Psi}}\left(
y,z\right)  -\sum\nolimits_{y}\mu\left(  y-x\right)  \phi_{L,{\Psi}%
}^{\mathrm{BM}}\left(  y,z\right)  \right\vert \leq C\left\vert \mu\right\vert
L^{-d-1/5}\,,
\]%
while
\begin{align}
& \sum\nolimits_{y}\mu\left(  y-x\right)  \phi_{L,{\Psi}}^{\mathrm{BM}}\left(
y,z\right)     =\sum\nolimits_{y}\mu\left(  y-x\right)  \left[  \phi
_{L,{\Psi}}^{\mathrm{BM}}\left(  y,z\right)  -\phi_{L,{\Psi}}^{\mathrm{BM}%
}\left(  x,z\right)  \right] \nonumber\\
&\quad  =\sum\nolimits_{y}\mu\left(  y-x\right)  \partial_{x}\phi_{L,{\Psi}%
}^{\mathrm{BM}}\left(  x,z\right)  \left[  y-x\right] \label{Eq_Harmonic}\\
&\quad 
\quad +\frac{1}{2}\sum\nolimits_{y}\mu\left(  y-x\right)  \partial_{x}^{2}%
\phi_{L,{\Psi}}^{\mathrm{BM}}\left(  x,z\right)  \left[  y-x,y-x\right]
 +R\left(  \mu,x,z\right)  ,\nonumber
\end{align}
where, due to Lemma \ref{Le_ThirdDerivative},
\begin{equation}
\left\vert R\left(  \mu,x,z\right)  \right\vert \leq C\left\vert
\mu\right\vert \left(  \frac{l}{L}\right)  ^{3}L^{-d} \label{eq-080605a}%
\end{equation}
uniformly in $x$ and $z$,
and  $\partial^{k}F\left[  u_{1},\ldots,u_{k}\right]  $
denotes the $k$-th derivative of a function $F$
in directions $u_{1},\ldots,u_{k}$. 
The first
summand on the right hand side of (\ref{Eq_Harmonic}) vanishes 
because $\mu$
has mean $0.$ The second vanishes because by the invariance under lattice
isometry of $\mu,$ the summand 
involves only the Laplacian of $\phi_{L,{\Psi}%
}^{\mathrm{BM}}\left(  \cdot,z\right)  ,$ which in turn 
vanishes because of
harmonicity of $\pi_{C_{L}\left(  0\right)  }^{\mathrm{BM}}\left(
x,\cdot\right)  $ in the $x$-variable. 
%The estimate (\ref{eq-080605a}) follows
%from Lemma \ref{Le_ThirdDerivative}. 
The proof of the proposition is complete.
\end{proof}

The next lemma gives a-priori estimates for coarse-grained walks. We use
$\hat{\pi}_{L}^{(i)}$, $i=1,2$, to denote the transitions of the coarse
grained random walk that uses the coarse graining $\mathcal{S}_{i}$, and
$\hat{g}_{L}^{(i)}$ to denote the corresponding Green's function. Note that
these quantities all depend on $L$ and $k_{0}$, but we suppress these from the
notation. Recall that $
\operatorname*{Sh} \nolimits_{L}
=
\operatorname*{Shell}%
\nolimits_{L}\left(  r\left(  L\right)  \right)$, c.f. (\ref{Def_sL&rL}).
\begin{lemma}
\label{Cor_Green} There exists a constant $C$ (independent of $k_{0}$!) such that:
\begin{enumerate}
\item[a)]
\[
\sup_{x\in V_{L}}\hat{g}_{L}^{(1)}\left(  x,
\operatorname*{Sh} \nolimits_{L}
\right)  \leq C.
\]
\item[b)] If $i=1$ and $r\left(  L\right)  \leq a\leq3s\left(  L\right)  $ or
$i=2$ and $a\leq3s\left(  L\right)  $ then,
\[
\sup_{x\in V_{L}}\hat{g}_{L}^{(i)}\left(  x,\operatorname*{Shell}%
\nolimits_{L}\left(  a,2a\right)  \right)  \leq C.
\]
\item[c)] For all $x,y\in V_{L}\setminus{\operatorname*{Shell}}_{L}(s(L))$,
and $i=1,2$,
\[
\hat{g}_{L}^{(i)}\left(  x,y\right)  \leq C\left\{
\begin{array}
[c]{ll}%
\frac{1}{s(L)^{2}[|x-y|\vee s(L)]^{d-2}}, & y\neq x\\
1, & y=x\,.
\end{array}
\right.
\]
\item[d)] For $i=1,2$,
\[
\sup_{x\in V_{L}}\hat{g}_{L}^{(i)}\left(  x,V_{L}\right)  \leq C\left(  \log
L\right)  ^{6}.
\]
\item[e)] For $i=1,2$,
\[
\sup_{x,x^{\prime}\in V_{L}:\left\vert x-x^{\prime}\right\vert \leq s\left(
L\right)  }\sum_{y\in V_{L}}\left\vert \hat{g}_{L}^{(i)}\left(  x,y\right)
-\hat{g}_{L}^{(i)}\left(  x^{\prime},y\right)  \right\vert \leq C\left(  \log
L\right)  ^{3}%
\]
\end{enumerate}
\end{lemma}
The proof is presented in Appendix \ref{App_B}.

Lemma \ref{Cor_Green} plays a crucial role in our smoothing procedure. 
As a
preparation, for $k\geq1$, set
\begin{equation}
\label{eq-080206}
B_{1}\left(  k\right)  \overset{\mathrm{def}}{=}\operatorname*{Shell}%
\nolimits_{L}\left(  \left(  4/3\right)  ^{k}r\left(  L\right)  \right)  .
\end{equation}
$B_{1}\left(  k\right)  \subset\operatorname*{Shell}\nolimits_{L}\left(
s\left(  L\right)  \right)  $ if $k\leq20\log\log L$. 
By Lemma \ref{Cor_Green}, 
there exists a constant $\bar c\geq1$ (again, independent of $k_0$!)
such that%
\begin{equation}
\sup_{x\in V_{L}}\hat{g}_{L}^{(1)} \left(  x,B_{1}\left(  k\right)  \right)
\leq\bar c \left\{
\begin{array}
[c]{cc}%
k\,, & \mathrm{if\ }k\leq20\log\log L\\
\left(  \log L\right)  ^{6} & \mathrm{if\ }k>20\log\log L
\end{array}
\right.  . \label{Est_BoundaryReach}%
\end{equation}
and, for any ball $V_{rs(L)}(z)\subset V_{L-s(L)}$,
$r\geq1$,
\begin{equation}
\sup_{x\in V_{L}}\hat{g}_{L}^{(1)} \left(  x,V_{rs(L)}(z) \right)  \leq\bar c
r^{d}\,. \label{Est_BoundaryReach1}%
\end{equation}
With $\bar{c}$ as in (\ref{Est_BoundaryReach}) and (\ref{Est_BoundaryReach1}),
we fix the constant $k_{0}$ large enough such that:
\begin{align}
k_{0}  &  \geq\bar{k}_{0}(1/200\bar{c}),\nonumber\\
\sup_{x\in 
\operatorname*{Sh} \nolimits_{L}
}P_{x}^{\mathrm{RW}}\left(  \tau_{V_{L}}<\tau_{V_{k_{0}%
r\left(  L\right)  }\left(  x\right)  }\right)   &  \geq9/10,\label{k_0_large}%
\\
\sup_{x\in 
\operatorname*{Sh} \nolimits_{L}
}\pi_{V_{k_{0}r\left(  L\right)  }\left(  x\right)  }\left(
x,V_{L}\right)   &  \leq17/32.\nonumber
\end{align}
That the two last estimates in (\ref{k_0_large}) hold for $k_{0}$ large
is obvious, for example from Donsker's invariance principle.

\section{Smoothed exits\label{Sect_Smooth}}
In this section, we provide estimates
on the quantity $D_{L,\Psi}(0)$. 
We use the perturbation expansion in (\ref{Pert1})
repeatedly. The main application is in comparing 
exit distributions, as follows.
If $V\subset\subset\mathbb{Z}^{d}$, and $\mathcal{S}$ is any coarse graining
scheme on $V$ (as in Definition \ref{Def_CoarseGrainingScheme}), we compare
the exit distribution of the RWRE $\Pi_{V}$ with the exit distribution
$\pi_{V}$ of simple random walk through this perturbation expansion, using
however coarse grained transitions inside $V:$ using (\ref{Green&Exit}) and
(\ref{EqualExits}) we get for $x\in V$%
\[
\left(  \Pi_{V}-\pi_{V}\right)  \left(  x,\cdot\right)  =\sum_{k=0}^{\infty
}\left(  \hat{g}_{\mathcal{S},V}\left[  \Delta_{\mathcal{S},V}\hat
{g}_{\mathcal{S},V}\right]  ^{k}\Delta_{\mathcal{S},V}\pi_{V}\right)  \left(
x,\cdot\right)  ,
\]
where%
\[
\Delta_{\mathcal{S},V}\overset{\mathrm{def}}{=}1_{V}\left(  \hat{\Pi
}_{\mathcal{S},V}-\hat{\pi}_{\mathcal{S},V}\right)  ,\ \hat{g}_{\mathcal{S}%
,V}\overset{\mathrm{def}}{=}g_{V}\left(  \hat{\pi}_{\mathcal{S},V}\right)  .
\]
Throughout this section, we consider only 
the coarse graining scheme $\mathcal{S=S} _{1}$ 
as in Definition \ref{Def_SmoothingScheme}. We
keep $L$ and $V_L$ fixed, and drop throughout the
$\mathcal{S},V$ subscripts, writing $\hat \Pi, \hat \pi,\Delta$ and $\hat g$ for
$\hat{\Pi}_{\mathcal{S},V},
\hat{\pi}_{\mathcal{S},V},
\Delta_{\mathcal{S},V}$ and $\hat{g}_{\mathcal{S},V}$. 
We use repeatedly the identity
\[
\hat{g}\left(  x,\cdot\right)  =\delta_{x,\cdot}+\hat{\pi}\hat{g}\left(
x,\cdot\right)  ,\ x\in V_L.
\]
Setting, for 
 $k\geq1$,
 \begin{equation}
	 \label{eq-120306e}
\zeta^{\left(  k\right)  }=\Delta^{k-1}\left(  \Delta\hat{\pi}\hat{g}\right)  ,
\end{equation}
we get%
\begin{equation}\label{eq-120306d}
\Pi_L-\pi_L=\hat{g}\sum_{m=1}^{\infty}\sum_{k_{1},\ldots,k_{m}=1}^{\infty
}\zeta^{\left(  k_{1}\right)  }\cdot\ldots\cdot\zeta^{\left(  k_{m-1}\right)
}\Delta^{k_{m}}\pi_L
\overset{\mathrm{def}}{=}
{\mathcal{R}}_L\,.
\end{equation}
Remark that we can replace in $\zeta^{\left(  k\right)  }$ the second part:%
\[
\left(  \Delta\hat{\pi}\hat{g}\right)  \left(  x,y\right)  =\sum_{z}\left(
\Delta\hat{\pi}\right)  \left(  x,z\right)  \left(  \hat{g}\left(  z,y\right)
-\hat{g}\left(  x,y\right)  \right)  ,
\]
i.e., we gain a discrete derivative in the Green function.

We can now describe informally our basic strategy. 
When analyzing the term
$D_{L,\Psi}(0)$, boundary effect are not essential, and one can
consider all steps to be coarse-grained  (some extra care is however
needed near the boundary, which leads to the specific form
of the coarse graining scheme ${\mathcal {S}}_1$,  
but we gloss over these details in the description 
that follows).
Note that the steps of the coarse-grained random walk are essentially
in the scale $L/(\log L)^3$. In this scale, most 
$x\in V_L$ are good, that is the individual steps of the coarse-grained
random walk are 
controlled by the good event in the induction hypothesis. 
Consider the linear
term in (\ref{eq-120306d}), that is the term with $m=1$, which turns out
to be the dominant term in the expansion.
Suppose first all $x\in V_L$ are good, and consider
the term with $k_1=0$. In this case, each term is smoothed at scale
$L$ from the right, and its variational norm is bounded
by $ o((\log L)^{-3}) O((\log L)^{-9})$. A-priori estimates
on the coarse-grained simple random walk yield that the
sum over the coarse grained Green function $\hat g$ is 
$O((\log L)^{-6})$. This would look alarming, as multiplying these
gives rise to an error which is only
 $o((\log L)^{-6})$, which
could
result in non-propagation of the induction hypothesis. However,
one can use the fact that the individual contributions from sites
distance by $\rho_{1,L}$ are independent, and of zero mean due
to the isotropy assumption. Averaging over this sum of essentially
independent random variables improves the estimate from the worst-case
value of  $o((\log L)^{-6})$ back to
the desired value of  $o((\log L)^{-9})$, see 
the proof of Proposition \ref{prop-erwin160605}.
The terms with $k_1\geq 1$ are handled similarly, using now
the part of the induction hypothesis involving $D_{L,0}(0)$
to control the extra powers of $\Delta$ and ensure the convergence 
of the series. A similar strategy is applied to the ``non-linear''
terms with $m>1$. Boundary terms are handled by using the fact that 
the coarse grained random walk is unlikely to stay at distance
less than $r(L)$ from the boundary for many steps.

A major complication in handling the perturbation expansion is
the presence of
``bad regions''. 
The advantage of 
the coarse graining scheme $\mathcal{S=S}_{1}$  is that it is unlikely
to have  more than one ``bad region'', and that this
single bad region can be handled by an appropriate surgery,
once appropriate 
Green function estimates for the RWRE in a ``good environment''
are derived, see Section \ref{Subsect_greengood}.

We now turn to the actual proof, and
write $B_{L}^{\left(  i\right)  }$, $i=1,2,3,4,$ 
for the collection of
points which are bad on level $i,$ and in the right scale,
with respect to the coarse graining scheme $\mathcal{S}_{1}$. That is,
for $i=1,2,3$,
\begin{align}
	\label{eq-120306b}
B_{L}^{\left(  i\right)  }=&
\{x\notin
\operatorname*{Sh}\nolimits_{L}: 
D_{r,h_{L}(x)}\left(  x\right)  >\left(  \log L\right)  ^{-9+\frac{9(i-1)}{4}}
\,\mbox{\rm for some }\ 
r\in\lbrack h_{L}(x),2h_{L}(x)],
\nonumber\\
&
D_{r,h_{L}(x)}\left(
x\right)  \leq\left(  \log L\right)  ^{-9+\frac{9i}{4}}
\mbox{\rm for all}\ 
r\in\lbrack h_{L}(x),2h_{L}(x)]\,,\
D_{r,0}\left(
x\right)  \leq\delta
\}\,,
\end{align}
and
\begin{align}
	\label{eq-120306c}
B_{L}^{\left(  4\right)  }=&
\{x\notin
\operatorname*{Sh}\nolimits_{L}: 
D_{r,h_{L}(x)}\left(  x\right)  >\left(  \log L\right)  ^{-\frac{9}{4}}
\ \mbox{\rm or}\
D_{r,0}\left(  x\right)  >\delta
\,,\nonumber \\
& \mbox{\rm for some }\ 
r\in\lbrack h_{L}(x),2h_{L}(x)]
%&\mbox{\rm or}\
%D_{r,0}\left(  x\right)  >\delta
\}\; \bigcap
\{x \in
\operatorname*{Sh}\nolimits_{L}: 
D_{k_0r\left(  L\right) ,0 }\left(  x\right)
\geq\delta\}\,.
\end{align}
%
%For $i=1,2,3$,
%$B_{L}^{\left(  i\right)  },$ are the set of points $x\notin
%\operatorname*{Shell}\nolimits_{L}\left(  r\left(  L\right)  \right)  $ such
%that for \textit{some} $r\in\lbrack h_{L}(x),2h_{L}(x)]$, one has
%$D_{r,h_{L}(x)}\left(  x\right)  >\left(  \log L\right)  ^{-11.25+2.25i},\ $
%but for all $r\in\lbrack h_{L}(x),2h_{L}(x)],$ $D_{r,h_{L}(x)}\left(
%x\right)  \leq\left(  \log L\right)  ^{-9+2.25 i},$ and $D_{r}^{0}\left(
%x\right)  \leq\delta.$ $B_{L}^{\left(  4\right)  }$ is the set of points $x$
%which for $d_{L}\left(  x\right)  >r\left(  L\right)  $ have the property that
%for some $r$, $h_{L}\left(  x\right)  \leq r\leq2h_{L}\left(  x\right)  ,$
%$D_{r,h_{L}(x)}\left(  x\right)  >\left(  \log L\right)  ^{-2.25},$ or
%$D_{r}^{0}\left(  x\right)  >\delta,$ and for $d_{L}\left(  x\right)  \leq
%r\left(  L\right)  $ satisfy 
%$D_{k_0r\left(  L\right)  }^{0}\left(  x\right)
%\geq\delta$. 
We also write
\begin{equation}
	B_{L}\overset{\mathrm{def}}{=}\bigcup_{i=1}^4 B_{L}^{\left(  i\right)}
	\,,\
\operatorname*{Good}\nolimits_{L}\overset{\mathrm{def}}{=}\left\{
B_{L}=\emptyset\right\}  . 
%\label{Def_GoodL}%
\label{Def_BL}%
\end{equation}

As mentioned in the beginning of this section,
a major complication in handling the perturbation expansion is
the ``bad regions''. 
The advantage of 
the coarse graining scheme $\mathcal{S=S}_{1}$  is that it is unlikely
to have essentially more than one ``bad region''.
To make this statement precise, note that
if $L_{1}\leq L\leq L_{1}\left(  \log L_{1}\right)  ^{2}$ then all the radii
involved in the definition of badness are smaller than $L_{1},$ if $L_{1}$ is
chosen large enough. Remark also that if $d_{L}\left(  x\right)  >r\left(
L\right)  ,$ then $h_{L}\left(  x+\cdot\right)  \in\mathcal{M}_{r}$ for
$h_{L}\left(  x\right)  \leq r\leq2h_{L}\left(  x\right)  ,$ 
and therefore, if
$L_1$ is large enough,
$\operatorname*{Cond}\left(  \delta,L_{1}\right)  $ holds, 
and $L_{1}\leq
L\leq L_{1}\left(  \log L_{1}\right)  ^{2},$ then%
\begin{equation}
\mathbb{P}\left(  x\in B_{L}\right)  \leq2\gamma s\left(  L\right)
\exp\left[  -\frac{10}{13}\left(  \log\frac{\gamma L}{\left(  \log L\right)
^{10}}\right)  ^{2}\right]  \leq\exp\left[  -0.7\left(  \log L\right)
^{2}\right]  \,. \label{BoundBadness}%
\end{equation}

The points $y$ whose random environment $\omega_{y}$ can influence the badness
of $x$ are evidently within radius $\rho_L(x)=\rho_{1,L}\left(  x\right)  $ 
from $x$, see (\ref{Def_Rho1}). If
$\left\vert x-y\right\vert >\rho_{L}\left(  x\right)  +\rho_{L}\left(
y\right)  ,$ then $\left\{  x\in B_{L}\right\}  $ and $\left\{  y\in
B_{L}\right\}  $ are independent. Therefore, if we define%
\begin{equation}
\operatorname*{TwoBad}\nolimits_{L}\overset{\mathrm{def}}{=}\bigcup_{x,y\in
V_{L}:\left\vert x-y\right\vert >\rho_{L}\left(  x\right)  +\rho
_{L}\left(  y\right)  }\left\{  x\in B_{L}\right\}  \cap\left\{  y\in
B_{L}\right\}  , \label{DefTwoBad}%
\end{equation}
then:
\begin{lemma}
\label{Le_TwoBad}Assume $L_{1}$ large enough, (\ref{BoundBad}) for $L_{1},$
and $L_{1}\leq L\leq L_{1}\left(  \log L_{1}\right)  ^{2}.$ Then%
\[
\mathbb{P}\left(  \operatorname*{TwoBad}\nolimits_{L}\right)  \leq\exp\left[
-1.2\left(  \log L\right)  ^{2}\right]  .
\]
\end{lemma}

Next, we regard $\hat{\Pi}$ as a field $\left(
\hat{\Pi}\left(  x,\cdot\right)  \right)  _{x\in V_{L}}$ of
random transition probabilities. We defined the \textquotedblleft
goodified\textquotedblright\ transition probabilities%
\begin{equation}
	\label{eq-160306a}
\operatorname*{gd}\left(  \hat{\Pi}\right)  \left(
x,\cdot\right)  \overset{\mathrm{def}}{=}\left\{
\begin{array}
[c]{cc}%
\hat{\Pi}\left(  x,\cdot\right)  & \mathrm{if\ }x\notin
B_{L}\\
\hat{\pi}\left(  x,\cdot\right)  & \mathrm{if\ }x\in B_{L}%
\end{array}
\right.  .
\end{equation}
This field might no longer come from an i.i.d. RWRE, but nevertheless, we have
the property that $\operatorname*{gd}\left(  \hat{\Pi}_L\right)
\left(  x,\cdot\right)  $ and $\operatorname*{gd}\left(  \hat{\Pi
}_L\right)  \left(  y,\cdot\right)  $ are independent provided
$\left\vert x-y\right\vert >\rho_{L}\left(  x\right)  +\rho_{L}\left(
y\right)  .$ If $X$ is a random variable depending on $\omega$ only trough 
$\hat{\Pi}_L$ we define $\operatorname*{gd}\left(  X\right)  $
by replacing $\hat{\Pi}_L$ by $\operatorname*{gd}\left(
\hat{\Pi}_L\right)  .$

We next take ${\Psi}\in\mathcal{M}_{L},$ and
set $\phi\overset{\mathrm{def}}{=}\phi_{L,{\Psi}},$ as in (\ref{Def_PhiLM}).
An easy consequence of our definitions and 
Lemma \ref{lem-130705} is the following.
\begin{lemma}
\label{lem-090805}
%There is a $\delta_0>0$ and a constant $C_g\geq 1$
If $\delta\leq(1/800\bar c)$ then, for all $x\in V_{L}$ and $k\geq2$,
\begin{equation}
\label{eq-090805a}\mathbf{1}_{\{B_{L}=\emptyset\}} \|\Delta^{k}(x,\cdot)\|_{1}
\leq\frac{1}{\bar c} \left(  \frac18\right)  ^{k} \,.
\end{equation}
\end{lemma}
\begin{proof}
Since $\max_{x\in V_{L}}\left\Vert \Delta(x,\cdot)\right\Vert _{1}\leq2$ and
$\bar{c}\geq1$, it is enough to prove that
\[
\mathbf{1}_{\{B_{L}=\emptyset\}}\sum_{z\in V_{L}}|\Delta^{2}(x,z)|\leq\left(
\frac{1}{64\bar{c}}\right)  \,.
\]
If $x\not\in\operatorname*{Sh}_L$ 
then, on the event $\left\{  B_{L}=\emptyset\right\}  $,
$\Vert\Delta(x,\cdot)\Vert_{1}\leq\delta$ and hence $\Vert\Delta^{2}%
(x,\cdot)\Vert_{1}\leq2\delta\leq1/64\bar{c}$ due to our choice of $\delta$.
On the other hand, if 
$x\in\operatorname*{Sh}_L$ 
then on the event $\left\{
B_{L}=\emptyset\right\}  $,
\begin{align}
&  \sum_{z\in V_{L}}|\Delta^{2}(x,z)|=\sum_{z\in V_{L}}\left\vert
\sum\nolimits_{y\in V_{L}}\Delta(x,y)\Delta(y,z)\right\vert \label{eq-090805b}%
\\
&  \leq2\left\vert \sum_{y\in{\operatorname*{Sh}}_{L}}\Delta
(x,y)\right\vert +\left\vert \sum_{y\in V_{L}\setminus{\operatorname*{Sh}%
}_{L}}\Delta(x,y)\right\vert \max_{y\in V_{L}\setminus
{\operatorname*{Sh}}_{L}}\sum_{z\in V_{L}}|\Delta(y,z)|\nonumber\\
&  \leq2(\delta+\frac{1}{200\bar{c}})+2\delta=4\delta+\frac{1}{100\bar{c}%
}<\frac{1}{64\bar{c}}\,,\nonumber
\end{align}
where Lemma \ref{lem-130705} and $k_{0}\geq\bar{k}_{0}(1/200\bar{c})$ were
used in the next to last inequality.
\end{proof}

In what follows, we will always consider $\delta\leq1/800\bar c$.
\subsection{The linear part \label{SubSect_GoodLinear}}
For $x\in V_{L},\ B\subset V_{L},$ set%
\begin{align}
\xi_{x}^{\left(  k\right)  }\left(  B,z\right)   &  =\sum_{y\in B}\hat
{g}\left(  x,y\right)  \left(  \Delta^{k}\pi_L\hat{\pi}_{{\Psi}}\right)  \left(
y,z\right) \label{XiB}\\
&  =\sum_{y\in B}\sum_{y^{\prime}\in V_{L}}\hat{g}\left(  x,y\right)
\Delta^{k}\left(  y,y^{\prime}\right)  \left(  \phi\left(  y^{\prime
},z\right)  -\phi\left(  y,z\right)  \right)\,, \nonumber
\end{align}
where the last equality is because the total mass of $\Delta(y,\cdot)$
vanishes. 

We write $\xi_{x}^{\left(  k\right)  }\left(  z\right)  $ for $\xi
_{x}^{\left(  k\right)  }\left(  V_{L},z\right)$;
in the notation of (\ref{eq-120306d}),
$\xi_x^{\left(k\right)}\left( z\right)=\zeta^{\left(k\right)}\hat\pi_{\Psi}
(x,z)$.
Define%
\[
G_{L}\overset{\mathrm{def}}{=}\left\{  \sup_{x\in V_{L}}\sum\nolimits_{k\geq
1}\left\Vert \xi_{x}^{\left(  k\right)  }\right\Vert _{1}\leq\left(  \log
L\right)  ^{-37/4}\right\}  .
\]
$G_L $ is precisely the event that the $m=1$
term in the perturbation
expansion (\ref{eq-120306d}), smoothed by $\hat \pi_\Psi$, is ``small''.
\begin{proposition}
\label{prop-erwin160605} If $L$ is large enough, then%
\[
\mathbb{P}\left(  \left(  G_{L}\right)  ^{c}\cap\operatorname*{Good}%
\nolimits_{L}\right)  \leq\exp\left[  -\left(  \log L\right)  ^{17/8}\right]
.
\]
\end{proposition}
\begin{proof}
It suffices to prove that%
\[
\sup_{x\in V_{L}}\mathbb{P}\left(  \sum\nolimits_{k\geq1}\left\Vert \xi
_{x}^{\left(  k\right)  }\right\Vert _{1}\geq\left(  \log L\right)
^{-37/4},\ \operatorname*{Good}\nolimits_{L}\right)  \leq\exp\left[  -\left(
\log L\right)  ^{9/4}\right]\,.
\]%
Note that
\begin{align*}
&  \mathbb{P}\left(  \sum\nolimits_{k\geq1}\left\Vert \xi_{x}^{\left(
k\right)  }\right\Vert _{1}\geq\left(  \log L\right)  ^{-37/4}%
,\ \operatorname*{Good}\nolimits_{L}\right) \\
&  
=\mathbb{P}\left(  \sum\nolimits_{k\geq1}\left\Vert \operatorname*{gd}%
\left(  \xi_{x}^{\left(  k\right)  }\right)  \right\Vert _{1}\geq\left(  \log
L\right)  ^{-37/4},\ \operatorname*{Good}\nolimits_{L}\right) \\
&  \leq\mathbb{P}\left(  \sum\nolimits_{k\geq1}\left\Vert \operatorname*{gd}%
\left(  \xi_{x}^{\left(  k\right)  }\right)  \right\Vert _{1}\geq\left(  \log
L\right)  ^{-37/4}\right)  .
\end{align*}
For notation convenience, we drop the notation $\operatorname*{gd}\left(
\cdot\right)  ,$ and just use the fact that all $\hat{\Pi}$
involved satisfy the appropriate \textquotedblleft goodness\textquotedblright%
\ properties. (Remark that after \textquotedblleft
goodifications\textquotedblright, the distribution of $\hat{\Pi}
\left(  x,x+\cdot\right)  $ remains invariant under lattice
isometries, provided $d_{L}\left(  x\right)  >2s\left(  L\right)  .$)

We split $\xi_{x}^{\left(  k\right)  }$ into different parts. If 
$y\not\in \operatorname*{Sh}_L$
and $\Delta\left(  y,y^{\prime}\right)  >0,$
we have, since $\gamma\leq1/8$, that $\left\vert y-y^{\prime}\right\vert \leq
d_{L}\left(  y\right)  /4,$ i.e. $d_{L}\left(  y\right)  \leq\left(
4/3\right)  d_{L}\left(  y^{\prime}\right)  .$ Therefore, if $
y^{\prime}\in \operatorname*{Sh}_L  $ and $\Delta^{k}\left(
y,y^{\prime}\right)  >0,$ then $d_{L}\left(  y\right)  \leq\left(  4/3\right)
^{k}r\left(  L\right)  .$ Recall the set $B_{1}(k)$, c.f.
(\ref{eq-080206}),
 and the estimate
(\ref{Est_BoundaryReach}).
If $y\in B_{1}(k),$ and $\Delta^{k}\left(  y,y^{\prime}\right)  >0,$ we have%
$$\left\vert y-y^{\prime}\right\vert    \leq kk_{0}r\left(  L\right)
+3^{k}\max\left(  r\left(  L\right)  ,d_{L}\left(  y\right)  \right) 
  \leq\left(  kk_{0}+4^{k}\right)  r\left(  L\right)  ,
$$
and applying (\ref{OneDerivative}), we see that for $y\in B_{1}(k),$ and
$y^{\prime}$ such that 
$\Delta^{k}\left(  y,y^{\prime}\right)  >0,$ we have%
\[
\left\vert \phi\left(  y,z\right)  -\phi\left(  y^{\prime},z\right)
\right\vert \leq C\left(  kk_{0}+4^{k}\right)  L^{-d}\left(  \log L\right)
^{-10}.
\]
By Lemma \ref{lem-090805},
we have $\left\Vert \Delta^{k}\left(  y,\cdot\right)  \right\Vert _{1}\leq
8^{-k}.$ Combining all these estimates with parts b) and d) of 
Lemma \ref{Cor_Green}, we have%
\begin{equation}
\left\Vert \xi_{x}^{(k)}\left(  B_{1}(k)\right)  \right\Vert _{1}\leq
C\left\{
\begin{array}
[c]{ll}%
8^{-k}\left(  kk_{0}+4^{k}\right)  \left(  \log L\right)  ^{-10}, &
\mathrm{if\ }k\leq20\log\log L,\\
8^{-k}\left(  kk_{0}+4^{k}\right)  \left(  \log L\right)  ^{-4}\,, &
\mathrm{if\ }k>20\log\log L.
\end{array}
\right.  \label{Est_B_1(k)}%
\end{equation}
(We emphasize our convention regarding constants, and in particular the fact
that $C$ does not depend on $x$.) Hence,
\begin{equation}
\sup_{x}\sum_{k\geq1}\left\Vert \xi_{x}^{\left(  k\right)  }\left(
B_{1}(k)\right)  \right\Vert _{1}\leq C\left(  \log L\right)  ^{-10}%
\leq\left(  \log L\right)  ^{-37/4}/3. \label{LinearEst1}%
\end{equation}

Next, let%
\[
B_{2}(k)\overset{\mathrm{def}}{=}\operatorname*{Shell}\nolimits_{L}\left(
\left(  4/3\right)  ^{k}r\left(  L\right)  ,\left(  5/4\right)  ^{k}2s\left(
L\right)  \right)  .
\]
If $y\in B_{2}(k)$ and $\Delta^{k}\left(  y,y^{\prime}\right)  >0,$ we have
$d_{L}\left(  y^{\prime}\right)  >r\left(  L\right)  ,$ and we get, using the
fact that for $x\not\in \operatorname*{Sh}_L$
one can write
$\pi_L\left(  x,\cdot\right)  =\left(  \hat{\pi}\pi_L\right)  \left(
x,\cdot\right)  ,$%
\[
\xi_{x}^{\left(  k\right)  }\left(  B_{2}(k),z\right)  =\sum_{y\in B_{2}%
(k)}\sum_{y^{\prime}\in V_{L}}\hat{g}\left(  x,y\right)  D_{k}\left(
y,y^{\prime}\right)  \left(  \phi\left(  y,z\right)  -\phi\left(  y^{\prime
},z\right)  \right)  ,
\]
where%
\begin{equation}
D_{k}\overset{\mathrm{def}}{=}\Delta^{k}\hat{\pi}, \label{Dk}%
\end{equation}
and, on a ``good'' environment,
\begin{align}
\sup_{y\in B_{2}(k)}\left\Vert D_{k}\left(  y,\cdot\right)  \right\Vert _{1}
&  \leq\sup_{y\in B_{2}(k)}\left\Vert \Delta^{k-1}\left(  y,\cdot\right)
\right\Vert _{1}\sup_{x:d_{L}\left(  x\right)  >r\left(  L\right)  }\left\Vert
\Delta\hat{\pi}\left(  x,\cdot\right)  \right\Vert _{1}\label{ESt_Dk}\\
&  \leq C8^{-k}\left(  \log L\right)  ^{-9}.\nonumber
\end{align}
Using Lemma \ref{Cor_Green} b), we have $\sup_{x}\hat{g}\left(
x,\operatorname*{Shell}\nolimits_{L}\left(  3s\left(  L\right)  \right)
\right)  \leq C\log\log L.$ Put%
\[
A_{j}\overset{\mathrm{def}}{=}\operatorname*{Shell}\nolimits_{L}\left(
\left(  2+\left(  j-1\right)  /4\right)  s\left(  L\right)  ,\left(
2+j/4\right)  s\left(  L\right)  \right)  ,j\geq1.
\]
Starting from a point in $A_{j},$ $j\geq3,$ the coarse grained 
simple random walk has
a probability $\geq1/C$ to reach $A_{j-2}$ in one step. Starting from
$A_{j-2},$ an ordinary random walk has a probability $\geq1/C$ to leave
$V_{L+k_{0}r\left(  L\right)  }$ before reaching $A_{j},$ and therefore, the
coarse grained simple random walk
leaves $V_{L}$ before reaching $A_{j}$ with at least the
same probability. Therefore $\sup_{x}\hat{g}\left(  x,A_{j}\right)  \leq Cj,$
and thus,%
\[
\sup_{x}\hat{g}\left(  x,B_{2}(k)\right)  \leq C\left(  \left(  \frac{5}%
{4}\right)  ^{2k}+\log\log L\right)  \leq C\left(  2^{k}+\log\log L\right)
\,.
\]
If $y\in B_{2}(k),$ and $\Delta^{k}\left(  y,y^{\prime}\right)  >0$, then
$\left\vert y-y^{\prime}\right\vert \leq2ks\left(  L\right)  ,$ and therefore,
\[
\left\vert \phi\left(  y,z\right)  -\phi\left(  y^{\prime},z\right)
\right\vert \leq CkL^{-d}\left(  \log L\right)  ^{-3},
\]
again by (\ref{OneDerivative}). Therefore, we get$,$%
\begin{align}
\left\Vert \xi_{x}^{\left(  k\right)  }\left(  B_{2}(k),\cdot\right)
\right\Vert _{1}  &  \leq Ck\left(  \log L\right)  ^{-12}\left[  4^{-k}%
+8^{-k}\log\log L\right]  ,\nonumber\\
\sup_{x}\sum_{k\geq1}\left\Vert \xi_{x}^{\left(  k\right)  }\left(
B_{2}(k),\cdot\right)  \right\Vert _{1}  &  \leq\left(  \log L\right)
^{-37/4}/3. \label{LinearEst2}%
\end{align}

Let $B_{3}(k)\overset{\mathrm{def}}{=}V_{L}\backslash\left(  B_{1}(k)\cup
B_{2}(k)\right)  .$ Given $j\in\mathbb{Z},$ let%
\[
I_{j}\overset{\mathrm{def}}{=}\left\{  jks\left(  L\right)  +1,\ldots,\left(
j+1\right)  ks\left(  L\right)  \right\}  .
\]
Then for $\mathbf{j}\in\mathbb{Z}^{d},$ put $W_{\mathbf{j},k}\overset
{\mathrm{def}}{=}B_{3}(k)\cap I_{j_{1}}\times\cdots\times I_{j_{d}},$ 
with
$\operatorname*{diameter}\left(  W_{\mathbf{j}}\right)  \leq\sqrt
{d}ks\left(  L\right).$
Let $J_k$
be the set of $\mathbf{j}$'s for which these sets are not empty. We 
subdivide
$J_k$ into 
subsets $J_{1,k},\ldots,J_{K\left(  d,k\right),k  }$ such that for any
$1\leq \ell \leq K\left(  d,k\right)  ,$%
\begin{equation}
\mathbf{j},\mathbf{j}^{\prime}\in J_{\ell,k},\ \mathbf{j}\neq\mathbf{j}^{\prime
}\Longrightarrow d\left(  W_{\mathbf{j},k},W_{\mathbf{j}^{\prime},k}\right)
>ks\left(  L\right)  . \label{Box_Interdistance}%
\end{equation}
%We also have $\operatorname*{diam}\left(  W_{\mathbf{j}}\right)  \leq\sqrt
%{d}ks\left(  L\right)  .$

We set, recalling (\ref{Dk}),%
\begin{equation}
\xi_{x,\mathbf{j}}^{\left(  k\right)  }\left(  z\right)  \overset
{\mathrm{def}}{=}\sum_{y\in W_{\mathbf{j},k}}\sum_{y^{\prime}\in V_{L}}\hat
{g}\left(  x,y\right)  D_{k}\left(  y,y^{\prime}\right)  \left(  \phi\left(
y,z\right)  -\phi\left(  y^{\prime},z\right)  \right)  . \label{Xi_Lj}%
\end{equation}
We fix for the moment $k$ and $x.$ If $t>0$, and%
\begin{equation}
\sum\nolimits_{\mathbf{j}}\mathbb{E}\xi_{x,\mathbf{j}}^{\left(  k\right)
}\left(  z\right)  \leq t/2, \label{Annealed}%
\end{equation}
and we have%
\begin{multline*}
\mathbb{P}\left(  \left\Vert \xi_{x}^{\left(  k\right)  }\left(
B_{3}(k),\cdot\right)  \right\Vert _{1}\geq t\right)  \leq\mathbb{P}\left(
\left\vert \sum\nolimits_{\mathbf{j}}\left(  \xi_{x,\mathbf{j}}^{\left(
k\right)  }\left(  z\right)  -\mathbb{E}\xi_{x,\mathbf{j}}^{\left(  k\right)
}\left(  z\right)  \right)  \right\vert \geq t/2\right) \\
\leq K\left(  d,k\right)  \max_{1\leq \ell \leq K\left(  d,k\right)
}\mathbb{P}%
\left(  \left\vert \sum\nolimits_{\mathbf{j}\in J_{\ell,k}}\left(  \xi
_{x,\mathbf{j}}^{\left(  k\right)  }\left(  z\right)  -\mathbb{E}%
\xi_{x,\mathbf{j}}^{\left(  k\right)  }\left(  z\right)  \right)  \right\vert
\geq t/\left(  2K\left(  d,k\right)  \right)  \right)  .
\end{multline*}
The random variables $\xi_{x,\mathbf{j}}^{\left(  k\right)  }\left(  z\right)
-\mathbb{E}\xi_{x,\mathbf{j}}^{\left(  k\right)  }\left(  z\right)  $,
$\mathbf{j}\in J_{\ell,k}$, are independent and centered, due to
(\ref{Box_Interdistance}), and we are going to estimate their sup-norm.
We have by (\ref{OneDerivative}) that
$|\phi\left(  y,z\right)  -\phi\left(  y^{\prime},z\right)  |\leq
C k\left(  \log L\right)  ^{-3}L^{-d}$ for
$y,y^{\prime}$ for which $D_{k}\left(  y,y^{\prime}\right)  \neq0.$ According
to Lemma \ref{Cor_Green} c), we have%
\[
\hat{g}\left(  x,W_{\mathbf{j},k}\right)  \leq Ck^{d}\left(  1+\frac{d\left(
x,W_{\mathbf{j},k}\right)  }{s\left(  L\right)  }\right)  ^{-d+2}.
\]
Substituting that into (\ref{Xi_Lj}), we get%
\[
\left\Vert \xi_{x,\mathbf{j}}^{\left(  k\right)  }\left(  z\right)
\right\Vert _{\infty}\leq Ck^{d+1}8^{-k}\left(  1+\frac{d\left(
x,W_{\mathbf{j},k}\right)  }{s\left(  L\right)  }\right)  ^{-d+2}L^{-d}\left(
\log L\right)  ^{-12}.
\]
By Hoeffding's inequality (see e.g. \cite[(1.23)]{Ledoux} ), we have for
$1\leq \ell \leq K\left(d,k\right)  $%
\begin{align*}
&  \mathbb{P}\left(  \left\vert \sum\nolimits_{\mathbf{j}\in J_{\ell,k}}\left(
\xi_{x,\mathbf{j}}^{\left(  k\right)  }\left(  z\right)  -\mathbb{E}%
\xi_{x,\mathbf{j}}^{\left(  k\right)  }\left(  z\right)  \right)  \right\vert
\geq\frac{2^{-k}L^{-d}}{2K\left(d,k\right)  \left(  \log L\right)  ^{37/4}%
}\right) \\
&  \leq2\exp\left[  -\frac{1}{C}\frac{\left(  \log L\right)  ^{-37/2}%
}{k^{2d+2}4^{-2k}\left(  \log L\right)  ^{-24}\sum\nolimits_{r=1}^{C\left(
\log L\right)  ^{3}}r^{-d+3}}\right]  \leq2\exp\left[  -\frac{1}{C}%
\frac{\left(  \log L\right)  ^{5/2}}{k^{2d+2}4^{-2k}}\right]  \,,
\end{align*}
where we used $d\geq3$ in the last inequality. The upshot of this estimate is
that provided (\ref{Annealed}) holds true with $t=t(k)=
2^{-k}L^{-d}\left(  \log
L\right)  ^{-37/4},$ we have%
\begin{align*}
\sup_{x}\mathbb{P}\left(  \sum\nolimits_{k\geq1}\left\Vert \xi_{x}^{\left(
k\right)  }\left(  B_{3}\right)  \right\Vert _{1}\geq\left(  \log L\right)
^{-37/4}\right)   &  \leq2\sum_{k\geq1}K(d,k)
\exp\left[  -\frac{1}{C}\frac{\left(
\log L\right)  ^{5/2}}{k^{2d+2}4^{-2k}}\right] \\
&  \leq\exp\left[  -\left(  \log L\right)  ^{17/8}\right]  ,
\end{align*}

It remains to prove (\ref{Annealed}) with this $t$. Write
\[
\sum_{\mathbf{j}}\mathbb{E}\xi_{x,\mathbf{j}}^{\left(  k\right)  }\left(
z\right)  =\sum_{y\in B_{3}}\sum_{y^{\prime}\in V_{L}}\hat{g}\left(
x,y\right)  \mathbb{E}\left(  D_{k}\left(  y,y^{\prime}\right)  \right)
\left(  \phi\left(  y,z\right)  -\phi\left(  y^{\prime},z\right)  \right)  .
\]
For every $y,$ $y^{\prime}\mapsto\mathbb{E}\left(  D_{k}\left(  y,y^{\prime
}\right)  \right)  $ is a signed measure with total mass $0,$ which is
invariant under lattice isometries. Furthermore%
\[
\sum_{y^{\prime}}\left\vert \mathbb{E}\left(  D_{k}\left(  y,y^{\prime
}\right)  \right)  \right\vert \leq C8^{-k}\left(  \log L\right)  ^{-9}.
\]
Applying (\ref{ThreeDerivatives}), we get%
\begin{align*}
&  \left\vert \sum\nolimits_{y^{\prime}}\mathbb{E}\left(  D_{k}\left(
y,y^{\prime}\right)  \right)  \left(  \phi\left(  y,z\right)  -\phi\left(
y^{\prime},z\right)  \right)  \right\vert \\
&  \leq C8^{-k}\left(  \log L\right)  ^{-9}\left(  L^{-d-1/4}+\left(
\frac{Lk\left(  \log L\right)  ^{-3}}{L}\right)  ^{3}L^{-d}\right)  \leq
C4^{-k}\left(  \log L\right)  ^{-18}L^{-d},
\end{align*}
uniformly in $y\in B_{3}(k),$ and $k.$ By Lemma \ref{Cor_Green} d), we have%
\[
\sup_{x}\sum_{y\in B_{3}(k)}\hat{g}\left(  x,y\right)  \leq C\left(  \log
L\right)  ^{6}.
\]
From this (\ref{Annealed}) follows.
\end{proof}

\subsection{The non-linear part, good environment \label{SubSect_Nonlinear}}
%Recall the random variable 
%$ D_{L,\Psi}\left(  0\right)$, c.f.
%(\ref{Def_DL}).
\begin{proposition}
\label{Prop_Main_no_bad}If $L$ is large enough and $\Psi\in\mathcal{M}_{L},$
then, with
$ D_{L,\Psi}\left(  0\right)$ as in
(\ref{Def_DL}),
\[
\mathbb{P}\left(  D_{L,\Psi}\left(  0\right)  \geq\left(  \log L\right)
^{-9};\operatorname*{Good}\nolimits_{L}\right)  \leq\exp\left[  -\left(  \log
L\right)  ^{17/8}\right]  .
\]
\end{proposition}
\begin{proof}
	We recall the abbreviation $\operatorname*{Sh}_{L}
	\overset{\mathrm{def}}{=}%
\operatorname*{Shell}\nolimits_{L}\left(  r\left(  L\right)  \right),  $ 
c.f. (\ref{Def_sL&rL}).
By Proposition \ref{prop-erwin160605}, it suffices to estimate 
on $G_{L}\cap\operatorname*{Good}\nolimits_{L}$ the expression
$\|{\mathcal{R}}_L \hat \pi_{\Psi}\|_1$, c.f. (\ref{eq-120306d}), where
\begin{equation}
	{\mathcal{R}}_{L}\hat\pi_{\Psi}
\overset{\mathrm{def}}{=}\sum_{m=1}^{\infty}\sum_{k_{1},\ldots,k_{m}%
=0}^{\infty}\left(  \hat{g}\Delta^{k_{1}}\Delta\hat{\pi}
1_{V_{L}}
\right)  \cdot
\ldots\cdot\left(  \hat{g}\Delta^{k_{n}}\Delta\hat{\pi}1_{V_{L}}\right)  \sum
_{k=1}^{\infty}\left(  \hat{g}\Delta^{k}\phi\right)  , \label{NonLinearSplit}%
\end{equation}
%on $G_{L}\cap\operatorname*{Good}\nolimits_{L}.$ 
and $\phi=\pi_L\hat \pi_{\Psi}$. The last factor in the
right hand side of (\ref{NonLinearSplit}) is
$\sum_{k=1}^{\infty}\xi^{\left(  k\right)  }$ of the last section, and
therefore, it suffices to 
show that on $\operatorname*{Good}\nolimits_{L},$
%the other factors are staying below $1,$ for instance%
\begin{equation}
\sup_{x}\sum_{k\geq0}\left\Vert \left(  \hat{g}\Delta^{k}\Delta\hat{\pi
}1_{V_{L}}
\right)  \left(  x,\cdot\right)  \right\Vert _{1}\leq15/16. \label{Est_15/16}%
\end{equation}
%Recall that by \cite[Theorem 1.7.1]{lawler},\ref{check!!},
%$\|\hat \pi(x,\cdot)-\hat\pi(x',\cdot)\|_1\leq |x-x'|/L$. 
Using  the definition of 
$\operatorname*{Good}\nolimits_{L}$ 
in the first inequality and Lemma \ref{Cor_Green} d) 
together with Lemma \ref{lem-090805} in the second, we get
	$$\sup_{y\notin \operatorname*{Sh}_{L}}
	\left\Vert \left(  \Delta\hat{\pi}\right)  \left(
y,\cdot\right)  \right\Vert _{1}  \leq C\left(  \log L\right)  ^{-9}\,,\
\sum_{k\geq0}\sup_{x}\left\Vert \left(  \hat{g}\Delta^{k}\right)  \left(
x,\cdot\right)  \right\Vert _{1}    \leq C\left(  \log L\right)  ^{6}.
$$
Therefore, we have%
\[
\sum_{k\geq0}\sup_{x}\left\Vert \sum\nolimits_{y\notin \operatorname*{Sh}_{L}}
\left(  \hat
{g}\Delta^{k}\right)  \left(  x,y\right)  \left(  \Delta\hat{\pi}\right)
\left(  y,\cdot\right)  \right\Vert _{1}\leq1/16,
\]
if $L$ is large enough, and in order to prove (\ref{Est_15/16}) it therefore
suffices to prove%
\[
\sum_{k\geq0}\sup_{x}\left\Vert \sum\nolimits_{y\in \operatorname*{Sh}_{L}}
\left(  \hat
{g}\Delta^{k}\right)  \left(  x,y\right)  \left(  \Delta\hat{\pi}1_{V_{L}}
\right)
\left(  y,\cdot\right)  \right\Vert _{1}\leq7/8.
\]
As in the proof of proposition (\ref{prop-erwin160605}), if $\Delta
^{k}(z,y)>0$ for $y\in \operatorname*{Sh}_{L}$ 
then $z\in B_{1}(k)$. Hence, using
(\ref{Est_BoundaryReach}) and Lemma \ref{lem-090805},
%together with our choice
%of $k_{0}$ in the second inequality,
\begin{align*}
	&  \sum_{k\geq1}\sup_{x}\left\Vert \sum\nolimits_{y\in 
	\operatorname*{Sh}_{L}}\left(  \hat
{g}\Delta^{k}\right)  \left(  x,y\right)  \left(  \Delta\hat{\pi}\right)
\left(  y,\cdot\right)  \right\Vert _{1}
%\\
%& 
\leq\sum_{k\geq1}\sup_{x}\hat{g}(x,B_{1}(k))\sup_{z\in B_{1}(k)}\left\Vert
\Delta^{k+1}(z,\cdot)\right\Vert _{1}\\
&  \leq\sum_{k=2}^{20\log\log L+1}k\left(  \frac{1}{8}\right)  ^{k}%
+\sum_{k\geq20\log\log L+2}(\log L)^{6}\left(  \frac{1}{8}\right)  ^{k}<\frac
{1}{8}\,.%
\end{align*}
Therefore, it suffices to prove%
\begin{equation}
	\sup_{x}\left\Vert \sum\nolimits_{y\in \operatorname*{Sh}_{L}}
	\hat{g}\left(  x,y\right)
	\left(  \Delta\hat{\pi}1_{V_{L}}
\right)  \left(  y,\cdot\right)  \right\Vert _{1}%
\leq3/4. \label{Est_3/4}%
\end{equation}

From the second part 
of (\ref{k_0_large}) it follows that%
\[
\sup_{x\in V_{L}}\hat{g}\left(  x,\operatorname*{Sh}\nolimits_{L}\right)  
\leq10/9\,,\
\sup_{x\in \operatorname*{Sh}_{L}}
\pi_{V_{k_{0}r\left(  L\right)  }\left(  x\right)  }\left(
x,V_{L}\right)  \leq 1/10.
\]
%Furthermore we have assumed $\delta\leq1/32,$ so that 
By the third part
of (\ref{k_0_large}), and the choice $\delta<1/800<1/32$,
we get%
\[
\sup_{x\in \operatorname*{Sh}_{L}}
\Pi_{V_{k_{0}r\left(  L\right)  }\left(  x\right)  }\left(
x,V_{L}\right)  \leq \delta+17/32\leq 9/16.
\]
Combining that, we get%
\begin{align*}
	&  \sup_{x}\left\Vert \sum\nolimits_{y\in \operatorname*{Sh}_{L}}
	\hat{g}\left(  x,y\right)
\left(  \Delta\hat{\pi}1_{V_{L}}\right)  \left(  y,\cdot\right)  \right\Vert
_{1}\\
\leq   &  \sup_{x}\sum\nolimits_{y\in \operatorname*{Sh}_{L}}
\hat{g}\left(  x,y\right)
\Pi_{V_{k_{0}r\left(  L\right)  }\left(  y\right)  }\left(  y,V_{L}\right)
+\sup_{x}\sum\nolimits_{y\in \operatorname*{Sh}_{L}}
\hat{g}\left(  x,y\right)  \pi
_{V_{k_{0}r\left(  L\right)  }\left(  y\right)  }\left(  y,V_{L}\right) \\
\leq &  \frac{10}{9}\cdot \frac{9}{16}+
 \frac{10}{9}\cdot \frac{1}{10} <
%\frac{5}{8}+\frac{1}{9}<
\frac{3}{4},
\end{align*}
proving (\ref{Est_3/4}).
We conclude that
%on $G_{L}\cap\operatorname*{Good}\nolimits_{L},$
%\[
$\sup_{x\in V_{L}}\left\Vert 
	{\mathcal{R}}_{L}\hat\pi_{\Psi}
\left(  x,\cdot\right)  \right\Vert _{1}\leq
C\left(  \log L\right)  ^{-37/4}%
$
on $G_{L}\cap\operatorname*{Good}\nolimits_{L}.$
%\]
\end{proof}

\subsection{Green function estimates in a goodified environment
\label{Subsect_greengood}}
Before proceeding to analyze environments where bad regions are
present, we consider first ``goodified'' transition
kernels $\operatorname*{gd}\left(  \hat{\Pi}\right)$, 
c.f. (\ref{eq-160306a}).
We write
$\tilde{G}_{L}$ for the Green function corresponding to
this transition kernel. The goal of this section is to 
derive some estimates on $\tilde{G}_L$, which will be useful
in handling the event
$\left(  \operatorname*{Good}_{L}\cup\operatorname*{TwoBad}_{L}\right)
^{c}$.

Recall the range $\rho=\rho_{1,L}$, c.f. (\ref{Def_Rho1}), and consider the
collection
\begin{equation}
	\label{eq-160306b}
	\mathcal{D}_{L}=\left\{ V_{5\rho\left(
	x\right)  }\left(  x\right)  \,, x\in V_{L} \right\}\,.
\end{equation}
\begin{lemma}
\label{lem-220705a}
There exists a constant $c_{0}$ such that for 
all 
%$L$ large, and
$D\in\mathcal{D}_{L}$, $D\cap {\operatorname*{Shell}}_{L}(L/2)\neq \emptyset$,%
\begin{equation}
\tilde{G}_{L}(0,D)\leq c_{0}\left[  \frac{\mbox{\rm diam}(D)^{d-2}\left(
\max_{y\in D}d_{L}(y)\vee s(L)\right)  }{L^{d-1}}\right]  \,.
\label{eq-220705b}%
\end{equation}
Further, there exists a constant $c_{1}\geq1$ such that
for all  
$D\in\mathcal{D}_{L}$, 
\begin{equation}
\sup_{y\in V_{L}}\tilde{G}_{L}(y,D)\leq c_{1}\,. \label{eq-110805d}%
\end{equation}
\end{lemma}
\begin{proof}
[Proof of Lemma \ref{lem-220705a}]
%The proof is similar to what has been done in
%controlling $\Pi_L$.
We begin by establishing some auxiliary estimates for the unperturbed Green
function $\hat g=\hat g_L$. We first show that there is a constant $C$ such that
for any 
$D\in\mathcal{D}_{L}$, 
\begin{equation}
\label{eq-120805d}\sup_{y\in V_{L}} \hat g(y,D)=
\sup_{y\in D} \hat g(y,D)
\leq C\,.
\end{equation}
%Indeed, in proving (\ref{eq-120805d}), it is enough to consider $y\in D$. Fix
%a constant $\beta$ to be chosen below (see (\ref{eq-120805g})). For $D$ such
%that $D\cap{\operatorname*{Shell}\nolimits}_{L}(\beta s(L)) \neq\emptyset$,
For $D$ such
that $D\cap{\operatorname*{Shell}\nolimits}_{L}(2s(L)) \neq\emptyset$,
the estimate (\ref{eq-120805d}) 
%(with $C=C(\beta)$ depending on the choice of
%$\beta$) 
is an immediate consequence of parts a) and b) of Lemma
%\ref{Cor_Green}. If $D\cap{\operatorname*{Shell}\nolimits}_{L}(\beta
\ref{Cor_Green}. If $D\cap{\operatorname*{Shell}\nolimits}_{L}(2
s(L))=\emptyset$, then
%and $y\in D$, then
	$$\max_{y,z\in D}
\hat{g}(y,z)   = 
\max_{y\in D} 
\hat{g}(y,y)    \leq
1+ \max_{z\in \partial 
V_{\gamma s(L)}(x)}\hat g(z,y)
\leq 1 +\frac{C}{s(L)^d}
\,,$$
%$<++>G_{\gamma s(L)}(y,D)+\\
where $C$ depends on $\gamma$ and the second inequality follows
from part c) of Lemma \ref{Cor_Green}.
Summing over $z\in D$ completes the proof of 
(\ref{eq-120805d}). 
%
%
%
%, using Lemma \ref{Le_MainExit} in the
%second inequality, and the choice of $\gamma$ implying $10\gamma\leq1$, see
%(\ref{Def_Gamma}), we may find a constant $C_{1}$ independent of $\beta$ such
%that
%\begin{align*}
%\max_{y\in D}\hat{g}_{L}(y,D)  &  \leq1+ \max_{y\in D}G_{\gamma s(L)}(y,D)+\\
%&  \quad\quad\sum_{x\in{\operatorname*{Shell}}_{L}(2s(L))} P_{y}^{\mathrm{R}%
%W}(X_{\sigma^{\prime} }=x)P_{x}^{\mathrm{R}W}(T_{D}<T_{V_{L}}) \max_{y\in
%D}\hat{g}_{L}(y,D)\\
%&  \leq C + C_{1} \frac{(\beta+3) s(L)^{2} s(L)^{d-2}} {((\beta-3) s(L))^{d}}
%\max_{y\in D}\hat{g}_{L}(y,D)\,.
%\end{align*}
%Choosing $\beta>3$ large enough such that
%\begin{equation}
%\label{eq-120805g}C_{\beta}:=C_{1} (\beta+3)/(\beta-3)^{d}<1\,,
%\end{equation}
%we find that
%\[
%\max_{y\in D}\hat{g}_{L}(y,D) \leq C+C_{\beta}\max_{y\in D}\hat{g}%
%_{L}(y,D)\,,
%\]
%from which the conclusion (\ref{eq-120805d}) follows.

We next note that, for any $z\in V_{L}$,
\[
\hat g(z,D)\leq P_{z}^{\mathrm{RW}}(T_{D}<\tau_{V_{L}}) \max_{w\in D}
\hat g(w,D)\,.
\]
Applying (\ref{eq-120805d}) and Lemma \ref{Le_MainExit}, we deduce that for
some constant $C_{0}$,
\begin{equation}
\label{eq-120805e}\hat g(z,D)\leq C_{0}\left[  \frac
{\mbox{\rm diam}(D)^{d-2} d_{L}(z) \max_{y\in D} d_{L}(y)}{d(z,D)^{d}}
\wedge 1\right]  \,.
\end{equation}

We now turn to proving (\ref{eq-110805d}).
Write the perturbation expansion
\begin{equation}
\label{eq-120805h}\tilde G_{L}(z,D)-\hat g(z,D)= \sum_{k\geq1}
\sum_{y,y^{\prime},w} \hat g(z,y)\Delta^{k}(y,y^{\prime})\hat\pi
(y^{\prime},w)\hat g(w,D) +\mbox{\rm NL}\,,
\end{equation}
where $\mbox{\rm NL}$ denotes the nonlinear term in the perturbation
expansion, that is
\begin{equation}
\mbox{\rm NL}= \sum_{m=2}^{\infty}\sum_{k_{1},\ldots,k_{m}=0}^{\infty} \left(
\hat{g}_{L}\Delta^{k_{1}}\Delta\hat{\pi}\right)  \cdot\ldots\cdot\left(
\hat{g}_{L}\Delta^{k_{m-1}}\Delta\hat{\pi}\right)  \left(  \hat{g}_{L}%
\Delta^{k_{m}}\Delta\hat g(\cdot,D)\right)  \,. \label{NonLinearSplit11}%
\end{equation}

We first handle the linear term in (\ref{eq-120805h}). 
Using (\ref{eq-120805d}), part d) of Lemma
\ref{Cor_Green}, and Lemma \ref{lem-090805},
we see that in a goodified environment,
\begin{equation}
\label{eq-120805a1}|\sum_{k\geq1} \sum_{y,y^{\prime},w: d_{L}(y^{\prime})\geq
k_{0}r(L)} \hat g(z,y)\Delta^{k}(y,y^{\prime})\hat\pi(y^{\prime},w)\hat
g(w,D)| \leq \frac{C(\log L)^{6-9}}{ 8^{k}}\,,
\end{equation}
and
\begin{equation}
\label{eq-120805a2}| \sum_{y,y^{\prime},w} \hat g(z,y)\Delta
^{k}(y,y^{\prime})\hat\pi(y^{\prime},w)\hat g(w,D)| \leq \frac{C(\log L)^{6}}
{8^{k}}\,.
\end{equation}
From (\ref{eq-120805a2}) it follows that
\begin{equation}
\label{eq-120805a3}|\sum_{k\geq20\log\log L} \sum_{y,y^{\prime},w} \hat
g(z,y)\Delta^{k}(y,y^{\prime})\hat\pi(y^{\prime},w)\hat g(w,D)| \leq C
(\log L)^{-9}\,.
\end{equation}
On the other hand, if $d_{L}(y^{\prime})\leq k_{0} r(L)$ and $\Delta
^{k}(y,y^{\prime}) >0$ then, as in the proof of Proposition
\ref{prop-erwin160605}, $d_{L}(y)\leq(4/3)^{k} k_{0} r(L)$. Using parts a),b)
of Lemma \ref{Cor_Green}, we get that for $k\leq20 \log\log L$,
\begin{equation}
\label{eq-120805a4}| \sum_{y,y^{\prime},w: d_{L}(y^{\prime})\leq k_{0}r(L)}
\hat g(z,y)\Delta^{k}(y,y^{\prime})\hat\pi(y^{\prime},w)\hat g(w,D)|
\leq Ck (1/8)^{k}%
\end{equation}
Combining (\ref{eq-120805a1}), (\ref{eq-120805a3}) and (\ref{eq-120805a4}), we
conclude that
\[
\sup_{z\in V_{L}} \sum_{k\geq1} \sum_{y,y^{\prime},w} \hat g%
(z,y)\Delta^{k}(y,y^{\prime})\hat\pi(y^{\prime},w)\hat g(w,D) \leq C\,.
\]
%
%We next divide $V_L$ into boxes $B_i$ of centers $x_i$,
%each of side $\rho_L(x_i)$.
The term involving $\mbox{\rm NL}$ is handled by recalling that
\[
\sup_{x}\sum_{k\geq0}\left\Vert \left(  \hat{g}_L
\Delta^{k}\Delta\hat{\pi
}\right)  \left(  x,\cdot\right)  \right\Vert _{1}\leq15/16\,,
\]
see (\ref{Est_15/16}). We then conclude, using (\ref{eq-120805d}), that
(\ref{eq-110805d}) holds.

To prove (\ref{eq-220705b}), our starting point is the perturbation expansion
(\ref{eq-120805h}). Again, the main contribution is the linear term. 
From (\ref{eq-120805a2}) one deduces that
%One has
%\[
%\sum_{y,y^{\prime},w} \hat g(0,y)\Delta^{k}(y,y^{\prime})\hat\pi
%(y^{\prime},w)\hat g(w,D) \leq C (\log L)^{6} (1/8)^{k}\,.
%\]
%Hence, 
there exists a constant $c_{d}$ such that for all $L$ large,
\begin{equation}
\label{eq-130805a}\sum_{k\geq c_{d} \log\log L} \sum_{y,y^{\prime},w} \hat
g(0,y)\Delta^{k}(y,y^{\prime})\hat\pi(y^{\prime},w)\hat g(w,D)
\leq\left(  \frac{r(L)}{L}\right)  ^{d-1}\,.
\end{equation}
We divide the sum in the linear term according to the location of $w$
with respect to $D$, writing
\begin{equation}
\label{eq-130805b}\sum_{y,y^{\prime},w} \hat g(0,y)\Delta^{k}(y,y^{\prime
})\hat\pi(y^{\prime},w)\hat g(w,D) =\sum_{y,y^{\prime}}\hat g(0,y)
\Delta^{k}(y,y^{\prime}) \sum_{j=1}^{2} \sum_{w\in B_{j}} \hat
\pi(y^{\prime},w)\hat g(w,D)\,,
\end{equation}
where
\[
B_{1}=\{w\in V_{L}: d(w,D)\leq L/8\}\,,\quad B_{2}=\{w\in V_{L}:
d(w,D)> L/8\}\,.
\]

Considering the term involving $B_{1}$, for $k<c_{d} \log\log L$ the
summation over $y$ extends over a subset of $V_{L}$ that is covered by at most
$Ck^{d}$ elements of ${\mathcal{D}}_{L}$, all inside ${\operatorname*{Shell}%
\nolimits}_{L}(3L/4)$. Thus, for such $k$, using (\ref{eq-120805e})
to bound $\hat g(0,y)$, (\ref{eq-120805d}) to bound $\hat g(w,D)$,
and Lemma \ref{lem-090805}, we get
\[
\sum_{y,y^{\prime}}\hat g(0,y) \Delta^{k}(y,y^{\prime}) \sum_{w\in
B_{1}} \hat\pi(y^{\prime},w)\hat g(w,D) \leq
%16 C C_0
C \left(  \frac{1+\gamma}{8}\right)  ^{k} k^{d} \frac{\mbox{\rm diam}(D)^{d-2}
\max_{y\in D} d_{L}(y)}{L^{d-1}}
\]
and hence
\begin{equation}
\label{eq-130805c}\sum_{k\leq c_{d} \log\log L} \sum_{y,y^{\prime}} \hat
g(0,y)\Delta^{k}(y,y^{\prime})\sum_{w\in B_{1}} \hat\pi(y^{\prime
},w)\hat g(w,D) \leq C \frac{\mbox{\rm diam}(D)^{d-2} \max_{y\in D}
d_{L}(y)}{L^{d-1}}\,.
\end{equation}
The term involving $w\in B_{2}$ is simpler: indeed, one has in that case
that $\hat g(w,D)$ satisfies, by (\ref{eq-120805e}), the required bound,
whereas for $k< c_d \log\log L$, using
(\ref{eq-120805d}),
\[
\sum_{y: \exists y^{\prime}\,\mbox{\rm with}\, \Delta^{k}(y,y^{\prime})
\hat\pi(y^{\prime},w)>0} \hat g(0,y) \leq C k^{d}\,,
\]
yielding
\begin{align}
\label{eq-130805d} &  \sum_{k\leq c_{d} \log\log L} \sum_{y,y^{\prime}} \hat
g(0,y)\Delta^{k}(y,y^{\prime})\sum_{w\in B_{2}} \hat\pi(y^{\prime
},w)\hat g(w,D)\nonumber\\
&  \leq C \sum_{k\leq c_{d} \log\log L} k^{d} (1/8)^{k} \frac
{\mbox{\rm diam}(D)^{d-2} \max_{y\in D} d_{L}(y)}{L^{d-1}}\,.
\end{align}
Combining (\ref{eq-130805a}), (\ref{eq-130805c}) and (\ref{eq-130805d})
results in the required control on the linear term in (\ref{eq-120805h}). The
nonlinear term is even simpler and similar to the handling of the nonlinear
term when estimating $\hat g(z,D)$.
\end{proof}

\subsection{Presence of bad regions\label{Subsect_BadPoints}}
On $\left(  \operatorname*{Good}_{L}\cup\operatorname*{TwoBad}_{L}\right)
^{c},$ it is clear that for some 
$D\in\mathcal{D}_{L}$, c.f.
(\ref{eq-160306b}),
we have
\begin{equation}
B_{L}\subset D
%V_{5\rho\left(  x\right)  }\left(  x\right)  .
\label{BL_included_in_ball}%
\end{equation}
%In particular, 
%for $D\in\mathcal{D}_{L},$, c.f.
%(\ref{eq-160306b}),
We write $\operatorname*{Bad}_{L}\left(  D\right)  $ for the event that
$\left\{  B_{L}\subset D\right\}  ,$ and $\operatorname*{Bad}_{L}^{(i)}\left(
D\right)  $ for the event that $\left\{  B_{L}^{(i)}\subset D\right\}  $,
$i=1,2,3,4$.
The main aim of this section is to prove the following. 
\begin{proposition}
\label{prop-erwin220705} There exists a $\delta_{0}\leq1/800\bar c$ such that
if $\delta<\delta_{0}$, and if $\operatorname*{Cond}\left(  L_{1}%
,\delta\right)  $ holds for a given $L_{1}$, and if $L\leq L_{1}\left(  \log
L_{1}\right)  ^{2}$ and ${\Psi}\in\mathcal{M}_{L}$, then,
for $i=1,2,3,4$,
$$\sup_{D\in\mathcal{D}_{L}}\mathbb{P}\left(  D_{L,\Psi}(0)
\geq\left(  \log L\right)  ^{-9+\frac{9(i-1)}{4}}
,\ \operatorname*{Bad}\nolimits_{L}^{\left(  i\right)  }\left(  D\right)
\right) 
\leq 
 \frac{\exp\left[  -\left(  \log L\right)  ^{2}\right]}{100}
\,.
$$
\end{proposition}
\begin{proof}
[Proof of Proposition \ref{prop-erwin220705}]
We start with the case when $D$ is \textquotedblleft not
near\textquotedblright\ the boundary, meaning that $D\subset V_{ L/2}$. We
write $D=V_{5\rho\left(  x_{0}\right)  }\left(  x_{0}\right)  =V_{5\gamma
s\left(  L\right)  }\left(  x_{0}\right)  .$ By Lemma \ref{Cor_Green} c), we
can find a constant $K$ (not depending on $L,x_{0}$), such that for any point
$x\notin\widetilde{D}\overset{\mathrm{def}}{=}V_{5K\gamma s\left(  L\right)
}\left(  x_{0}\right)  ,$ and all $L$ large, one has $\hat{g}\left(
x,D\right)  \leq1/10.$ We modify now the transition probabilities $\hat{\Pi
},\hat{\pi}$ slightly, when starting in $x\in D,$ by defining%
\begin{equation}
\widetilde{\Pi}\left(  x,\cdot\right)  \overset{\mathrm{def}}{=}\left\{
\begin{array}
[c]{cc}%
\operatorname*{ex}\nolimits_{\widetilde{D}}\left(  x,\cdot;\hat{\Pi}\right)  &
\mathrm{\ \mathrm{for\ }}x\in D\\
\hat{\Pi}\left(  x,\cdot\right)  & \mathrm{\ \mathrm{for\ }}x\notin D
\end{array}
\right.  , \label{Def_Enlargement}%
\end{equation}
and similarly we define $\widetilde{\pi}.$ (Remark that this destroys somewhat
the symmetry, when $x\neq x_{0},$ but this is no problem below). Clearly,
these transition probabilities have the same exit distribution from $V_{L}$ as
the one used before. If we write $\widetilde{g}$ for the Green's function on
$V_{L}$ of $\widetilde{\pi},$ we have $\widetilde{g}\left(  x,y\right)
=\hat{g}\left(  x,y\right)  $ for $y\notin\widetilde{D},$ and all $x,$ whereas
$\widetilde{g}\left(  x,y\right)  \leq\hat{g}\left(  x,y\right)  $ for
$y\in\widetilde{D}.$ In particular, we have%
\begin{equation}
\sup_{x\not \in \widetilde{D}}\widetilde{g}\left(  x,D\right)  \leq1/10.
\label{eq-260805a}%
\end{equation}
Writing down the perturbation expansion (\ref{eq-120306d}) using the kernels
$\tilde \Pi$ and $\tilde \pi$, we have%
\[
\left(  \left[  \Pi-\pi\right]  \hat{\pi}_{{\Psi}}\right)
=\sum_{m=1}^{\infty}\sum_{k_{1},\ldots,k_{m}=0}^{\infty}\left(  \widetilde
{g}\Delta^{k_{1}}\Delta\hat{\pi}\right)  \cdot\ldots\cdot\left(  \widetilde
{g}\Delta^{k_{m-1}}\Delta\hat{\pi}\right) 
\left(  \widetilde{g}\Delta^{k_{m}%
}\Delta\phi\right)  ,
\]
where $\Delta$ now uses the modified transitions, that is $\Delta
(x,y)=\tilde{\Pi}(x,y)-\tilde{\pi}(x,y)$, but remark that for $x\notin D,$
$\Delta\left(  x,\cdot\right)  $ is the same as before, and
that always $\Delta \tilde \pi=\Delta \hat \pi$. Also, $\phi$ is
modified accordingly.

We first estimate the part with $m=1$. In anticipation of what follows, we
consider an arbitrary starting point $x\in V_{L}$. Put $k=k_{1}+1.$ The part
of the sum%
\[
\sum_{y}\sum_{x_{1},\ldots,x_{k}}\widetilde{g}\left(  x,x_{1}\right)
\Delta\left(  x_{1},x_{2}\right)  \cdot\ldots\cdot\Delta\left(  x_{k}%
,y\right)  \phi\left(  y,\cdot\right)
\]
where all $x_{j}\notin D,$ is estimated in Section \ref{SubSect_GoodLinear},
and the probability that it exceeds $\left(  \log L\right)  ^{-9}/3$ is
bounded by $\exp\left[  -\left(  \log L\right)  ^{2}\right]  /100.$ If an
$x_{j}\in D,$ then the sum over $x_{j+1}$ extends only to points outside
$\widetilde{D},$ and therefore, the sum over $x_{j+1},x_{j+2},\ldots,x_{j+K}$
is running only over points outside $D.$ Therefore%
\begin{equation}
\label{eq-011005}\sup_{x_{j}\in D}\sum_{x_{j+1},\ldots x_{j+K}}\left\vert
\Delta\left(  x_{j},x_{j+1}\right)  \cdot\ldots\cdot\Delta\left(
x_{j+K},x_{j+K+1}\right)  \right\vert \leq2\delta^{K}.
\end{equation}
Further, let $j$ denote the smallest index such that $x_{j}\in D$. Let
\[
\mathcal{X}_{j}:=\left\{  x_{1}:\Delta(x_{1},x_{2})\cdots\Delta(x_{j-1}%
,x_{j})\right\}  >0\,.
\]
Then $\max_{x_{1}\in\mathcal{X}_{j}}d(x_{1},D)\leq5j\gamma s(L)$. For $j<(\log
L)^{2}$ it follows that $\mathcal{X}_{j}\subset V_{L-s(L)}$ and therefore, by
(\ref{Est_BoundaryReach1}),
%\ref{Cor_Green} c), for some constant $C$ (recall that $\gamma<1$!),
$\max_{x\in V_{L}}\tilde{g}(x,\mathcal{X}_{j})\leq C j^{d}$. Thus,
\begin{equation}
\label{eq-011005b}\left\vert \sum\nolimits_{x_{1},\ldots,x_{j}}\widetilde
{g}\left(  x,x_{1}\right)  \Delta\left(  x_{1},x_{2}\right)  \cdots
\Delta\left(  x_{j-1},x_{j}\right)  \right\vert \leq C\delta^{j-1}j^{d}\,.
\end{equation}
On the other hand, for $j\geq(\log L)^{2}$ one has
(recalling that $x_i\not\in D$ for $i<j$, and applying
part d) of Lemma \ref{Cor_Green} together with Lemma \ref{lem-090805}),
\[
|\sum_{x_{1},\ldots,x_{j}}\widetilde{g}\left(  x,x_{1}\right)  \Delta\left(
x_{1},x_{2}\right)  \cdots\Delta\left(  x_{j-1},x_{j}\right)  |\leq
C(1/8)^{j}(\log L)^{6}\,.
\]
Therefore, using (\ref{eq-011005}),
%for some constant $\bar c$,%
\[
\sum_{j=1}^{\infty}\left\vert \sum_{x_{1},\ldots,x_{j-1}\not \in D,x_{j}\in
D}\widetilde{g}\left(  x,x_{1}\right)  \Delta\left(  x_{1},x_{2}\right)
\cdots\Delta\left(  x_{j-1},x_{j}\right)  \right\vert \leq C\,.
\]
If $x_{k}\notin D,$ then on the event $B_L\subset D$, using
part a) of Proposition \ref{Prop_LipshitzPhi},
it holds that\\ $\left\Vert \sum_{y}\Delta\left(  x_{k},y\right)
\phi\left(  y,\cdot\right)  \right\Vert _{1}$ $
\leq C\left(  \log L\right)
^{-12}.$ On the other hand, if $x_{k}\in D,$ then%
\begin{equation}
\left\Vert \sum\nolimits_{y}\Delta\left(  x_{k},y\right)  \phi\left(
y,\cdot\right)  \right\Vert _{1}\leq C\gamma K\left(  \log L\right)
^{-12+2.25i}. \label{ExitFromBad}%
\end{equation}
Combining all the above, we conclude that for some constant $c_{2}$ it holds
that
\[
|\sum_{y,z}\sum_{x_{1},\ldots,x_{k}}\widetilde{g}\left(  x,x_{1}\right)
\Delta\left(  x_{1},x_{2}\right)  \cdot\ldots\cdot\Delta\left(  x_{k}%
,y\right)  \phi\left(  y,z\right)  |\leq c_{2}\gamma K(\log L)^{-12+
2.25i}\,.
\]
It follows that%
\begin{equation}
\label{eq-011005d}\left\Vert \sum\nolimits_{x_{1},\ldots,x_{k}}^{\prime
}\widetilde{g} \left(  0,x_{1} \right)  \Delta\left(  x_{1},x_{2}\right)
\cdot\ldots\cdot\left(  \Delta\phi\right)  \left(  x_{k},\cdot\right)
\right\Vert _{1}\leq\left(  \log L\right)  ^{-11.5+2.25i},
\end{equation}
where $\sum^{\prime}$ denotes summation where at least one $x_{j}$ is in $D.$
(We note that for $i=1,2,3$, one does not 
need to use the $K$-enlargement and
modification of the transition probabilities, as a factor
$\delta$ is caught for each $\|\Delta\|_1$.)

The case $m\geq2$ is handled with an evident modification of the above
procedure, using the estimate (\ref{eq-260805a}). Indeed, let $D^{\prime
}=\{z\in V_{L}:d(z,\tilde D)\leq2\gamma s(L)\}$. A repeat of the previous
argument shows that
\[
\sup_{x} \sum_{k=4}^{\infty}\sum_{x_{k}}|\sum_{\overset{ x_{1},\ldots
,x_{k-1}:}{ \exists j\leq k, x_{j}\in D^{\prime}}}\widetilde{g} \left(
x,x_{1} \right)  \Delta\left(  x_{1},x_{2}\right)  \cdot\ldots\cdot
\Delta\left(  x_{k-2},x_{k-1}\right)  \hat\pi\left(  x_{k-1},x_{k}\right)
|\leq C\delta
\]
while
\[
\sup_{x} \sum_{x_{3}}|\sum_{\overset{x_{1},x_{2}:}{ \exists j\leq3, x_{j}\in
D^{\prime}}}\widetilde{g} \left(  x,x_{1} \right)  \Delta\left(  x_{1}%
,x_{2}\right)  \hat\pi\left(  x_{2},x_{3}\right)  | \leq\left\{
\begin{array}
[c]{ll}%
\frac{2}{10}\,, & x\not \in D^{\prime}\\
C\,, & x\in D^{\prime}\,,
\end{array}
\right.
\]
and, by the computation in Section \ref{SubSect_Nonlinear}, c.f.
(\ref{Est_15/16}),
\[
\sup_{x} \sum_{k=3}^{\infty}\sum_{x_{k}\not \in D^{\prime}}|\sum
_{\overset{x_{1},\ldots,x_{k-1}:}{ x_{j}\not \in D^{\prime}}}\widetilde{g}
\left(  x,x_{1} \right)  \Delta\left(  x_{1},x_{2}\right)  \cdot\ldots
\cdot\Delta\left(  x_{k-2},x_{k-1}\right)  \hat\pi\left(  x_{k-1}%
,x_{k}\right)  | \leq\frac{15}{16}\,.
\]
Hence, we conclude that always,
\begin{equation}
\sup_{x}\sum_{k\geq0}\left\Vert \left(  \hat{g}\Delta^{k}\Delta\hat{\pi
}\right)  \left(  x,\cdot\right)  \right\Vert _{1}\leq C, \label{Est_15/16C}%
\end{equation}
and for all $\delta$ small,
\begin{equation}
\sup_{x}\sum_{k_{1},k_{2}\geq0} \left\Vert \left(  \hat{g}\Delta^{k_{1}}%
\Delta\hat{\pi}\right)  \left(  \hat{g}\Delta^{k_{2}}\Delta\hat{\pi}\right)
\left(  x,\cdot\right)  \right\Vert _{1}\leq\frac{16}{17}\,.
\label{Est_15/162}%
\end{equation}
Together with the computation for $m=1$, c.f. (\ref{eq-011005d}) when
$D^{\prime}$ is visited, and Proposition \ref{prop-erwin160605} when it is
not, this completes the proof of Proposition \ref{prop-erwin220705} in case
$D\subset V_{L/2}$.

We next turn to 
$D\cap {\operatorname*{Shell}}_{L}(L/2)\neq \emptyset$.
Recall the Green function $\tilde G_{L}$ of the goodified environment,
introduced above Lemma \ref{lem-220705a}. Let $\Pi^{g}_L$ denote the
exit distribution $\Pi_L$ 
from $V_L$ with the environment replaced by the goodified
environment. Let $\Delta^{g}=1_{D}(\Pi_L-\Pi^{g}_L)$. 
The perturbation expansion (\ref{Pert1}) then gives
\[
[\Pi_L-\Pi^{g}_L](z)=\sum\tilde G_{L}(0,y)\Delta^{g}(y,y^{\prime
})\Pi_L(y^{\prime},z)\,,
\]
and thus, using part a) of Lemma \ref{lem-220705a} in the second inequality,
\begin{equation}
\label{eq-220705c}\|\Pi_{{L}}-\Pi_{{L}}^{g}\|_1 \leq2 \tilde G_{L}(0,D)
\leq C 
%3\cdot10^{d} c_{0} 
\frac{s(L)^{d-2}}{L^{d-2}} \leq C (\log L)^{3(2-d)}
\,,
\end{equation}
This completes the proof in case $i=4$ (and also $i=1,2,3$ if $d\geq5$, although
we do not use this fact).

Consider next the case $i=1,2,3$ (and $d=3,4$).
Rewrite the perturbation expansion as
\begin{equation}
\lbrack\Pi_L-\Pi^{g}_L](z)=\sum_{k\geq1}\sum_{y}\tilde{G}%
_{L}\left(  \Delta^{g}\right)  ^{k}(0,y)\left(  \hat{\Pi}
\tilde{G}_{L}\Delta^{g}\Pi^{g}_L\right)  (y,z)\,
\end{equation}
In particular, using Lemma \ref{lem-090805} and part b) of Lemma
\ref{lem-220705a},
\begin{align}
\Vert\Pi_{L}-\Pi_{L}^{g}\Vert_1 &  \leq C\tilde{G}_{L}(0,D)\sum_{k\geq
1}(1/8)^{k}(\log L)^{-9+2.25i} \sup_{y^{\prime}\in V_{L}}\tilde{G}%
_{L}(y^{\prime},D)\nonumber\label{eq-110805b}\\
&  \leq C (\log L)^{3(2-d)}(\log L)^{-9+ 2.25i}\leq(\log L)^{-11.5+2.25 i}\,.
\end{align}
\end{proof}
\section{The non-smoothed exit estimate\label{Sect_NonSmooth}}
The aim of this section is to prove the following.
\begin{proposition}
\label{Prop_Nonsmooth} There exists $0<\delta_{0}\leq1/2$ such that for
$\delta\leq\delta_{0},$ there exist $L_{0}\left(  \delta\right)  $ and
$\varepsilon_{0}\left(  \delta\right)  $ such that if $L_{1}\geq L_{0}$ and
$\varepsilon\leq\varepsilon_{0}$, then $\operatorname*{Cond}\left(
L_{1},\delta\right)  ,$ and $L\leq L_{1}\left(  \log L_{1}\right)  ^{2}$ imply%
\[
\mathbb{P}\left(  \left\Vert \Pi_{L}\left(  0,\cdot\right)  -\pi_{L}\left(
0,\cdot\right)  \right\Vert _{1}\geq\delta\right)  \leq\frac{1}{10}\exp\left[
-\left(  \log L\right)  ^{2}\right]\,.
\]
\end{proposition}
Before starting the proof, we provide a sketch of the main idea. 
As with the smoothed estimates, the starting point is the perturbation
expansion (\ref{eq-120306d}). In contrast to the proof in Section
\ref{Sect_Smooth}, however, no smoothing is provided by the kernel
$\hat \pi_\Psi$, and hence the lack of control of the exit 
measure in the last step of the coarse-graining scheme 
$\mathcal{S}_1$ does not allow one 
to propagate the estimate on $D_{L,0}(0)$. 
This is why we need to work with the scheme $\mathcal{S}_2$ introduced 
in Definition \ref{Def_CoarseGrainingScheme}.
Using $\mathcal{S}_{2}$ means
that we refine the coarse graining 
scale up to the boundary, and when carrying out the perturbation
expansion, less smoothing is gained from the coarse graining
for steps near the boundary. 
The drawback of $\mathcal{S}_2$
is that the presence of many bad regions
close to the boundary is unavoidable. We will however show that these regions
are
rather sparse, so that with high enough probability, the RWRE
avoids the bad regions. As in Section \ref{Subsect_BadPoints},
this will be achieved by an appropriate
estimate on the Green function in a ``goodified'' environment.

\noindent
\begin{proof}
[Proof of Proposition \ref{Prop_Nonsmooth}]
We use the coarse graining scheme $\mathcal{S}_{2}$ from Definition
\ref{Def_CoarseGrainingScheme}, but we stick to the notations before, so
$\hat{\pi}=\hat{\pi}_{\mathcal{S}_{2},L}$, etc. Using $\mathcal{S}_{2}$ means
that we refine the coarsening scale up to the boundary. In particular, 
 $h_{L}\left(  x\right)  =\gamma d_{L}\left(  x\right)  $
for all $x$ with $d_{L}\left(  x\right)  \leq s\left(  L\right)  /2,$ and
$\hat{\pi}\left(  x,\cdot\right)  $ is obtained by averaging exit
distributions from balls with radii between $\gamma d_{L}\left(  x\right)  $
and $2\gamma d_{L}\left(  x\right)  $ ($\gamma$ from (\ref{Def_Gamma})). If
$d_{L}\left(  x\right)  <1/2\gamma,$ then there is no coarsening at all, and
$\hat{\pi}\left(  x,\cdot\right)  =p^{\mathrm{RW}}\left(  x,\cdot\right)  .$

To handle the presence of many bad regions near the boundary, we 
introduce the 
layers%
\begin{equation}
	\label{eq-110206b}
\Lambda_{j}\overset{\mathrm{def}}{=}\operatorname*{Shell}\nolimits_{L}\left(
2^{j-1},2^{j}\right)  ,
\end{equation}
for $j=1,\ldots,J_{1}\left(  L\right)  \overset{\mathrm{def}}{=}\left[
\frac{\log r\left(  L\right)  }{\log2}\right]  +1,$ so that
\begin{equation}
\operatorname*{Shell}\nolimits_{L}\left(  r\left(  L\right)  \right)
\subset\bigcup\nolimits_{j\leq J_1(L)}\Lambda_{j}\subset\operatorname*{Shell}%
\nolimits_{L}\left(  2r\left(  L\right)  \right)  . \label{Layers_in_Shells}%
\end{equation}
We subdivide each $\Lambda_{j}$ into 
subsets $D_{1}^{\left(  j\right)  },D_{2}^{\left( 
 j\right)  },\ldots,D_{N_{j}}^{\left(  j\right)  }$ of diameter
$\leq\sqrt{d}2^{j},$ where%
\begin{equation}
C^{-1} 
\left(  L2^{-j}\right)  ^{d-1}
\leq 
N_{j}\leq C\left(  L2^{-j}\right)  ^{d-1}. \label{Est_NoLayerBoxes}%
\end{equation}
The collection
of these subsets is denoted by $\mathcal{L}_{j}.$ $\mathcal{L}_{j}$ is
split into disjoint $\mathcal{L}_{j}^{\left(  1\right)  },\ldots
,\mathcal{L}_{j}^{\left(  R\right)  },$ such that for any $m$ one has
\begin{equation}
d\left(  D,D^{\prime}\right)  >5\gamma2^{j},\ \forall D,D^{\prime}%
\in\mathcal{L}_{j}^{\left(  m\right)  }, \label{MinimalDistance}%
\end{equation}%
\begin{equation}
N_{j}^{\left(  m\right)  }\overset{\mathrm{def}}{=}\left\vert \mathcal{L}%
_{j}^{\left(  m\right)  }\right\vert \geq N_{j}/2R. \label{Number_of_Cells}%
\end{equation}
We can do that in such a way that $R\in\mathbb{N}$ depends only on the
dimension $d$ (recall that $\gamma$ is fixed by 
(\ref{Def_Gamma}) once the dimension is fixed).

\begin{figure}
\begin{picture}(10,200)(-80,0)\input{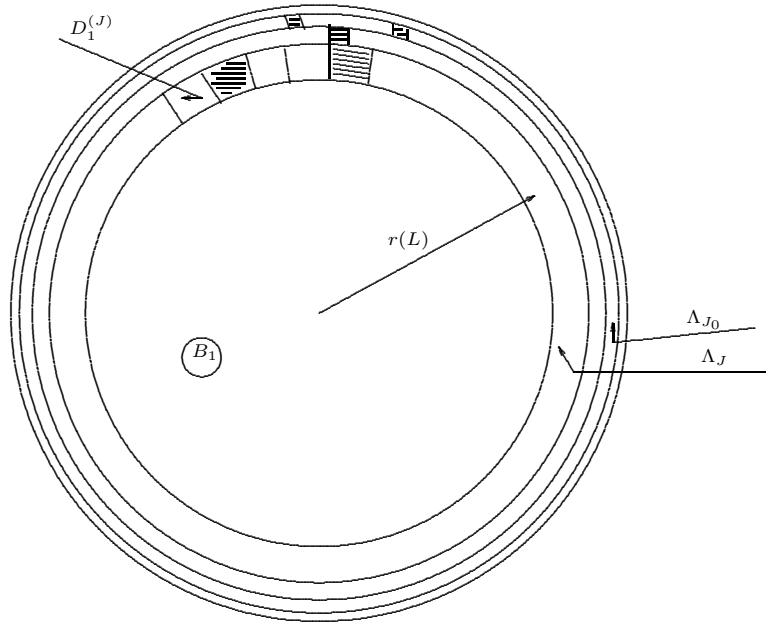}
\end{picture}
\caption{The layers $\Lambda_j$. Bad regions (excluding $B_1\neq\emptyset$)
are shaded.}
\end{figure}

For $x\in \cup_{j=1}^{J_1(L)}\Lambda_j$,  
we modify the definition of 
$\operatorname*{Good}_{L}$ in order to adapt it to the smoothing
scheme $\mathcal{S}_{2}$.
Thus, we set
%, 
%for $x\in \cup_{j=1}^{J_1(L)}\Lambda_j$:  
%for $i=1,2,3$, we keep
%$B_{L}^{\left(  i\right)  }$ as above (\ref{Def-BL}),
%while for $i=4$ we set
%are the set of points $x\notin
%\operatorname*{Shell}\nolimits_{L}\left(  r\left(  L\right)  \right)  $ such
%that for \textit{some} $r\in\lbrack h_{L}(x),2h_{L}(x)]$, one has
%$D_{r,h_{L}(x)}\left(  x\right)  >\left(  \log L\right)  ^{-11.25+2.25i},\ $
%but for all $r\in\lbrack h_{L}(x),2h_{L}(x)],$ $D_{r,h_{L}(x)}\left(
%x\right)  \leq\left(  \log L\right)  ^{-9+2.25 i},$ and $D_{r}^{0}\left(
%x\right)  \leq\delta.$ 
$\hat B_{L}^{\left(  4\right)  }$ to consist of the
union (over $j=1,\ldots,J_1(L)$) of 
points $x\in \Lambda_j$
which 
%for $d_{L}\left(  x\right)  >r\left(  L\right)  $ 
have the property that
$D_{\gamma d_L(x),0  }\left(  x\right)
\geq\delta$. We also write
$\widehat{ 
\operatorname*{Good}_{L}}=V_L\setminus \hat B_{L}^{\left(  4\right)  }$.

If $B\in\mathcal{L}_{j},$ we write $\operatorname*{Bad}\left(  B\right)  $ for
the event $\left\{  B\not\subset\widehat 
{\operatorname*{Good}_{L}}\right\}  .$ Remark
that%
\begin{equation}
	\label{eq-110206c}
\mathbb{P}\left(  \operatorname*{Bad}\left(  B\right)  \right)  \leq
C2^{\left(  d+1\right)  j}\exp\left[  -\frac{10}{13}
\log^{2}\left(  \gamma2^{j-1}\right)
\right]  
\leq\exp\left[  -j^{5/3}\right]  \overset{\mathrm{def}}{=}p_{j}.
\end{equation}
for $j\geq J_{0},$ $J_{0}$ appropriately chosen (depending on $d$).

We set%
\[
X_{j}^{\left(  m\right)  }\overset{\mathrm{def}}{=}\sum_{D\in\mathcal{L}%
_{j}^{\left(  m\right)  }}1_{\operatorname*{Bad}\left(  D\right)  }%
,\ X_{j}\overset{\mathrm{def}}{=}\sum_{m=1}^{R}X_{j}^{\left(  m\right)  }.
\]
Due to (\ref{MinimalDistance}), the events $\operatorname*{Bad}\left(
D\right)  ,\ D\in\mathcal{L}_{j}^{\left(  m\right)  },$ are independent.
Remark that $p_{j}<j^{-3/2}\leq1/2$ for all $j\geq2.$ From a standard coin
tossing estimate via Chebycheff's inequality, we get
\[
\mathbb{P}\left(  X_{j}^{\left(  m\right)  }\geq j^{-3/2}N_{j}^{\left(
m\right)  }\right)  \leq\exp\left[  -N_{j}^{\left(  m\right)  }I\left(
j^{-3/2}\mid p_{j}\right)  \right]  
\]
with $I\left(  x\mid p\right)  \overset{\mathrm{def}}{=}x\log\left(
x/p\right)  +\left(  1-x\right)  \log\left(  \left(  1-x\right)  /\left(
1-p\right)  \right)  $, and
\[
I\left(  j^{-3/2}\mid p_{j}\right)  \geq-\frac{3}{2}j^{-3/2}\log
j+j^{-3/2}j^{5/3}-\log2\geq2Rj^{1/7}%
\]
for $j\geq J_0$, if $J_{0}$ is large enough. Therefore%
\begin{align*}
\mathbb{P}\left(  X_{j}\geq j^{-3/2}N_{j}\right)   &  \leq R\max_{1\leq m\leq
R}\mathbb{P}\left(  X_{j}^{\left(  m\right)  }\geq j^{-3/2}N_{j}^{\left(
m\right)  }\right) \\
&  \leq R\exp\left[  -\left(  L2^{-j}\right)  ^{d-1}j^{1/7}\right]  \leq
R\exp\left[  -\frac{1}{C}\left(  \log L\right)  ^{20}j^{1/7}\right]
\end{align*}
for $J_{0} \leq j\leq J_1\left(  L\right)  ,$ $L$ large
enough (implied by $L_{0}$ large enough). Using this, we get for $L\geq L_0$,%
\[
\sum_{J_{0}  \leq j\leq J_1\left(  L\right)  }%
\mathbb{P}\left(  X_{j}\geq j^{-3/2}N_{j}\right)  \leq\frac{1}{20}\exp\left[
-\left(  \log L\right)  ^{2}\right]\,,
\]
increasing further $J_0$ and $L_0$ if necessary.
Setting%
\[
\operatorname*{ManyBad}\nolimits_{L}\overset{\mathrm{def}}{=}\bigcup
\nolimits_{J_{0}  \leq j\leq J_1\left(  L\right)  }\left\{
X_{j}\geq j^{-3/2}N_{j}\right\}  \cup\operatorname*{TwoBad}\nolimits_{L},
\]
we get%
\begin{align}
	\mathbb{P}\left(  \operatorname*{ManyBad}\nolimits_{L}\right)   &  \leq
\frac{1}{20}\exp\left[  -\left(  \log L\right)  ^{2}\right]  +\exp\left[
-1.2\left(  \log L\right)  ^{2}\right] \label{BadEventNonSmooth}\\
&  \leq\frac{1}{10}\exp\left[  -\left(  \log L\right)  ^{2}\right]  ,\nonumber
\end{align}
for all $L$ large enough (note that the choice
of $J_0$ and $L_0$ made above depended on the dimension only).
We now choose $\varepsilon_{0}
  >0$ small enough such that for $\varepsilon\leq\varepsilon
_{0},$ one has $X_{j}=0,$ deterministically, for $j<J_{0}
.$ 
%(As $\gamma$ is completely fixed in (\ref{Def_Gamma}), we
%usually don't explicitly indicate it in the notation).

We will show now
that if $\omega\notin\operatorname*{ManyBad}\nolimits_{L}$, then
$\left\Vert \Pi_{L}-\pi_{L}\right\Vert _{1}\leq\delta$, which
together with
(\ref{BadEventNonSmooth}),
will
prove
Proposition \ref{Prop_Nonsmooth}.
Toward this end, 
distinguish between two (disjoint) bad regions $B_{1},B_{2}\subset V_{L}.$
We set $\widetilde{B}_{L}\overset{\mathrm{def}}{=}B_{L}\backslash
\operatorname*{Sh}\nolimits_{L}  ,$ (
$B_{L}$ is as in (\ref{Def_BL})). Set%
\begin{equation}
	\label{eq-110206hh}
B_{2}^{\prime}\overset{\mathrm{def}}{=}\bigcup\left\{  D_{i}^{\left(
j\right)  }:\omega\in\operatorname*{Bad}\left(  D_{i}^{\left(  j\right)
}\right)  ,\ j=1,\ldots,J_1(L);\ i\leq N_{j}\right\}  .
\end{equation}
On the complement of $\operatorname*{TwoBad}\nolimits_{L}$ there exists
$x_{0}$ with $d_L(x_{0}) >r\left(  L\right)  ,$ such that
$\widetilde{B}_{L}\subset V_{5\rho\left(  x_{0}\right)  }\left(  x_{0}\right)
$. (See (\ref{BL_included_in_ball}). There is some ambiguity in choosing
$x_{0},$ but this of no importance. In particular,
$x_0$ is arbitrary if $\widetilde{B}_{L}=\emptyset$.) 
If $\left\vert x_{0}\right\vert \leq
L/2,$ we define $B_{1}\overset{\mathrm{def}}{=}V_{5\rho\left(  x_{0}\right)
}\left(  x_{0}\right)  =V_{5\gamma s\left(  L\right)  }\left(  x_{0}\right)
,$ and $B_{2}\overset{\mathrm{def}}{=}B_{2}^{\prime}.$ If $\left\vert
x_{0}\right\vert >L/2,$ we put $B_{1}\overset{\mathrm{def}}{=}\emptyset,$ and
$B_{2}\overset{\mathrm{def}}{=}B_{2}^{\prime}\cup V_{5\rho\left(
x_{0}\right)  }\left(  x_{0}\right)  .$ Of course, if $\widetilde{B}%
_{L}=\emptyset,$ then $B_{1}\overset{\mathrm{def}}{=}\emptyset,$ and
$B_{2}\overset{\mathrm{def}}{=}B_{2}^{\prime}$. Remark that $B_{1}$ and
$B_{2}$ are disjoint. We put $B\overset{\mathrm{def}}{=}B_{1}\cup B_{2},$ and
$G\overset{\mathrm{def}}{=}V_{L}\backslash B.$

In case $B_{1}=V_{5\gamma s\left(  L\right)  }\left(  x_{0}\right)  ,$
$\left\vert x_{0}\right\vert \leq L/2,$ we use the same (slight) modification
of $\hat{\Pi}\left(  y,\cdot\right)  ,\ \hat{\pi}\left(  y,\cdot\right)  $ for
$y\in V_{5\gamma s\left(  L\right)  }\left(  x_{0}\right)  $ as used in
Section \ref{Subsect_BadPoints}, i.e. we replace $\hat{\pi},\hat{\Pi}$ by
$\widetilde{\pi},\widetilde{\Pi}$ as defined in (\ref{Def_Enlargement}), but
we retain the \symbol{94}-notation for convenience.

We use a slight modification of the  perturbation expansion
(\ref{Pert1}). Again with $\Delta
\overset{\mathrm{def}}{=}\hat{\Pi}-\hat{\pi},$ we have%
\[
\Pi_{L}=\pi_{L}+\hat{g}1_{B}\Delta\Pi_{L}+\hat{g}1_{G}\Delta\Pi_{L}.
\]
Set $\gamma_{k}\overset{\mathrm{def}}{=}\hat{g}\left(  1_{G}\Delta\right)
^{k}.$ Then%
\begin{align*}
\gamma_{k}\Pi_{L}  &  =\hat{g}\left(  1_{G}\Delta\right)  ^{k}\Pi_{L}\\
&  =\hat{g}\left(  1_{G}\Delta\right)  ^{k}\pi_{L}+\hat{g}\left(  1_{G}%
\Delta\right)  ^{k}\hat{g}\Delta\Pi_{L}\\
&  =\hat{g}\left(  1_{G}\Delta\right)  ^{k}\pi_{L}+\hat{g}\left(  1_{G}%
\Delta\right)  ^{k}1_{B}\Delta\Pi_{L}+\hat{g}\left(  1_{G}\Delta\right)
^{k}\hat{\pi}\hat{g}\Delta\Pi_{L}+\gamma_{k+1}\Pi_{L}%
\end{align*}
Therefore, iterating, we get%
\begin{align*}
\Pi_{L}  &  =\pi_{L}+\hat{g}\sum_{k=0}^{\infty}\left(  1_{G}\Delta\right)
^{k}1_{B}\Delta\Pi_{L}+\hat{g}\sum_{k=1}^{\infty}\left(  1_{G}\Delta\right)
^{k}\hat{\pi}\hat{g}\Delta\Pi_{L}+\hat{g}\sum_{k=1}^{\infty}\left(
1_{G}\Delta\right)  ^{k}\pi_{L}\\
&  =\pi_{L}+\hat{g}\overline{\Gamma}1_{B}\Delta\Pi_{L}+\hat{g}\Gamma\hat{\pi
}\Pi_{L}.
\end{align*}
where $\Gamma\overset{\mathrm{def}}{=}\sum_{k=1}^{\infty}\left(  1_{G}%
\Delta\right)  ^{k},\ \overline{\Gamma}\overset{\mathrm{def}}{=}I+\Gamma.$
With the partition $B=B_{1}\cup B_{2},$ we get,
setting $\Xi
_{1}\overset{\mathrm{def}}{=}\hat{g}\overline{\Gamma}1_{B_{1}}\Delta,$
$\Xi_{2}\overset{\mathrm{def}}
{=}\hat{g}\overline{\Gamma}1_{B_{2}}\Delta$,%
\[
\Pi_{L}=\pi_{L}+\Xi_{1}\Pi_{L}+\Xi_{2}\Pi_{L}+\hat{g}\Gamma\hat{\pi}\Pi_{L},
\]
and by induction on $m\in\mathbb{N},$ replacing successively $\Pi_{L}$ in the
second summand%
\[
\Pi_{L}-\pi_{L}=\left(  \sum_{r=1}^{m}\Xi_{1}^{r}\right)  \pi_{L}+\left(
\sum_{r=0}^{m}\Xi_{1}^{r}\right)  \Xi_{2}\Pi_{L}+\left(  \sum_{r=0}^{m}\Xi
_{1}^{r}\right)  \hat{g}\Gamma\hat{\pi}\Pi_{L}+\Xi_{1}^{m+1}\Pi_{L}%
\]
i.e. with $m\rightarrow\infty$%
\begin{align}
\Pi_{L}-\pi_{L}  &  =\sum_{r=1}^{\infty}\Xi_{1}^{r}\pi_{L}+\left(  \sum
_{r=0}^{\infty}\Xi_{1}^{r}\right)  \Xi_{2}\Pi_{L}+\left(  \sum_{r=0}^{\infty
}\Xi_{1}^{r}\right)  \hat{g}\Gamma\hat{\pi}\Pi_{L}\label{PertExp_stopped}\\
& : =A_{1}+A_{2}+A_{3}\,.\nonumber
\end{align}
For $D\subset V_{L},$ we write%
\[
U_{k}\left(  D\right)  \overset{\mathrm{def}}{=}\left\{  y\in V_{L}:\exists
x\in D\ \mathrm{with\ }\Delta^{k}\left(  y,x\right)  >0\right\}  .
\]
We now prove that each of the 
three parts $A_{1},A_{2},A_{3}$ is bounded by
$\delta/3.$

\noindent\textbf{First summand} $A_{1}:$ This does not involve the bad regions
near the boundary, and we can apply the estimates from Section
\ref{Subsect_BadPoints}. There is nothing to prove if $B_{1}=\emptyset,$ so we
assume $B_{1}=V_{5\gamma s\left(  L\right)  }\left(  x_{0}\right)  ,$
$\left\vert x_{0}\right\vert \leq L/2.$ 
We have%
\begin{equation}
\sup_{x\in V_{L}}\left\vert \hat{g}\left( 
 1_{G}\Delta\right)  ^{k}\left(
x,B_{1}\right)  \right\vert \leq
\sup_{x\in V_{L}}
\delta^{k}\hat{g}\left(  x,U_{k}\left(
B_{1}\right)  \right)  \leq C\delta^{k}k^{d}, \label{NonSmooth1}%
\end{equation}
where the second inequality is due to part c) of
Lemma \ref{Cor_Green}, 
and therefore,%
\begin{equation}
\sum_{k=0}^{\infty}\left\Vert \hat{g}\left(  1_{G}\Delta\right)  ^{k}1_{B_{1}%
}\right\Vert _{1}\leq C. \label{NonSmooth2}%
\end{equation}
In the same way, we obtain, with
$K$ from Section \ref{Subsect_BadPoints},
\begin{equation}
\sum_{k=0}^{\infty}\sup_{x\notin V_{5K\gamma s\left(  L\right)  }}\left\Vert
\hat{g}\left(  1_{G}\Delta\right)  ^{k}1_{B_{1}}\left(  x,\cdot\right)
\right\Vert _{1}\leq\frac{1}{2}, \label{NonSmooth3}%
\end{equation}
by using (\ref{NonSmooth1}) for $k\geq1,$ and (\ref{eq-260805a}) for $k=0.$
 Furthermore\,,%
\begin{align}
\left\Vert \Xi_{1}\pi_{L}\right\Vert _{1}  &  \leq\sum_{k=0}^{\infty
}\left\Vert \hat{g}\left(  1_{G}\Delta\right)  ^{k}1_{B_{1}}\Delta\pi
_{L}\right\Vert _{1}\leq C\sum_{k=0}^{\infty}\left\Vert \hat{g}\left(
\cdot,U_{k}\left(  B_{1}\right)  \right)  \right\Vert _{\infty}2^{-k}%
\sup_{x\in B_{1}}\left\Vert \Delta\pi_{L}\left(  x,\cdot\right)  \right\Vert
_{1}\label{NonSmooth4}\\
&  \leq C\sum_{k=0}^{\infty}k^{d}2^{-k}\sup_{x\in B_{1}}\left\Vert \Delta
\pi_{L}\left(  x,\cdot\right)  \right\Vert _{1}\leq C\left(  \log L\right)
^{-3}.\nonumber
\end{align}
Using these inequalities, we get $\left\Vert A_{1}\right\Vert _{1}\leq
C\left(  \log L\right)  ^{-3}\leq C\left(  \log L_{0}\right)  ^{-3}\leq
\delta/3$ by choosing $L_{0}\left(  \delta\right)  $ large enough: When
estimating $\left\Vert \Xi_{1}^{r}\pi_{L}\right\Vert _{1}$ for $r\geq2,$ we
use (\ref{NonSmooth2}) for the first factor $\Xi_{1}$, (\ref{NonSmooth4}) for
the last $\Xi_{1}\pi_{L}$, and (\ref{NonSmooth3}) for the middle $\Xi
_{1}^{r-2}.$ The point is that $\left(  1_{B_{1}}\Delta\right)  \left(
x,y\right)  $ is $\neq0$ only if $y\notin V_{5K\gamma s\left(  L\right)
}\left(  x_{0}\right)  ,$ and so we can use (\ref{NonSmooth3}) for this part.

\noindent\textbf{Second summand }$A_{2}$: We drop here the $\Pi_{L}$-factor,
using the trivial estimate $\left\Vert \Pi_{L}\left(  x,\cdot\right)
\right\Vert _{1}\leq2$. If $r=0,$ one has to estimate $\left\Vert \Xi
_{2}\left(  0,\cdot\right)  \right\Vert _{1}$ where $B_{2}$ consists of the
bad regions in the layers $\mathcal{L}_{j}$, and the possible one bad ball
from $\widetilde{B}_{L}$ which is outside $V_{L/3}.$ In case $r\geq1,$ when
$B_{1}\neq\emptyset,$ we have $B_{2}=B_{2}^{\prime}$, which is at distance
$\geq L/3$ from $B_{1}.$ Therefore, in case $r=0,$ we have to estimate%
\begin{equation}
\left\Vert \hat{g}\left(  1_{G}\Delta\right)  ^{k}1_{B_{2}}\left(
0,\cdot\right)  \right\Vert _{1} \label{NonSmooth5}%
\end{equation}
(the last $\Delta$ is of no help, and we drop it), and in case $r\geq1,$ using
(\ref{NonSmooth2}) and (\ref{NonSmooth3})%
\[
C2^{-r}\sup_{\left\vert x\right\vert \leq2L/3}\left\Vert \hat{g}\left(
1_{G}\Delta\right)  ^{k}1_{B_{2}}\left(  x,\cdot\right)  \right\Vert _{1},
\]
but in this case, we have $B_{2}\subset\operatorname*{Shell}\nolimits_{L}%
\left(  2r\left(  L\right)  \right)  .$ The estimate of the second case is
entirely similar to the estimate of (\ref{NonSmooth5}), and we therefore 
provide the details only of the proof of the latter.

We split the parts coming from the different bad regions. For a bad
region $D_{i}^{\left(  j\right)  }$ in layer $\mathcal{L}_{j},$ we have%
\[
\left\Vert \hat{g}\left(  1_{G}\Delta\right)  ^{k}1_{D_{i}^{\left(  j\right)
}}\left(  0,\cdot\right)  \right\Vert _{1}\leq C2^{-k}\hat{g}\left(
0,U_{k}\left(  D_{i}^{\left(  j\right)  }\right)  \right)  .
\]
It suffices to estimate $\hat{g}\left(  0,U_{k}\left(  D_{i}^{\left(
j\right)  }\right)  \right)  $ very crudely. Points in $U_{k}\left(
D_{i}^{\left(  j\right)  }\right)  $ are at distance of
at  most $r_{j,k}=2^{j}\left(
1-2\gamma\right)  ^{-k}$ from $D_{i}^{\left(  j\right)  }.$ We first
consider $k$'s only such that 
$V_L\setminus \operatorname*{Shell}\nolimits_{L}\left(
s\left(  L\right)  \right)  $ is not touched, 
which is the case if
$k\leq20\log\log L$ ($L$ large enough). 
Then, for some $y$ with $0.5 r_{j,k}\leq d_L(y)\leq 2.5 r_{j,k}$,
 $U_{k}\left(
D_{i}^{\left(  j\right)  }\right)\subset B(y,r_{j,k})=:B_{j,k}$.
Applying Lemma 
\ref{Le_MainExit}, we see that 
the probability that 
simple random walk started at the origin hits $B_{j,k}$ before
$\tau_L$ is bounded above by 
$$C2^{(d-1)j}(1-2\gamma)^{-(d-1)k} L^{-d+1}\leq
C2^{(d-1)j} \left(\frac{3}{2}\right)^k L^{-d+1}\,,$$
where in the last inequality, we have used the definition
of $\gamma$, c.f. (\ref{Def_Gamma}).
Combined with Lemma \ref{Cor_Green} b), we conclude that
for any $r$ such that $\Lambda_r\cap
 U_{k}\left(
D_{i}^{\left(  j\right)  }\right)\neq \emptyset$,
it holds that
\[
\hat{g}\left(  0,U_{k}\left(  D_{i}^{\left(  j\right)  }\right) 
\cap \Lambda_r  \right)  \leq
C 2^{\left(  d-1\right)  j}\left(  \frac{3}{2}\right)
^{k}L^{-d+1}.
\]
The number of layers $r$ touched is
bounded by $2(1+k)$,
 and thus we conclude that
%for each $\Lambda_{r}$ which intersects $U_{k}\left(
%D_{i}^{\left(  j\right)  }\right)  $, a very crude estimate gives%
%\begin{equation}
%\label{eq-260805b}\left\vert U_{k}\left(  D_{i}^{\left(  j\right)  }\right)
%\cap\Lambda_{r}\right\vert \leq C2^{r}2^{\left(  d-1\right)  j}\left(
%1-2\gamma\right)  ^{-\left(  d-1\right)  k}\leq C2^{r}2^{\left(  d-1\right)
%j}\left(  \frac{3}{2}\right)  ^{k}%
%\end{equation}
%where in the last inequality, we have used (\ref{Def_Gamma}). 
%Using Lemma
%\ref{Le_MainExit}, we see that%
\[
\hat{g}\left(  0,U_{k}\left(  D_{i}^{\left(  j\right)  }\right)  \right)  \leq
C\left(  1+k\right)  2^{\left(  d-1\right)  j}\left(  \frac{3}{2}\right)
^{k}L^{-d+1}.
\]
Therefore, using $\omega\notin\bigcup\nolimits_{J_{0}
\leq j\leq J_1\left(  L\right)  }\left\{  X_{j}\geq j^{-3/2}N_{j}\right\}  ,$ we
have the estimates%
\begin{align*}
\sum_{k\leq10\log\log L}\left\Vert \hat{g}\left(  1_{G}\Delta\right)
^{k}1_{B_{2}^{\prime}\cap\Lambda_{j}}\left(  0,\cdot\right)  \right\Vert _{1}
&  \leq Cj^{-3/2},\\
\sum_{k\leq10\log\log L}\left\Vert \hat{g}\left(  1_{G}\Delta\right)
^{k}1_{B_{2}^{\prime}}\left(  0,\cdot\right)  \right\Vert _{1}  &  \leq
CJ_{0}^{-1/2}.
\end{align*}
For the sum over $k>20\log\log L,$ we simply estimate $\hat{g}\left(
0,U_{k}\left(  B_{2}^{\prime}\right)  \right)  \leq\hat{g}\left(
0,V_{L}\right)  \leq C\left(  \log L\right)  ^{6}$ and we therefore get%
\begin{align}
\sum_{k}\left\Vert \hat{g}\left(  1_{G}\Delta\right)  ^{k}1_{B_{2}%
\cap\operatorname*{Shell}\nolimits_{L}\left(  r\left(  L\right)  \right)
}\left(  0,\cdot\right)  \right\Vert _{1}  &  \leq C\left(  J_{0}%
^{-1/2}+\left(  \log L\right)  ^{6}2^{-20\log\log L}\right) \label{NonSmooth6}%
\\
&  \leq C\left(  J_{0}^{-1/2}+\left(  \log L\right)  ^{-7}\right)  \leq
\delta/6\nonumber
\end{align}
for all $L\geq L_0$,
by choosing $J_{0}=J_0(\delta)  $ and $L_{0} =L_0(\delta) $
large enough (again, depending only on $d$ and $\delta$).

It remains to add the part of $B_{2}$ outside $B_{2}^{\prime}.$ This is
(contained in) a ball $V_{5\gamma\rho\left(  x_{0}\right)  }\left(
x_{0}\right)  $ with $\left\vert x_{0}\right\vert >L/2.$
\[
\hat{g}\left(  0,U_{k}\left(  V_{5\gamma\rho\left(  x_{0}\right)  }\left(
x_{0}\right)  \right)  \right)  \leq\hat{g}\left(  0,U_{k}\left(  V_{5\gamma
s\left(  L\right)  }\left(  x_{0}\right)  \right)  \right)  \leq\hat{g}\left(
0,V_{\left(  5+2k\right)  \gamma s\left(  L\right)  }\left(  x_{0}\right)
\right)  .
\]
As $\left\vert x_{0}\right\vert \geq L/2,$ we have $V_{\left(  5+2k\right)
\gamma s\left(  L\right)  }\left(  x_{0}\right)  \cap V_{L/3}=\emptyset$
provided $k\leq\left(  \log L\right)  ^{3}/C,$ and $V_{\left(  5+2k\right)
\gamma s\left(  L\right)  }\left(  x_{0}\right)  $ can be covered by $\leq
Ck^{d}$ balls $V_{s\left(  L\right)  }\left(  y\right)  ,$ $\left\vert
y\right\vert \geq L/3.$ By Lemma \ref{Le_MainExit}, one has $\hat{g}\left(
0,V_{s\left(  L\right)  }\left(  y\right)  \right)  \leq C\left(  \log
L\right)  ^{-3}.$ (This remains true also if $V_{s\left(  L\right)  }\left(
y\right)  $ intersects $\operatorname*{Shell}\nolimits_{L}\left(  s\left(
L\right)  \right)  ,$ as is easily checked). Therefore, for $k\leq\left(  \log
L\right)  ^{3}/C,$ we have%
\[
\hat{g}\left(  0,U_{k}\left(  V_{5\gamma\rho\left(  x_{0}\right)  }\left(
x_{0}\right)  \right)  \right)  \leq Ck^{d}\left(  \log L\right)  ^{-3},
\]
and therefore,%
\begin{align*}
&  \sum_{k}\left\Vert \hat{g}\left(  1_{G}\Delta\right)  ^{k}1_{V_{5\gamma
\rho\left(  x_{0}\right)  }\left(  x_{0}\right)  }\left(  0,\cdot\right)
\right\Vert _{1}\\
\leq &  C\sum_{k\leq\left(  \log L\right)  ^{3}/C}2^{-k}k^{d}\left(  \log
L\right)  ^{-3}+C\sum_{k>\left(  \log L\right)  ^{3}/C}2^{-k}\left(  \log
L\right)  ^{6} \leq\delta/6,
\end{align*}
provided $L_{0}$ is large enough. Combining this with (\ref{NonSmooth6})
proves $\left\Vert A_{2}\right\Vert _{1}\leq\delta/3.$

\noindent\textbf{Third summand} $A_{3}.$ By the same argument as in the
discussion of $A_{2},$ it suffices to consider $r=0,$ and we 
drop $\Pi_{L}.$ Then, %
\begin{align}
&  \sum_{k\geq1}\left\Vert \sum_{x\notin\operatorname*{Shell}\nolimits_{L}%
\left(  r\left(  L\right)  \right)  }\hat{g}\left(  1_{G}\Delta\right)
^{k-1}\left(  0,x\right)  \left(  1_{G}\Delta\hat{\pi}\right)  \left(
x,\cdot\right)  \right\Vert _{1}\nonumber\\
&  \leq\sum_{k\geq1}^{\infty}2^{-k+1}\hat{g}\left(  0,V_{L}\right)
\sup_{x\notin\operatorname*{Shell}\nolimits_{L}\left(  r\left(  L\right)
\right)  }\left\Vert 1_{G}\Delta\hat{\pi}\left(  x,\cdot\right)  \right\Vert
_{1}\label{NonSmooth7}\\
&  \leq C\left(  \log L\right)  ^{-3}\leq\delta/9\nonumber
\end{align}
if $L_{0}$ is large enough.
 For $J_{0}  \leq j\leq
J_{1}\left(  L\right)  $,%
\begin{align*}
\left\Vert \sum_{x\in{\Lambda}_{j}}\hat{g}\left(  1_{G}\Delta\right)
^{k-1}\left(  0,x\right)  \left(  1_{G}\Delta\hat{\pi}\right)  \left(
x,\cdot\right)  \right\Vert _{1}  &  \leq2^{-k+1}\hat{g}\left(  0,U_{k}\left(
\Lambda_{j}\right)  \right)  \sup_{x\in\Lambda_{j}}\left\Vert 1_{G}\Delta
\hat{\pi}\left(  x,\cdot\right)  \right\Vert _{1}\\
&  \leq Cj^{-9}2^{-k+1}\hat{g}\left(  0,U_{k}\left(  \Lambda_{j}\right)
\right)  ,
\end{align*}
and it is evident from  part b) of
Lemma \ref{Cor_Green} 
%Lemma \ref{Le_MainExit} 
that $\sum_{k\geq1}2^{-k+1}%
\hat{g}\left(  0,U_{k}\left(  \Lambda_{j}\right)  \right)  \leq C.$ Therefore,%
\begin{equation}
\sum_{k\geq1}\left\Vert \sum_{J_{0}  \leq j\leq
J_{1}\left(  L\right)  }\sum_{x\in{\Lambda}_{j}}\hat{g}\left(  1_{G}%
\Delta\right)  ^{k-1}\left(  0,x\right)
  \left(  1_{G}\Delta\hat{\pi}\right)
\left(  x,\cdot\right)  \right\Vert _{1}\leq C\left(
J_{0}
\right)^{-8}\leq\delta/9, \label{NonSmooth8}%
\end{equation}
if $J_{0}$ is chosen large enough
(again, independently of $\varepsilon_0$!). 
On the other hand, putting
$\hat{\Lambda}\overset{\mathrm{def}}%
{=}\bigcup\nolimits_{j\leq J_{0}  }\Lambda_{j}$,
\begin{align*}
&\sum_{k\geq1}\left\Vert \sum\nolimits_{x\in\hat{\Lambda}}\hat{g}\left(
1_{G}\Delta\right)^{k-1}\left(  0,x\right)  \left(  1_{G}\Delta\hat{\pi
}\right)  \left(  x,\cdot\right)  \right\Vert _{1}\\
&  \leq C\sum_{k\geq
1}2^{-k+1}\hat{g}\left(  0,U_{k}\left(  \hat{\Lambda}\right)  \right)
\sup_{x\in\hat{\Lambda}}\left\Vert \Delta\left(x,\cdot\right)
\right\Vert_{1}
  \leq C\left(  J_{0}\right)  \sup_{x\in\hat{\Lambda}}\left\Vert
\Delta\left(  x,\cdot\right)  \right\Vert _{1}\leq\delta/9
\end{align*}
if $\varepsilon\leq\varepsilon_{0}\left(  \delta\right)$ 
and $\varepsilon_0(\delta)$ is taken small enough.
Combining this
with (\ref{NonSmooth7}) and (\ref{NonSmooth8}) proves $\left\Vert
A_{3}\right\Vert _{1}\leq\delta/3.$
Substituting the estimate on $\|A_i\|_1$, $i=1,2,3$,
into (\ref{PertExp_stopped})
and using (\ref{BadEventNonSmooth})
completes the proof of
Proposition \ref{Prop_Nonsmooth}.
\end{proof}

\section{Proof of Proposition \ref{Prop_Main}}
\label{sec-finalpush}
We just have to collect the estimates we have obtained so far. We take
$\delta_{0}$ small enough as in the
conclusion of
Propositions \ref{prop-erwin220705} and
\ref{Prop_Nonsmooth}, and for $\delta\leq\delta_{0},$ we choose
$L_{0}$ large enough, also according to these propositions.

For $L_{1}\geq L_{0}$ we assume $\operatorname*{Cond}\left(  \delta
,L_{1}\right)  ,$ and take and $L\leq L_{1}\left(  \log L_{1}\right)  ^{2}.$
For $i=1,2,3,$ and $\Psi\in\mathcal{M}_{L},$ we have according to Lemma
\ref{Le_TwoBad} and Proposition \ref{Prop_Main_no_bad}%
\begin{align*}
b_{i}\left(  L,{\Psi},\delta\right)   &  \leq\mathbb{P}\left(  D_{L,\Psi
}\left(  0\right)  >\left(  \log L\right)  ^{-11.25-2.25i}\right)  \\
&  \leq\mathbb{P}\left(  D_{L,\Psi}\left(  0\right)  >\left(  \log L\right)
^{-11.25-2.25i},\left(  \operatorname*{TwoBad}\nolimits_{L}\right)  ^{c}%
\cap\left(  \operatorname*{Good}\nolimits_{L}\right)  ^{c}\right)  \\
&  +\mathbb{P}\left(  D_{L,\Psi}\left(  0\right)  >\left(  \log L\right)
^{-9},\operatorname*{TwoBad}\nolimits_{L}\cap\operatorname*{Good}%
\nolimits_{L}\right)  +\mathbb{P}\left(  \operatorname*{TwoBad}\nolimits_{L}%
\right)  \\
&  \leq\mathbb{P}\left(  D_{L,\Psi}\left(  0\right)  >\left(  \log L\right)
^{-11.25-2.25i},\left(  \operatorname*{TwoBad}\nolimits_{L}\right)  ^{c}%
\cap\left(  \operatorname*{Good}\nolimits_{L}\right)  ^{c}\right)  \\
&  +\exp\left[  -1.2\left(  \log L\right)  ^{2}\right]  +\exp\left[  -\left(
\log L\right)  ^{17/8}\right]  .
\end{align*}
We therefore only have to estimate the first summand.
Using the notation of Sections \ref{Subsect_greengood}
and \ref{Subsect_BadPoints},%
\begin{align*}
& \mathbb{P}\left(  D_{L,\Psi}\left(  0\right)  >\left(  \log L\right)
^{-11.25-2.25i},\left(  \operatorname*{TwoBad}\nolimits_{L}\right)  ^{c}%
\cap\left(  \operatorname*{Good}\nolimits_{L}\right)  ^{c}\right)  \\
& \leq\sum_{D\in\mathcal{D}_{L}}\sum_{j}\mathbb{P}\left(  \left\Vert \left(
\left[  \Pi_{V_{L}}-\pi_{V_{L}}\right]  \hat{\pi}_{{\Psi}}\right)  \left(
0,\cdot\right)  \right\Vert _{1}\geq\left(  \log L\right)  ^{-11.25+2.25i}%
,\ \operatorname*{Bad}\nolimits_{L}^{\left(  j\right)  }\left(  D\right)
\right)  \\
& \leq\sum_{D\in\mathcal{D}_{L}}\sum_{j\leq i}\mathbb{P}\left(  \left\Vert
\left(  \left[  \Pi_{V_{L}}-\pi_{V_{L}}\right]  \hat{\pi}_{{\Psi}}\right)
\left(  0,\cdot\right)  \right\Vert _{1}\geq\left(  \log L\right)
^{-11.25+2.25j},\ \operatorname*{Bad}\nolimits_{L}^{\left(  j\right)  }\left(
D\right)  \right)  \\
& \quad\quad +\sum_{D\in\mathcal{D}_{L}}\sum_{j>i}\mathbb{P}\left(
\operatorname*{Bad}\nolimits_{L}^{\left(  j\right)  }\left(  D\right)
\right)  \\
& \leq\frac{4\left\vert \mathcal{D}_{L}\right\vert }{100}\exp\left[  -\left(
\log L\right)  ^{2}\right]  +\left\vert \mathcal{D}_{L}\right\vert \exp\left[
-\left[  1-\left(  4-i-1\right)  /13\right]  \left(  \log\frac{L}{\left(  \log
L\right)  ^{10}}\right)  ^{2}\right]  \\
& \leq\frac{1}{8}\exp\left[  -\left[  1-\left(  4-i\right)  /13\right]
\left(  \log L\right)  ^{2}\right]  .
\end{align*}
Combining these estimates, we get for $i=1,2,3$ and $L$ large enough,%
\[
b_{i}\left(  L,{\Psi},\delta\right)  \leq\frac{1}{4}\exp\left[  -\left[
1-\left(  4-i\right)  /13\right]  \left(  \log L\right)  ^{2}\right]\,.
\]

For $i=4,$ we have%
\[
b_{4}\left(  L,{\Psi},\delta\right)  \leq\mathbb{P}\left(  D_{L,\Psi}\left(
0\right)  >\left(  \log L\right)  ^{-2.25}\right)  +\mathbb{P}\left(
\left\Vert \Pi_{L}\left(  0,\cdot\right)  -\pi_{L}\left(  0,\cdot\right)
\right\Vert _{1}\geq\delta\right)  .
\]
The second summand is estimated by Proposition \ref{Prop_Nonsmooth}, and the
first in the same way as the $b_{i},\ i\leq3.$
This completes the proof of Proposition \ref{Prop_Main}.

\section{Proof of Theorem \ref{Th_Main1}\label{proof-Main1}}
The proof is based on a modification of 
the computations in
Section \ref{Sect_NonSmooth}. 
We begin by an auxiliary  definition.
%Fix some $\delta\in(0,\delta_0)$. 
In what follows, $c_d$ is a constant large enough so as to satisfy 
\begin{equation}
	\label{eq-110206a} \log c_d>4d\,.
\end{equation}
%Recall the notation $\Lambda_j$, c.f. (\ref{eq-110206b}). 
%For $x\in \Lambda_j$, set $U_1(x)=\partial V_{\gamma d_L(x)}(x)$, and 
%inductively, for $k\geq 2$, set $U_k(x)=\cup_{y\in U_{k-1}(x)} \partial
%V_{c_d^{k} 2^j}(y)$. 
%For a nearest neighbor path $X_n$
%started at
%$y\in U_1(x)$, 
%let $\eta_2=\min\{n: X_n\in \partial V_{c_d^2 2^J}(y)\}$ and, for
%$k=2,\ldots,K$ and $z\in V_L$, let $\eta_k(z)=\min\{n: X_n\in \partial
%V_{c_d^k 2^J}(z)\}$ and define succesively $\eta_k=\eta_k(X_{\eta_{k-1}})$.
For any $x\in V_L$ and random walk $\{X_n\}$ with $X_0=x$, set
$\eta(x)=\min\{n>0: |X_n-x|>d_L(x)\}$.
We fix $\delta$ small enough such that $\delta<\delta_0$ and
\begin{equation}
	\label{eq-110206z}
\bar c_d	
\overset{\mathrm{def}}{=}\max_{y\in V_L} 
	p^{\mathrm{RW},y}(\tau_L>\eta(y))+\delta<1\end{equation}
for all $L$ large enough (this is possible by the
	invariance principle for simple random walk). 
	%We let $\bar c_d$ 
	%denote the value of the maximum in (\ref{eq-110206z}).
We then chose $\varepsilon_0$ so as the conclusion of Theorem
\ref{Th_Main} is satisfied with this value of $\delta$.
For $x\in \Lambda_j$, set $U_1(x)=\partial V_{\gamma d_L(x)}(x)$, and 
inductively, for $k\geq 2$, set $U_k(x)=\cup_{y\in U_{k-1}(x)} \partial
V_{c_d^{k} 2^j}(y)$. 
\begin{definition}
	\label{def-110206a}
A point $x\in \Lambda_j$
is $K$-good if $D_{\gamma d_L(x),0}(x)\leq \delta$ and, 
for any $k=1,\ldots,K$,  
for all $y\in U_k(x)$, $D_{c_d^{k+1} 2^j,0}(y)\leq \delta$.
\end{definition}
We note that  there exists
some constant $C$
depending only  on $d$ such that 
for all $x\in \Lambda_j$, 
\begin{equation}
	\label{eq-110206gg}
	\mathbb{P}(x\,
	\mbox{\rm is not $K$-good})\leq 
	\exp\left(-\frac{10}{13} (\log (\gamma  2^{j-1}))^2\right)+
	\sum_{k=1}^K C(c_d^k 2^j)^{d}
\exp\left(-\frac{10}{13} (\log ( c_d^k 2^j))^2\right)\,.
\end{equation}
For $J>J_0$, set
$$\tilde B_{L,J,K}^{(4)}
\overset{\mathrm{def}}{=}
\left(\cup_{j=J+1}^{J_1(L)} \{x\in \Lambda_j:
D_{\gamma d_L(x),0}(x)
\geq \delta\}\right)\bigcup\left(\{x\in \Lambda_J: \mbox{\rm $x$ is
not $K$-good}\}\right)\,,$$
and 
$\widetilde{ 
\operatorname*{Good}}_{L,J,K}
\overset{\mathrm{def}}{=}
V_L\setminus \tilde B_{L,J,K}^{\left(  4\right)  }$.

If $B\in\mathcal{L}_{j},$ $j\geq J$,
we write $\widetilde{\operatorname*{Bad}}_{L,J,K}\left(  B\right)  $ for
the event $\left\{  B\not \subset\widetilde 
{\operatorname*{Good}}_{L,J,K}\right\}.$ Remark
that for $J>J_1(L)$, 
$\widetilde{\operatorname*{Bad}}_{L,J,K}\left(  B\right)  =
 \operatorname*{Bad}\left(  B\right) $. 
 By combining the computation in (\ref{eq-110206c}) 
 with (\ref{eq-110206gg}),
 and using our choice for the value of $c_d$,
 we get that
there exists a $J_2\geq J_0$ such that for all $J>J_2$, all $K$,
and all $L$ large enough,
\begin{equation}
	\label{eq-110206d}
	\mathbb{P}\left(  \widetilde{\operatorname*{Bad}}_{L,J,K}
	\left(  B\right)  \right)  \leq p_J\,.
\end{equation}

We set next
\[
\tilde
X_{j}\overset{\mathrm{def}}{=}\sum_{D\in\mathcal{L}_{j}}
1_{\widetilde{\operatorname*{Bad}}_{L,J,K}\left(  D\right)  }\]
and
\[
\widetilde{\operatorname*{ManyBad}}_{L,J,K}
\overset{\mathrm{def}}{=}\bigcup
\nolimits_{J_{0}\left(  \gamma\right)  \leq j\leq J\left(  L\right)  }\left\{
\tilde X_{j}\geq j^{-3/2}N_{j}\right\}  
\cup\operatorname*{TwoBad}\nolimits_{L}\,.
\]
Arguing as in the computation leading to 
(\ref{BadEventNonSmooth}) (except that ${\cal L}_J$ is divided
into more sets to achieve independence, and the number of such sets depends
on $K$), we conclude that for each $J,K$ there is an $L_2=L_2(J,K)$ such that
for all $L>L_2(J,K)$,
\begin{equation}
	\mathbb{P}\left(  
	\widetilde{\operatorname*{ManyBad}}_{L,J,K}\right)     \leq
\frac{1}{10}\exp\left[  -\left(  \log L\right)  ^{2}\right]\,. 
 \label{BadEventNonSmoothtilde}
\end{equation}
Next, replace $B_2'$ in 
(\ref{eq-110206hh}) by considering  
in the union there only $j\in [J,J_1(L)]$, and replacing
$\operatorname*{Bad}\left(  D_{i}^{\left(  j\right)
}\right)$ by 
$\widetilde{\operatorname*{Bad}}_{J,K,L}\left(  D_{i}^{\left(  j\right)
}\right)$ (note that this influences only the layers $\Lambda_j$ with $j\leq J$).
We rewrite the expansion (\ref{PertExp_stopped})
\begin{align}
\Pi_{L}-\pi_{L}  &  =\sum_{r=1}^{\infty}\Xi_{1}^{r}\pi_{L}+\left(  \sum
_{r=0}^{\infty}\Xi_{1}^{r}\right)  \Xi_{2}\Pi_{L}+
\left(  \sum_{r=0}^{\infty
}\Xi_{1}^{r}\right)  \hat{g}\Gamma\hat{\pi}\Pi_{L}\label{PertExp_stopped1}\\
& : =A_{1}+A_{2}+A_{3}\,.\nonumber
\end{align}
except that now $\tilde B_2'$ is used instead of $B_2'$ in the definition
of $\Xi_{2}$. By repeating the computation leading to  
(\ref{NonSmooth4}) and 
(\ref{NonSmooth6}) (using 
 (\ref{BadEventNonSmoothtilde}) instead of  
 (\ref{BadEventNonSmooth}) in the latter), we conclude that
 \begin{equation}
	 \label{eq-110206w}
	\|A_1\|+\|A_2\|\leq CJ^{-1/2}\,,
\end{equation}
for $L$ large. To analyze $A_3$, 
 we write, with obvious notation, $A_3=\sum_{r=0}^\infty A_3^{(r)}$, and argue
 as in Section \ref{Sect_NonSmooth} that it is enough to consider
 $A_3^{(0)}$. We then write
 %it is again
 %enough to focus on the term  $r=0$, and we write
 %we split it as
 \begin{eqnarray}
	 \label{eq-110206kk1}
	 A_3^{(0)}(\cdot)&=&
 \sum_{k\geq1} \sum_{x\notin\operatorname*{Shell}\nolimits_{L}
\left(  2^J  \right)  }\hat{g}\left(  1_{G}\Delta\right)
^{k-1}\left(  0,x\right)  \left(  1_{G}\Delta\hat{\pi}\right)  \left(
x,z\right) \Pi_L (z,\cdot)\\
&&
\!\!\!\!\!\!\!\!
+
 \sum_{k\geq1} \sum_{x\in\operatorname*{Shell}\nolimits_{L}
\left(  2^J  \right)  }\hat{g}\left(  1_{G}\Delta\right)
^{k-1}\left(  0,x\right)  \left(  1_{G}\Delta\hat{\pi}\right)  \left(
x,z\right) \Pi_L (z,\cdot)=A_3^{(0),1}+A_3^{(0),2}\,.
\nonumber
\end{eqnarray}
We already have from 
Section \ref{Sect_NonSmooth} that on the event
	$\left(\widetilde{\operatorname*{ManyBad}}_{L,J,K}\right)^c$, 
	it holds that
	$\Vert A_3^{(0),1}\Vert_1\leq C J^{-1/2}$. Note that all the estimates
	so far held for any large fixed $J,K$, as long as $L$ is large enough
	(large enough depending on the choice of $J,K$).
	It remains only
	to analyze $\|A_3^{(0,1)} \hat \pi_s\Vert_1$ for $s$ independent
	of $L$, and when doing that, we can 
	choose $K$ in any way that does not depend on $L$. As a preliminary step
	in the choice of  $K$, we have the following lemma
	(recall the constant $\bar c_d$, c.f. (\ref{eq-110206z})):
\begin{lemma}
	\label{lem-021106a}
	For all $J,K$, and 
	any  $K$-good
$x\in \Lambda_J$, it holds that
for all $L$ large enough,
\begin{equation}
	\label{eq-110206ll}
	\sum_{y: |x-y|>c_d^{K+2} 2^J} | \Delta \hat \pi \Pi_L(x,y)| \leq 
	(\bar c_d)^ K
\end{equation}
\end{lemma}
\begin{proof}
	[Proof of Lemma \ref{lem-021106a}]
The proof is similar to the argument in Lemma 
\ref{lem-130705}. Consider a RWRE $X_n$ started at
$y\in U_1(x)$. Let $\eta_1=\min\{n: X_n\in \partial V_{c_d 2^J}(y)\}$,
and, for
$k=2,\ldots,K$, define successively
$\eta_k(y)=\min\{n>\eta_{k-1}(y): X_n\in \partial
V_{c_d^k 2^J}(
X_{\eta_{k-1}})
)\}$.
%and define succesively $\eta_k=\eta_k(X_{\eta_{k-1}})$.
%and let $\eta_K$ 
%denote the hitting succesive hitting times of $ U_k(x)$, $k\geq 2$.
The sum in (\ref{eq-110206ll}) is then bounded above by 
\begin{equation}
	\label{eq-110206zz}
	\max_{y\in U_1(x)} 
p_\omega^y(\tau_L>\eta_K(y))\leq \prod_{k=1}^K \max_{y\in U_k(x)} 
p_\omega^y(\tau_L>\eta(y))\,.
\end{equation}
Since $x$ is $K$-good, we have that for $y\in U_k(x)$,
$p_\omega^y(\tau_L>\eta_k(y))\leq \delta+p^{\mathrm{RW},y}(\tau_L>\eta(y))\leq
\bar c_d$. Substituting in (\ref{eq-110206zz}), the lemma follows.
\end{proof}

We next recall, c.f. Lemma \ref{Le_Lipshitz_Pi}, that
$$\sup_{x_1,x_2: |x_1-x_2|\leq D}\|\hat \pi_s(x_1,\cdot)-\hat \pi_s(x_2,\cdot)\|
\leq  \frac{CD \log s}{s}\,.$$
In particular, for any fixed $K,J$ with $J>J_2$,
	using (\ref{eq-110206ll}), and the fact that
	$\sum_z A_3^{(0),1}\hat\pi_s(z)=0$, we get
$$	\Vert A_3^{(0),1} \hat\pi_s\Vert_1\leq 
\delta (\bar c_d)^K+
	C\frac{c_d^{K+1} \log s}{s}\,.$$
and thus
\begin{equation}
	\label{eq-110206y}
	\Vert A_3 \hat \pi_s\Vert 
\leq \delta (\bar c_d)^K+
C\frac{(c_d)^{K+1}2^J \log s}{s}+CJ^{-1/2}\,.
\end{equation}
Combining (\ref{eq-110206y}) and 
	 (\ref{eq-110206w}),
	 we conclude that on the event 
	$\left(\widetilde{\operatorname*{ManyBad}}_{L,J,K}\right)^c$,
	 $$D_{L,\Psi_s}(0)
\leq \delta (\bar c_d)^K+
C\frac{2(c_d)^{K+1}2^J \log s}{s}+CJ^{-1/2}\,.$$
Choosing $J$ large such that $CJ^{-1/2}<\delta/3$
and $K$ large enough such that $(\bar c_d)^K<\delta/3$, and
$s$ large enough such that $2C (c_d)^{K+1}2^J\log s/s<\delta/3$,
and using  (\ref{BadEventNonSmoothtilde}),
it follows that\\
$\limsup_{L\to\infty} L^r b(L,\Psi_s,\delta)$ $=0$\,.
This completes the proof of the theorem. 

\appendix{}

\section{Proofs of the random walk results \label{App_A}}

%\begin{proof}
%[Proof of Lemma \ref{Le_MainExit}]
We begin by stating and proving some auxiliary estimates. If $A\subset
\subset\mathbb{Z}^{d},$ $x\in A,\ y\in\partial A,$ then by time reversibility
of simple random walk and transience, 
%the usual time
%reversal, 
one gets%
\begin{equation}
\label{Est_NonLog2}
P_{x}\left(  X_{\tau_{A}}=y\right)   
%&  =\sum_{y^{\prime}\in A,\ \left\vert
%y-y^{\prime}\right\vert =1}\left(  2d\right)  ^{-1}g_{A}\left(  y^{\prime
%},x\right) \label{Est_NonLog2}\\
  \leq C \sum_{y^{\prime}\in A,\ \left\vert y-y^{\prime}\right\vert =1}
  %\left(
%2d\right)  ^{-1}
P_{y^{\prime}}\left(  T_{x}<\tau_{A}\right)  .\nonumber
\end{equation}
Throughout this appendix, we write $\tau\overset{\mathrm{def}}{=}\tau_{V_{L}}$.
Since we do not deal with the RWRE in this appendix, we
consistently write $P_x$ for $P^{\mathrm{RW}}_x$.
\begin{lemma}
\label{Le_NonLog2}Let $x\in V_{L},\ y\in\partial V_{L}.$ Then,
%(with
%$\tau\overset{\mathrm{def}}{=}\tau_{V_{L}}$)
for some $\bar c_{1}\geq1$,%
\[
P_{x}\left(  X_{\tau}=y\right)  \leq\bar c_{1}d_{L}\left(  x\right)  ^{-d+1}.
\]

\end{lemma}
\begin{proof}
Let $r\overset{\mathrm{def}}{=}d_{L}\left(  x\right)  .$ We may assume that
$r\geq4.$ Put $r^{\prime}\overset{\mathrm{def}}{=}\left[  r/2\right]  -1.$
Then $V_{r^{\prime}}\left(  x\right)  \subset V_{L-r^{\prime}}.$ If
$y^{\prime}$ is any neighbor of $y$ in $V_{L}$ then, by part c) of 
Lemma \ref{Le_Lawler_Exit},%
\[
P_{y^{\prime}}\left(  T_{\partial V_{r^{\prime}}\left(  x\right)  }<\tau
\right)  \leq P_{y^{\prime}}\left(  T_{V_{L-r^{\prime}}}<\tau\right)
\leq\frac{C}{r}.
\]
Furthermore uniformly in $z\in\partial V_{r^{\prime}}\left(  x\right)  , $%
\[
P_{z}\left(  T_{x}<\tau\right)  \leq P_{z}\left(  T_{x}<\infty\right)
\leq C(r^{\prime})^{-d+2}\leq Cr^{-d+2}.
\]
Using the Markov property and (\ref{Est_NonLog2}) proves the claim.
\end{proof}

\begin{lemma}
\label{Le_NonLog3}Let $x\in V_{L},\ y\in \partial V_{L}$ and set $t\overset{%
\mathrm{def}}{=}\left\vert x-y\right\vert .$ Then for some $\bar{c}_{2}\geq 1
$,%
\begin{equation*}
P_{x}\left( X_{\tau }=y\right) \leq \bar{c}_{2}\frac{d_{L}\left( x\right) }{t%
}\sup_{x^{\prime }\in \partial V_{t/3}\left( y\right) \cap
V_{L}}P_{x^{\prime }}\left( X_{\tau }=y\right) .
\end{equation*}
\end{lemma}
\begin{proof}
The bound is evident if $r\overset{\mathrm{def}}{=}d_{L}\left( x\right) \geq
t/10.$ Therefore, we assume $r<t/10.$ We choose a point $x^{\prime }\notin
V_{L}$ such that $V_{t/8}\left( x^{\prime }\right) \cap V_{L}=\emptyset ,$
and $\left\vert x-x^{\prime }\right\vert \leq t/8+2r.$ Then $x\in
V_{t/4}\left( x^{\prime }\right) $ and $\left\vert x^{\prime }-y\right\vert
\geq 3t/4.$ Therefore%
\begin{equation}
\partial V_{t/4}\left( x^{\prime }\right) \cap V_{t/3}\left( y\right)
=\emptyset .  \label{Disjoint_balls}
\end{equation}
A walk starting in $x$ has to reach $\partial V_{t/4}\left( x^{\prime
}\right) $ before it can reach $y,$ and therefore, if it reaches $%
V_{t/8}\left( x^{\prime }\right) $ before reaching $\partial V_{t/4}\left(
x^{\prime }\right) $ it exits $V_{L}$ before it reaches $y.$ After it has
reached $\partial V_{t/4}\left( x^{\prime }\right) \cap V_{L}$ it still has
to reach $V_{t/3}\left( y\right) \cap V_{L}$ before it can reach $y,$ moving
inside $V_{L}.$ So we get%
\begin{eqnarray*}
P_{x}\left( X_{\tau }=y\right)  &\leq &CP_{x}\left( \tau _{V_{t/4}\left(
x^{\prime }\right) }<T_{V_{t/8}\left( x^{\prime }\right) }\right) \sup_{z\in
\partial V_{t/4}\left( x^{\prime }\right) \cap V_{L}}P_{z}\left( X_{\tau
}=y\right)  \\
&\leq &\frac{Cr}{t}\sup_{z\in \partial V_{t/4}\left( x^{\prime }\right) \cap
V_{L}}P_{z}\left( X_{\tau }=y\right) \leq \frac{Cr}{t}\sup_{z\in \partial
V_{t/3}\left( y\right) \cap V_{L}}P_{z}\left( X_{\tau }=y\right) ,
\end{eqnarray*}%
the second inequality using Lemma \ref{Le_Lawler_Exit} c).
\end{proof}

\begin{lemma}
\label{Le_NonLog4}With $x,y,t$ as above, and $\bar c_{1},\bar c_{2}$ from the
previous lemmas,%
\[
P_{x}\left(  X_{\tau}=y\right)  \leq\bar c_{1} \bar c_{2}^{d}3^{\left(
d-1\right)  ^{2}+(d-1)}\frac{d_{L}\left(  x\right)  }{t^{d}}.
\]
\end{lemma}
\begin{proof}
Put $\eta\overset{\mathrm{def}}{=}3^{-d+1}\bar c_{2}^{-1},$ and set $\bar
K\overset{\mathrm{def}}{=}\bar c_{1}\eta^{-d+1}=\bar c_{1}\bar c_{2}%
^{d-1}3^{\left(  d-1\right)  ^{2}}.$ Using Lemma \ref{Le_NonLog3}, it suffices
to prove
\begin{equation}
\sup_{x\in\partial V_{r}(y)\cap V_{L}}P_{x}
\left(  X_{\tau}=y\right)  \leq\bar
Kr^{-d+1}. \label{Est_NonLog1}%
\end{equation}
As $\bar K\geq 9^{\left(  d-1\right)  }$, there is nothing to prove if
$r\leq9.$ Assume that we have proved (\ref{Est_NonLog1}) for all $r\leq
r_{0},$ and assume $r_{0}<r\leq2r_{0}.$ Then for $d_{L}\left(  x\right)  >\eta
r,$ we have by Lemma \ref{Le_NonLog2} that%
\[
P_{x}\left(  X_{\tau}=y\right)  \leq \bar
c_{1}\eta^{-d+1}r^{-d+1}=\bar Kr^{-d+1},
\]
and for $d_{L}\left(  x\right)  \leq\eta r,$ by Lemma \ref{Le_NonLog3} and the
fact that $r/3\leq r_{0}$,%
\[
P_{x}\left(  X_{\tau}=y\right)  \leq \bar c_{2}\eta\bar K\left(  \frac{r}%
{3}\right)  ^{-d+1}=\bar Kr^{-d+1}.
\]
Therefore, the lemma is proved by induction.
\end{proof}

\noindent
\begin{proof}
[Proof of Lemma \ref{Le_MainExit}]If $\left\vert x-y\right\vert \leq
d_{L}\left(  y\right)  /2,$ then $d_{L}\left(  x\right)  \geq d_{L}\left(
y\right)  /2,$ and in this case, we can simply use
part c) of Lemma \ref{Le_Lawler_Exit} to conclude that%
\[
P_{x}\left(  T_{V_{a}\left(  y\right)  }<\tau\right)  \leq P_{x}\left(
T_{V_{a}\left(  y\right)  }<\infty\right)  \leq C\left(  \frac{a}{\left\vert
x-y\right\vert }\right)  ^{d-2}\leq C\frac{a^{d-2}d_{L}\left(  y\right)
d_{L}\left(  x\right)  }{\left\vert x-y\right\vert ^{d}}.
\]
Therefore, we may assume $\left\vert x-y\right\vert >d_{L}\left(  y\right)
/2.$ Furthermore, it suffices to consider the case $1\leq a\leq d_{L}\left(
y\right)  /5,$ simply because for $d_{L}\left(  y\right)  /5<a\leq
5d_{L}\left(  y\right)  $, we get an upper bound with replacing $a$ by
$d_{L}\left(  y\right)  /5.$ Assume that we have proved the bound for
$a=d_{L}\left(  y\right)  /5.$ Then we get for $a<d_{L}\left(  y\right)  /5$%
\[
P_{x}\left(  T_{V_{a}\left(  y\right)  }<\tau\right)  \leq C\frac{d_{L}\left(
y\right)  ^{d-1}d_{L}\left(  x\right)  }{\left\vert x-y\right\vert ^{d}%
}\left(  \frac{a}{d_{L}\left(  y\right)  }\right)  ^{d-2}\leq C\frac
{a^{d-2}d_{L}\left(  y\right)  d_{L}\left(  x\right)  }{\left\vert
x-y\right\vert ^{d}}.
\]
We therefore see that it suffices to prove the bound for $a=d_{L}\left(
y\right)  /5.$

Let $y^{\prime}\in\partial V_{L}$ be a point closest to $y.$ There exists
$\delta>0,$ such that%
\[
\inf_{x^{\prime}\in V_{a}\left(  y\right)  }P_{x^{\prime}}\left(  X_{\tau}\in
V_{a}\left(  y^{\prime}\right)  \right)  \geq\delta.
\]
Evidently, $\inf_{z\in V_{a}\left(  y^{\prime}\right)  \cap\partial V_{L}%
}\left\vert x-z\right\vert \geq\left\vert x-y\right\vert /2,$ and therefore,
by Lemma \ref{Le_NonLog4},%
\[
\sup_{z\in V_{a}\left(  y^{\prime}\right)  \cap\partial V_{L}}P_{x}\left(
X_{\tau}=z\right)  \leq C\frac{d_{L}\left(  x\right)  }{\left\vert
x-y\right\vert ^{d}}.
\]
Consequently%
\begin{align*}
\frac{d_{L}\left(  x\right)  a^{d-1}}{\left\vert x-y\right\vert ^{d}}  &
\geq\frac{1}{C}P_{x}\left(  X_{\tau}\in V_{a}\left(  y^{\prime}\right)
\right)  \geq\frac{1}{C}P_{x}\left(  X_{\tau}\in V_{a}\left(  y^{\prime
}\right)  ,\ T_{V_{a}\left(  y\right)  }<\tau\right) \\
&  =\frac{1}{C}\sum_{x^{\prime}\in V_{a}\left(  y\right)  }P_{x}\left(
T_{V_{a}\left(  y\right)  }<\tau,\ X_{T_{V_{a}\left(  y\right)  }}=x^{\prime
}\right)  P_{x^{\prime}}\left(  X_{\tau}\in V_{a}\left(  y^{\prime}\right)
\right) \\
&  \geq\frac{\delta}{C}P_{x}\left(  T_{V_{a}\left(  y\right)  }%
<\tau\right)  .
\end{align*}
This proves the claim.
\end{proof}
Before presenting the proofs of Lemmas \ref{Le_Approx_Phi_by_BM} and
\ref{Le_ThirdDerivative}, we 
%introduce some notation and 
state and prove some
additional auxiliary estimates. 
We define the Brownian analogue $\hat\pi_\Psi^{\mathrm{BM}}$ 
of $\hat \pi_\Psi$,
c.f. (\ref{170306bb}), by
\[
\hat{\pi}_{\Psi}^{\mathrm{BM}}(x,dz)\overset{\mathrm{def}}{=}\int\frac{1}{m_{x}%
}\varphi\left(  t/m_{x}\right)  \pi_{C_{t}}^{\mathrm{BM}}(x,dz)dt.
\]
$\hat{\pi}_{\Psi}^{\mathrm{BM}}(x,dz)$ has a density with respect to Lebesgue
measure which, by an abuse of notation, we write as $\hat{\pi}_{\Psi}%
^{\mathrm{BM}}(x,z).$
\begin{lemma}
\label{Le_Lipshitz_Pi} There is a constant $C$ such that for any $L$ large
enough, any $\Psi\in\mathcal{M}_{L}$, any $x,x^{\prime},z,z^{\prime}\in
\mathbb{Z}^{d}$, it holds that
\begin{align}
&  \hat{\pi}_{\Psi}\left(  x,z\right)  \leq CL^{-d}\,,\quad\hat{\pi}%
_{\Psi}^{\mathrm{BM}}(x,z)\leq CL^{-d}\,.\label{eq-200305b}\\
&  |\hat{\pi}_{\Psi}\left(  x,z\right)  -\hat{\pi}_{\Psi}\left(  x^{\prime
},z\right)  |\leq C|x-x^{\prime}|L^{-(d+1)}\log L\,,\label{eq-200305ab}\\
&  |\hat{\pi}_{\Psi}^{\mathrm{BM}}(x,z)-\hat{\pi}_{\Psi}^{\mathrm{BM}}(x^{\prime
},z)|\leq C|x-x^{\prime}|L^{-(d+1)}\log L\,,\label{eq-200305bBM}\\
&  |\hat{\pi}_{\Psi}\left(  x,z\right)  -\hat{\pi}_{\Psi}\left(  x,z^{\prime
}\right)  |\leq C|z-z^{\prime}|L^{-(d+1)}\log L\,,\label{eq-200305c}\\
&  |\hat{\pi}_{\Psi}^{\mathrm{BM}}(x,z)-
\hat{\pi}_{\Psi}^{\mathrm{BM}}(x,z^{\prime
})|\leq C|x-x^{\prime}|L^{-(d+1)}\log L\,. \label{eq-200305cBM}%
\end{align}
Further, for $1<a<b<2,$ and $aL\leq\left\vert x-z\right\vert \leq bL,$%
\begin{equation}
\hat{\pi}_{\Psi}\left(  x,z\right)  \geq C\left(  a,b\right)  ^{-1}L^{-d}.
\label{Est_PiHat3}%
\end{equation}
%and, for any $A\in \partial V_L$, $\beta>1/2$, $y\in V_L$ with
%$d(y,\partial V_L)>L^\beta$ and $L>L_0$,
%\begin{equation}
%\label{eq-200305a}
%\sum_{\scriptscriptstyle y'\in A}
%\pi_{\scriptscriptstyle V_{L}}\left(  y,y^{\prime}\right)
%\leq
%\int_{\scriptscriptstyle  d(y',A)\leq L^\beta}
%\pi_{\scriptscriptstyle L}^{\mathrm{BM}}
%\left(  y,dy^{\prime}\right)
%\left(1+\frac{C(\beta)}{L^{(\beta-1/2)/4}}
%\right)+\frac{C(\beta)}{L^{d+1}}\,.
%\end{equation}Similarly, for any $A\in \partial C_L$, $\beta>1/2$,
%$y\in C_L$ with
%$d(y,\partial C_L)>L^\beta$ and $L>L_0$,
%\begin{equation}
%\label{eq-050505d}
%\int_{\scriptscriptstyle  A}
%\pi_{\scriptscriptstyle L}^{\mathrm{BM}}
%\left(  y,dy^{\prime}\right)
%\leq
%\sum_{\scriptscriptstyle y': d(y',A)\leq L^\beta}
%\pi_{\scriptscriptstyle V_{L}}\left(  y,y^{\prime}\right)
%\left(1+\frac{C(\beta)}{L^{(\beta-1/2)/4}}
%\right)+\frac{C(\beta)}{L^{d+1}}\,.
%\end{equation}
%\begin{align}
%\hat{\pi}_{L}\left(  x\right)   &  \leq CL^{-d},\label{Est_PiHat1}\\
%\left\vert \hat{\pi}_{L}\left(  x\right)  -\hat{\pi}_{L}\left(  x^{\prime
%}\right)  \right\vert  &  \leq C|x-x^{\prime}|L^{-d-1}\log
%L\,,\label{Est_PiHat2}%
%\end{align}
%For the Brownian motion case, one has%
%\begin{equation}
%\left\vert \hat{\pi}_{L}^{\mathrm{BM}}(x)-\hat{\pi}_{L}^{\mathrm{BM}%
%}(x^{\prime})\right\vert \leq C|x-x^{\prime}|L^{-d-1}\log L\label{Est_PiHatBM}%
%\end{equation}
%for all $x,x^{\prime}\in\mathbb{R}^{d}.$

\end{lemma}
\begin{proof}
[Proof of Lemma \ref{Le_Lipshitz_Pi}]The estimates (\ref{eq-200305ab}) and
(\ref{Est_PiHat3}) are immediate from Lemmas \ref{Le_Lawler_Exit} and
\ref{Le_ExitsBM}, and the definition of $\hat{\pi}_{\Psi}$.

%We prove only the statements involving
%the random walk, the estimate (\ref{Est_PiHatBM}) can be proved similarly or
%by directly using the Poisson formula, see e.g. \cite{Lawler}.

We turn to the proof of (\ref{eq-200305b}) and (\ref{eq-200305c}). It clearly
suffices to consider only the cases $\left\vert x-x^{\prime}\right\vert =1$ or
$\left\vert z-z^{\prime}\right\vert =1$. Note first that
\begin{align*}
&  |\hat{\pi}_{\Psi}\left(  x,z\right)  -\hat{\pi}_{\Psi}\left(  x^{\prime
},z\right)  |=\left[  1-\frac{m_{x}}{m_{x^{\prime}}}\right]  \hat{\pi}%
_{\Psi}\left(  x,z\right) \\
&  \quad+\frac{1}{m_{x^{\prime}}}\int_{\mathbb{R}^{+}}\left[  \varphi\left(
\frac{t}{m_{x}}\right)  -\varphi\left(  \frac{t}{m_{x^{\prime}}}\right)
\right]  \pi_{\scriptscriptstyle V_{t}(x)}(x,z)dt\,\\
&  \quad+\frac{1}{m_{x^{\prime}}}\int_{\mathbb{R}^{+}}\varphi\left(  \frac
{t}{m_{x^{\prime}}}\right)  \left[  \pi_{\scriptscriptstyle V_{t}(x)}%
(x,z)-\pi_{\scriptscriptstyle V_{t}(x^{\prime})}(x^{\prime},z)\right]
dt\overset{\mathrm{def}}{=}I_{1}+I_{2}+I_{3}\,.
\end{align*}
Since $\Psi\in\mathcal{M}_{L}$, it holds that $\left[  1-\frac{m_{x}%
}{m_{x^{\prime}}}\right]  \leq CL^{-1}|x-x^{\prime}|$, and hence, using
(\ref{eq-200305ab}), it holds that
\begin{equation}
I_{1}\leq CL^{-d}\frac{|x-x^{\prime}|}{L}\,. \label{eq-220305b}%
\end{equation}
Similarly, using the smoothness of $\varphi$ and the estimates $m_{x^{\prime}%
}\geq L/2$ and $\pi_{V_{t}(x)}(x,z)\leq CL^{1-d}$, see Lemma
\ref{Le_Lawler_Exit} a), one gets
\begin{equation}
I_{2}\leq CL^{-d}\frac{|x-x^{\prime}|}{L}\,. \label{eq-220305c}%
\end{equation}
By translation invariance of simple random walk, we have that $\pi_{V_{r}%
(x)}(x,z)=\pi_{V_{r}}(0,z-x)$. Thus, both (\ref{eq-200305b}) and
(\ref{eq-200305c}) will follow if we can show, for $|x-x^{\prime}|=1$ and
$y=x$ or $x^{\prime}$, the estimate
\begin{equation}
\left\vert \int_{\mathbb{R}^{+}}\varphi\left(  \frac{t}{m_{y}}\right)  \left[
\pi_{\scriptscriptstyle V_{t}}(0,z-x)-\pi_{\scriptscriptstyle V_{t}%
}(0,z-x^{\prime})\right]  dt\right\vert \leq CL^{-d}\,. \label{eq-070605a}%
\end{equation}
Of course, we may assume that $|x-z|$ is of order $L$. Note that the
integration in (\ref{eq-070605a}) is over the union of two intervals, each of
length at most $\sqrt{d}$. Hence,
%Let $J\overset{\mathrm{def}}{=}\left\{
%t>0:x\in\partial V_{t}\right\}  .$ $J$ is an interval of length at maximum
%$\sqrt{d}.$ In particular,
due to the smoothness of $\varphi$, (\ref{eq-070605a}) will follow if we can
show that
\begin{equation}
\left\vert \int_{\mathbb{R}^{+}}\left[  \pi_{\scriptscriptstyle V_{t}%
}(0,z-x)-\pi_{\scriptscriptstyle V_{t}}(0,z-x^{\prime})\right]  dt\right\vert
\leq CL^{-d}\,. \label{eq-070605aa}%
\end{equation}
Let $J\overset{\mathrm{def}}{=}\left\{  t>0:x-z\in\partial V_{t}\right\}  .$
$J$ is an interval of length at most $\sqrt{d}.$ For $t\in J,$ we set%
\[
t^{\prime}=t^{\prime}\left(  t\right)  \overset{\mathrm{def}}{=}\left\vert
x^{\prime}-t\frac{z-x}{\left\vert z-x\right\vert }\right\vert .
\]
Evidently, $dt^{\prime}/dt=1+O\left(  L^{-1}\right)  ,$ and if we set
$J^{\prime}\overset{\mathrm{def}}{=}\left\{  t>0:x^{\prime}-z\in\partial
V_{t^{\prime}}\right\}  ,$ then $J^{\prime}$ is an interval of the same length
as $J,$ up to $O\left(  L^{-1}\right)  $, and further $|J\Delta J^{\prime
}|=O\left(  L^{-1}\right)  $. Therefore, if we prove%
\begin{equation}
\left\vert \int_{J\cap J^{\prime}}[\pi_{\scriptscriptstyle V_{t}(x)}%
(x,z)-\pi_{\scriptscriptstyle V_{t^{\prime}}(x^{\prime})}(x^{\prime
},z)]dt\right\vert \leq CL^{-d}\log L, \label{eq-070605b}%
\end{equation}
the estimate (\ref{eq-070605a}) will follow. To abbreviate notation, we write
$V$ for $V_{t}\left(  x\right)  ,$ and $V^{\prime}$ for $V_{t^{\prime}}\left(
x^{\prime}\right)  .$ A first exit decomposition yields%
\begin{equation}
\pi_{V}(x,z)\leq\pi_{V^{\prime}}(x,z)+\sum_{y\in V\backslash V^{\prime}}%
P_{x}^{\mathrm{RW}}\left(  T_{y}<\tau_{V}\right)  \pi_{V}\left(  y,z\right)  .
\label{eq-080605b}%
\end{equation}
%By \cite{Lawler}, we have $\pi_{V^{\prime}}(x,z)=
%\pi_{V^{\prime}}(x^{\prime
%},z)+O\left(  L^{-d}\right)  .$
We have two simple geometric facts:

\begin{itemize}
\item
\[
\bigcup\nolimits_{t\in J\cap J^{\prime}}\left(  V\backslash V^{\prime}\right)
\subset x+\operatorname*{Shell}\nolimits_{L}\left(  C\right)  .
\]

\item For any $y\in x+\operatorname*{Shell}\nolimits_{L}\left(  C\right)  $%
\[
\int_{J\cap J^{\prime}}1_{\left\{  y\in V\backslash V^{\prime}\right\}
}dt\leq C\frac{\left\vert y-z\right\vert }{L}.
\]

\end{itemize}

Using this together with
%Lemma \ref{Le_Est_Exits}, we get%
%By \cite{Lawler}, we have
$\pi_{V^{\prime}}(x,z)=\pi_{V^{\prime}}(x^{\prime},z)+O\left(  L^{-d}\right)
,$ see \cite[Theorem 1.7.1]{Lawler}, we deduce from (\ref{eq-080605b}) that
\begin{align*}
\int_{J\cap J^{\prime}}\pi_{\scriptscriptstyle V_{t}(x)}(x,z)dt  &  \leq
\int_{J\cap J^{\prime}}\pi_{\scriptscriptstyle V_{t^{\prime}}(x^{\prime}%
)}(x^{\prime},z)dt+O\left(  L^{-d}\right)  +CL^{-d} \!\!\!\!\!\!\!\sum_{y\in
x+\operatorname*{Shell}\nolimits_{L}\left(  C\right)  }\left\vert
y-z\right\vert ^{-d}\frac{\left\vert y-z\right\vert }{L}\\
&  \leq\int_{J\cap J^{\prime}}\pi_{\scriptscriptstyle V_{t^{\prime}}%
(x^{\prime})}(x^{\prime},z)dt+O\left(  L^{-d}\log L\right)
\end{align*}
The inequality in the opposite direction is proved in the same way. This
proves (\ref{eq-070605aa}) and completes the proof of (\ref{eq-200305b}) and
(\ref{eq-200305c}).

The estimates (\ref{eq-200305bBM}) and (\ref{eq-200305cBM}) can be obtained
either by repeating the argument above, replacing the random walk by Brownian
motion, or by applying the Poisson formula \cite[(1.43)]{Lawler}. We omit
further details.
\end{proof}

In order to prove Lemma \ref{Le_Approx_Phi_by_BM} we need also the following
technical result:
\begin{lemma}
\label{Le_ComparisonBrown}There exists a constant $C=C(\beta,\epsilon)$ such
that for any $A\in\partial V_{L}$, $\beta>6\epsilon>0$, $y\in V_{L}$ with
$d(y,\partial V_{L})>L^{\beta}$ and $L>L_{0}$,
\begin{equation}
\sum_{y^{\prime}\in A}\pi_{{L}}\left(  y,y^{\prime
}\right)  \leq\int_{d(y^{\prime},A)\leq L^{\beta}}\pi_{\scriptscriptstyle L}%
^{\mathrm{BM}}\left(  y,dy^{\prime}\right)  \left(  1+\frac{C(\beta,\epsilon)
)}{L^{\beta-5\epsilon}}\right)  +\frac{C(\beta,\epsilon)} {L^{d+1}}\,.
\label{eq-200305a}%
\end{equation}
and for any $A^{\prime}\in\partial C_{L}$ and $z\in V_{L}$ with $d(z,\partial
C_{L})>L^{\beta}$,
\begin{equation}
\int_{A^{\prime}}\pi_{\scriptscriptstyle L}^{\mathrm{BM}}\left(  z,dy^{\prime
}\right)  \leq\sum_{y^{\prime}:d(y^{\prime},A)\leq L^{\beta}}\pi
_{{L}}\left(  z,y^{\prime}\right)  \left(  1+\frac
{C(\beta,\epsilon)}{L^{\beta-5\epsilon}}\right)  + \frac{C(\beta,\epsilon
)}{L^{d+1}}\,. \label{eq-050505d}%
\end{equation}
Finally, for any $x,z\in\mathbb{Z}^{d}$ and $\Psi\in\mathcal{M}_{L}$,
\begin{equation}
\left\vert \hat{\pi}_{\Psi}\left(  x,z\right)  -\hat{\pi}_{\Psi}^{\mathrm{BM}%
}(x,z)\right\vert \leq\frac{C}{L^{d+1/4}}. \label{eq-200305cc}%
\end{equation}
\end{lemma}
\begin{proof}
[Proof of Lemma \ref{Le_ComparisonBrown}]We first prove (\ref{eq-200305a}).
Set $A_{\beta}=\{y^{\prime}\in\partial C_{L}:\,d(y^{\prime},A)\leq L^{\beta
}\}$. Pick $\epsilon\in(0,\beta)$
%Set
%$\epsilon=(\beta-1/2)/4$, $\hat{\beta}=1/2+\epsilon/2$,
%noting that $\beta-\hat{\beta}=7(\beta
%-1/2)/8>(\beta-1/2)/4$, and
and set $L^{\prime}=L+L^{\epsilon}$ and $L^{\prime\prime}=L+L^{2\epsilon}$.
Let $A_{\beta}^{\prime}$ be the image of $A_{\beta}$ in $\partial
C_{L^{\prime}}$ under the map $x\mapsto(L^{\prime}/L)x$. Then, one has (with
$\hat{y}=L^{\prime}y/L$),
\begin{equation}
\int_{\scriptscriptstyle A_{\beta}}\pi_{\scriptscriptstyle{L}}^{\mathrm{BM}%
}\left(  y,dy^{\prime}\right)  =\int_{\scriptscriptstyle A_{\beta}^{\prime}%
}\pi_{\scriptscriptstyle{L^{\prime}}}^{\mathrm{BM}}\left(  \hat{y},dy^{\prime
}\right)  \,. \label{eq-050405bold}%
\end{equation}
Note further, using the Poisson formula \cite[(1.43)]{Lawler}, that
\begin{align}
\int_{\scriptscriptstyle A_{\beta}^{\prime}}\pi_{\scriptscriptstyle{L^{\prime
}}}^{\mathrm{BM}}\left(  \hat{y},dy^{\prime}\right)   &  =\int
_{\scriptscriptstyle A_{\beta}^{\prime}}\frac{d\pi
_{\scriptscriptstyle{L^{\prime}}}^{\mathrm{BM}}\left(  \hat{y},\cdot\right)
}{d\pi_{\scriptscriptstyle{L^{\prime}}}^{\mathrm{BM}}\left(  {y},\cdot\right)
}\pi_{\scriptscriptstyle{L^{\prime}}}^{\mathrm{BM}}\left(  {y},dy^{\prime
}\right) \label{eq-noclosing1}\\
&  =\int_{\scriptscriptstyle A_{\beta}^{\prime}}\frac{\left(  (L^{\prime}%
)^{2}-|\hat{y}|^{2}\right)  |y^{\prime}-y|^{d}}{\left(  (L^{\prime}%
)^{2}-|y|^{2}\right)  |y^{\prime}-\hat{y}|^{d}}\pi
_{\scriptscriptstyle{L^{\prime}}}^{\mathrm{BM}}\left(  {y},dy^{\prime}\right)
\nonumber
\end{align}
An explicit computation, using that $|y|\leq L-L^{\beta}$ and that
$1>\beta>\epsilon>0$, reveals that
\[
\left\vert \log\frac{\left(  (L^{\prime})^{2}-|\hat{y}|^{2}\right)
|y^{\prime}-y|^{d}}{\left(  (L^{\prime})^{2}-|y|^{2}\right)  |y^{\prime}%
-\hat{y}|^{d}}\right\vert \leq CL^{\epsilon-\beta}\,.
\]
Substituting in (\ref{eq-noclosing1}) one finds that
\begin{equation}
\int_{\scriptscriptstyle A_{\beta}}\pi_{\scriptscriptstyle{L}}^{\mathrm{BM}%
}\left(  y,dy^{\prime}\right)  \geq\int_{\scriptscriptstyle A_{\beta}^{\prime
}}\pi_{\scriptscriptstyle{L^{\prime}}}^{\mathrm{BM}}\left(  {y},dy^{\prime
}\right)  \left(  1-C(\beta,\epsilon)L^{-\beta+2\epsilon}\right)  \,.
\label{eq-050405b}%
\end{equation}
Recall that $\pi_{L}^{\mathrm{BM}}$ is unchanged if one replaces the Brownian
motion by a Brownian motion of covariance $I_{d}/\sqrt{d}$. Let $W_{t}^{y}$ be
such a Brownian motion started at $y$, and recall that by \cite[Corollary
1]{zaitsev}, there exists a constant $C_{0}$ such that for every integer $n$,
one may construct $\{W_{t}^{x}\}$ in the same space as $\{X_{n}\}$ such that
\begin{equation}
\label{eq-040605a}P_{x}(\max_{0\leq m\leq n}|X_{m}- W_{m}^{x}|>C_{0}\log n)
\leq\frac{C_{0}}{n^{d+1}}\,.
\end{equation}
%(\ref{eq-040605a}) and the
%constant $C_{0}$ defined there. The latter inequality, and
Standard estimates involving the maximum of the increments of the Brownian
motion, imply that one may construct the Brownian motion $W_{t}^{y}$ and the
random walk $X_{n}$ on the same space such that, with
\[
D\overset{\mathrm{def}}{=}\{\sup_{\scriptscriptstyle0\leq t\leq L^{2+\epsilon
/100}}\left\vert X_{\left[  t\right]  }-W_{t}^{y}\right\vert \leq4C_{0}\log
L\}\,,
\]
one has
\begin{equation}
P_{y}(D^{c})\leq\frac{2C_{0}}{n^{d+1}}\,. \label{eq-040605anew}%
\end{equation}
Set $\tau\overset{\mathrm{def}}{=}\min\{n:X_{n}\in\partial V_{L}\}$,
$\tau^{\prime}\overset{\mathrm{def}}{=}\inf\{t:W_{t}^{y}\in\partial
C_{L^{\prime}}\}$, $\tau^{\prime\prime}\overset{\mathrm{def}}{=}\min
\{n:X_{n}\in\partial V_{L^{\prime\prime}}\}$, and $B\overset{\mathrm{def}}%
{=}\left\{  \left(  \tau^{\prime}\vee\tau^{\prime\prime}\right)  \leq
L^{2+\epsilon/100}\right\}  $. Standard estimates imply that if $X_{0}=y$ then
$P(B^{c})$ decays like a stretched exponential, and in particular
$P(B^{c})\leq L^{-d-1}$ for large $L.$ Note that on $D\cap B$, one has that
$\tau<\tau^{\prime}<\tau^{\prime\prime}$.
Now, defining $G_{\beta}^{\prime}=\{z\in\mathbb{Z}^{d}:d(z,(A_{\beta}^{\prime
})^{c}\cap\partial C_{L})<4C_{0}\log L\}$, and setting $T_{G_{\beta}^{\prime}%
}=\inf\{n:X_{n}\in G_{\beta}^{\prime}\}$,
\begin{align}
P\left(  W_{\tau^{\prime}}^{y}\in A_{\beta}^{\prime}\right)   &  \geq
P_{y}(X_{\tau}\in A,W_{\tau^{\prime}}\in A_{\beta}^{\prime}%
)\label{eq-080605dd}\\
&  \geq P_{y}(X_{\tau}\in A,W_{\tau^{\prime}}\in A_{\beta}^{\prime},B\cap
D)-\frac{1}{L^{d+1}}\nonumber\\
&  \geq P_{y}(X_{\tau}\in A)-P_{y}(X_{\tau}\in A,W_{\tau^{\prime}}%
\not \in A_{\beta}^{\prime},B\cap D)-\frac{2}{L^{d+1}}\nonumber\\
&  \geq P_{y}^{\mathrm{RW}}(X_{\tau}\in A)-P_{y}^{\mathrm{RW}}(X_{\tau}\in
A,T_{G_{\beta}^{\prime}}<\tau^{\prime\prime})-\frac{2}{L^{d+1}}\nonumber
\end{align}
Using the Markov property, one has
\begin{align*}
&P_{y}^{\mathrm{RW}}(X_{\tau}\in A,T_{G_{\beta}^{\prime}}<\tau^{\prime\prime})
  \leq P_{y}^{\mathrm{RW}}(X_{\tau}\in A)\sup_{z\in A}P_{z}^{\mathrm{RW}%
}(T_{G_{\beta}^{\prime}}<\tau^{\prime\prime})\\
&  \leq\sup_{z\in A}\sum_{z^{\prime}\in G_{\beta}^{\prime}}P_{z}^{\mathrm{RW}%
}(T_{z^{\prime}}<\tau^{\prime\prime})
  \leq\sup_{z\in A}C\sum_{z^{\prime}\in G_{\beta}^{\prime}}\frac
{L^{3\epsilon}\log^{d+2}L}{|z^{\prime}-z|^{d}}
  \leq CL^{5\epsilon-\beta}\,,
\end{align*}
where the next to last inequality is due to Lemma \ref{Le_MainExit}.
Substituting in (\ref{eq-080605dd}), one completes the proof of
(\ref{eq-200305a}). The 
reverse inequality (\ref{eq-050505d}) is proved similarly.

It remains to prove (\ref{eq-200305cc}). Fix $\alpha=2/3$, $\beta=1/3$, and
$\epsilon=1/60$. Note that with $\mathcal{D}=C_{L^{\alpha}}(z)$, using
(\ref{eq-200305c}),
\begin{equation}
\hat{\pi}_{\Psi}\left(  x,z\right)  \leq\frac{1}{|\mathcal{D}|} \sum_{z^{\prime
}\in\mathcal{D}}\hat{\pi}_{\Psi}\left(  x,z^{\prime}\right)  +CL^{-d-1+\alpha
}\log L\,. \label{eq-050505e}%
\end{equation}
Next, note that
\begin{align*}
\sum_{z^{\prime}\in\mathcal{D}} \hat{\pi}_{\Psi}\left(  x,z^{\prime}\right)   &
=\int dt\varphi_{m_{x}}(t)\sum_{z^{\prime}\in\mathcal{D}} {\pi}%
_{\scriptscriptstyle V_{t}(x)}\left(  x,z^{\prime}\right) \\
&  \leq\int dt\varphi_{m_{x}}(t)\int_{\scriptscriptstyle C_{ L^{\alpha
}+L^{\beta}}(z)}{\pi}_{t}^{\mathrm{BM}}\left(  x,dz^{\prime}\right)  \left(
1+\frac{C}{L^{\beta-5\epsilon}}\right)  +\frac{C|\mathcal{D}|}{L^{d+1}}\\
&  \leq\hat{\pi}_{\Psi}^{\mathrm{BM}}(x,\mathcal{D})
\left(  1+\frac{C}{L^{\beta-5
\epsilon} }\right)  +CL^{-d}|C_{L^{\alpha}+L^{\beta}}(z)\setminus
C_{L^{\alpha}}(z) |+\frac{C|\mathcal{D}|}{L^{d+1}}\\
&  \leq|\mathcal{D}|\hat{\pi}_{\Psi}^{\mathrm{BM}}(x,z)\left(  1+\frac
{C}{L^{\beta-5\epsilon}}\right)  +\frac{C|\mathcal{D}|\log L}{L^{d+1-\alpha}%
}+\frac{C|\mathcal{D}|} {L^{\alpha-\beta-d}}\,.
\end{align*}
Substituting in (\ref{eq-050505e}), one gets
\[
\hat{\pi}_{\Psi}\left(  x,z\right)  \leq\hat{\pi}_{\Psi}^{\mathrm{BM}}%
(x,z)+CL^{-d-1/4}\,.
\]
The reverse equality is proved similarly. This completes the proof of
(\ref{eq-200305cc}) and of the lemma
\end{proof}

\noindent
\begin{proof}
[Proof of Lemma \ref{Le_Approx_Phi_by_BM}]Fix $\alpha=2/3,\beta=1/3$. Set
$\eta\overset{\mathrm{def}}{=}d(y,\partial V_{L})$, and let $y_{1}\in\partial
V_{L}$ be such that $\eta=|y-y_{1}|$. Consider first $\eta\leq L^{\beta+1/15}
$. Then, using (\ref{eq-080605gg}) and (\ref{eq-200305ab}) in the first
inequality and (\ref{eq-200305b}) in the second,
\begin{align*}
\phi_{L,\Psi}\left(  y,z\right)   &  \leq\sum_{y^{\prime}\in\partial
V_{L}:|y^{\prime}-y_{1}|<L^{\alpha}}\pi_{\scriptscriptstyle V_{L}}\left(
y,y^{\prime}\right)  \hat{\pi}_{\Psi}\left(  y^{\prime},z\right)  +\frac{C\log
L}{L^{d+\alpha-\beta}}\\
&  \leq\hat{\pi}_{\Psi}\left(  y_{1},z\right)  \sum_{y^{\prime}\in\partial
V_{L}:|y^{\prime}-y_{1}|<L^{\alpha}}\pi_{\scriptscriptstyle V_{L}}\left(
y,y^{\prime}\right)  +\frac{C}{L^{d+1/5}}\,.
\end{align*}
Consequently,
\[
\phi_{L,\Psi}\left(  y,z\right)  \leq\hat{\pi}_{\Psi}\left(  y_{1},z\right)
+\frac{C}{L^{d+1/5}}\,.
\]
Applying now (\ref{eq-200305ff}) in the first inequality and (\ref{eq-200305b}%
) in the second, we conclude that
\begin{align*}
\phi_{L,\Psi}\left(  y,z\right)   &  \leq\hat{\pi}_{\Psi}\left(  y_{1},z\right)
\int_{y^{\prime}\in\partial V_{L}:|y^{\prime}-y_{1}|<L^{\alpha}}%
\pi_{\scriptscriptstyle L}^{\mathrm{BM}}\left(  y,dy^{\prime}\right)
+\frac{C}{L^{d+1/5}}\\
&  \leq\int_{y^{\prime}\in\partial V_{L}:|y^{\prime}-y_{1}|<L^{\alpha}}%
\hat{\pi}_{\Psi}\left(  y^{\prime},z\right)  \pi_{\scriptscriptstyle L}%
^{\mathrm{BM}}\left(  y,dy^{\prime}\right)  +\frac{C}{L^{d+1/5}}\,.
\end{align*}
An application of (\ref{eq-200305cc}) then implies that for $\eta\leq
L^{\beta+1/15}$,
\[
\phi_{L,\Psi}\left(  y,z\right)  \leq\phi_{L,\Psi}^{\mathrm{BM}}\left(  y,z\right)
+CL^{-d-1/5}%
\]
where, as in our convention, the constant $C$ is uniform in the choice of
$y,z$. The reverse inequality is obtained using the same steps.

Consider next $\eta>L^{\beta+1/15}$. Fix strictly positive constants $c_{j} $,
$j=1,\ldots,4$, depending on $d,\alpha$ only, and a sequence of disjoint sets
$A_{i}\subset\partial V_{L}$, $i=1,\ldots,k_{L}$ with $\cup_{i=1}^{k_{L}}%
A_{i}=\partial V_{L}$, $c_{1}L^{\alpha(d-1)}\leq|A_{i}|\leq c_{2}%
L^{\alpha(d-1)}$, $\mbox{\rm diam}(A_{i})\leq c_{3}L^{\alpha}$, $d(y_{1}%
,\partial A_{1}\cap\partial V_{L})\geq\mathrm{diam}(A_{1})/4$, and $|\partial
A_{i}|\cap\partial V_{L}\leq c_{4}L^{\alpha(d-2)}$ (such a collection of
\textquotedblleft cube-like\textquotedblright\ $A_{i}$ can clearly be found).
We also set $A_{i}^{\beta}=\{y\in\mathbb{R}^{d}:\,d(y,A_{i})\leq L^{\beta}\}$
and for $i\geq2$, fix an arbitrary $y_{i}\in A_{i}$. We then have
\begin{align*}
\phi_{L,\Psi}\left(  y,z\right)   &  =\sum_{i=1}^{k_{L}}\sum_{y^{\prime}\in
A_{i}}\pi_{\scriptscriptstyle V_{L}}\left(  y,y^{\prime}\right)  \hat{\pi}%
_{\Psi}\left(  y^{\prime},z\right) \\
&  \leq\sum_{i=1}^{k_{L}}\hat{\pi}_{\Psi}\left(  y_{i},z\right)  \sum_{y^{\prime
}\in A_{i}}\pi_{\scriptscriptstyle V_{L}}\left(  y,y^{\prime}\right)
+\frac{C\log L}{L^{d+1-\alpha}}\,,
\end{align*}
where (\ref{eq-200305b}) was used in the last inequality. Consequently, using
(\ref{eq-200305a}),
\begin{equation}
\phi_{L,\Psi}\left(  y,z\right)  \leq\sum_{i=1}^{k_{L}}\hat{\pi}_{\Psi}\left(
y_{i},z\right)  \int_{A_{i}^{\beta}}\pi_{\scriptscriptstyle{L}}^{\mathrm{BM}%
}\left(  y,dy^{\prime}\right)  \left(  1+\frac{C}{L^{1/4}}\right)  +\frac
{C}{L^{d+1/5}}\,. \label{eq-200305e}%
\end{equation}
%
%Recall that, for some constant $c_5$,
%$$
%\pi_{\scriptscriptstyle {L}}^{\mathrm{BM}}
%\left(  y,y^{\prime}\right)\leq \frac{c_5}{L^{d-1}}\,.$$
%Consequently,
Let $\{\tilde{A}_{i}\subset\partial C_{L}\}_{i=1}^{\scriptscriptstyle k_{L}} $
be a collection of measurable disjoint sets with $\cup\tilde{A}_{i}=\partial
C_{L}$, $\tilde{A}_{1}=A_{1}^{\beta}\cap\partial C_{L}$, and $\tilde{A}%
_{i}\subset A_{i}^{\beta}$. Using (\ref{eq-200305ff}) and $d(y,\partial
C_{L})\geq L^{\beta+1/15}/2$, one gets
\[
\int_{A_{i}^{\beta}}\pi_{\scriptscriptstyle{L}}^{\mathrm{BM}}\left(
y,dy^{\prime}\right)  \leq\int_{\tilde{A}_{i}}\pi_{\scriptscriptstyle{L}%
}^{\mathrm{BM}}\left(  y,dy^{\prime}\right)  \left(  1+C\frac{|(A_{i}^{\beta
}\cap\partial C_{L})\setminus\tilde{A}_{i}|}{|A_{i}^{\beta}\cap\partial
C_{L}|}\right)  \,.
\]
Substituting in (\ref{eq-200305e}) we get%
\[
\phi_{L,\Psi}\left(  y,z\right)  \leq\sum_{i=1}^{k_{L}}\hat{\pi}_{\Psi}\left(
y_{i},z\right)  \int_{\tilde{A}_{i}}\pi_{\scriptscriptstyle{L}}^{\mathrm{BM}%
}\left(  y,dy^{\prime}\right)  \left(  1+CL^{-1/5}\right)  +\frac{C}%
{L^{d+1/5}}\,.
\]
Hence, recalling (\ref{eq-200305ab}), (\ref{eq-200305b}), and
(\ref{eq-200305cc}), we get
\begin{align*}
\phi_{L,\Psi}\left(  y,z\right)   &  \leq\sum_{i=1}^{k_{L}}\int_{\tilde{A}_{i}%
}\hat{\pi}_{\Psi}\left(  y^{\prime},z\right)  \pi_{\scriptscriptstyle{L}%
}^{\mathrm{BM}}\left(  y,dy^{\prime}\right)  +\frac{C}{L^{d+1/5}}\\
&  \leq\sum_{i=1}^{k_{L}}\int_{\tilde{A}_{i}}\hat{\pi}_{\Psi}^{\mathrm{BM}%
}\left(  y^{\prime},z\right)  \pi_{\scriptscriptstyle{L}}^{\mathrm{BM}}\left(
y,dy^{\prime}\right)  +\frac{C}{L^{d+1/5}}=\phi_{\scriptscriptstyle L,\Psi}%
^{\mathrm{BM}}(y,z)+\frac{C}{L^{d+1/5}}\,.
\end{align*}
The reverse inequality is obtained by a similar argument.
\end{proof}

\noindent
\begin{proof}
[Proof of Lemma \ref{Le_ThirdDerivative}]We write $\pi_{t}^{\mathrm{BM}%
}\left(  w,z\right)  $ as the density on $\partial C_t(w)$
%with respect to Lebesgue's measure 
of
the measure $\pi_{C_{t}\left(  w\right)  }^{\mathrm{BM}}\left(  w,dz\right)
$. Set $g(w,z)=\int\pi_{t}^{\mathrm{BM}}\left(  w,z\right)  \varphi_{m_{w}%
}\left(  t\right)  dt\,.$ Then,
\[
%u(y,z)  \overset{\mathrm{def}}{=}%
\phi_{L,\Psi}^{\mathrm{BM}}\left(  y,z\right)  =\int_{\partial C_{L}\left(
0\right)  }\pi_{C_{L}\left(  0\right)  }^{\mathrm{BM}}\left(  y,dw\right)
g(w,z)\,.
\]
For $\bar z\in \partial C_1(0)$, 
set
%$u(y,z)=\phi_{L,\Psi}^{\mathrm{BM}}\left(  Ly,Lz\right)  $
%satisfies the equation
\[
\left\{
\begin{array}
[c]{ll}%
\frac{1}{2}\Delta_{\bar y}u(\bar y,\bar z)=0, & \bar y\in C_{1}(0)\,,\\
u(\bar y,\bar z)=g(L\bar y,L\bar z), & \bar y\in\partial C_{1}(0)\,.
\end{array}
\right.
\]
Then,
$\phi_{L,\Psi}^{\mathrm{BM}}\left(  y,z\right)=u(y/L,\bar z)  $ with
$\bar z=z/L$ 
and hence
\begin{equation}
\label{eq-150306a}
\left|\frac{\partial^i \phi_{L,\Psi}^{\mathrm{BM}}(y,z)}{\partial y^i}\right|
=\frac{1}{L^{i}}\left| \frac{\partial^i u(\bar y,\bar z)}{\partial \bar y^i}
\right|\,.
\end{equation} 
Write 
$$\|u(\bar y,\bar z)\|_k=
\sum_{j=0}^k 
\sup_{\bar y, \bar z} \left|\frac{ \partial^j u(\bar y,\bar z)}{\partial
\bar y^j}\right|\,.$$
By \cite[Theorem 6.3.2]{krylov}, 
\begin{equation}
\|u\left(  \bar y,\bar z\right)
\|_3 \leq C  \| g(\bar w,\bar z)\|_4\,.
  \label{eq-090605a}%
\end{equation}
By the smoothness of $\varphi$ and the translation invariance and scaling
properties of the Brownian motion, 
and applying
\cite[Theorem 2.10]{gilbarg}, one gets that
$$\| g(\bar w,\bar z)\|_4\leq L^{-d}\,.$$
Substituting in (\ref{eq-090605a}) and using (\ref{eq-150306a}), the
lemma follows.
\end{proof}

\section{A local CLT and proof of Lemma \ref{Cor_Green}}

\label{App_B}

We need a number of properties for simple random walk, and coarse-grained
random
walks, which can readily be obtained from known results. We keep $L$ and $V_{L}$
fixed through this section, and don't emphasize them in the notation. $\pi$ is
$\pi_{V_{L}},$ the exit distribution of simple random walk from $V_{L}.$ Since
the proofs are very similar, and for concreteness, we prove all results for
the smoothing scheme $\mathcal{S}_{1}$ and only sketch the necessary changes
for the scheme $\mathcal{S}_{2}$. 
%That is, we take:
%as in Section
%\ref{SubSect_Smoothing} without
%\textquotedblleft bad\textquotedblright%
%\ points:%
%\[
%s_{x}\overset{\mathrm{def}}{=}\left\{
%\begin{array}
%[c]{cc}%
%\delta_{k_{0} r\left(  L\right)  } & \mathrm{if\ }x\in{\operatorname*{Shell}%
%\nolimits}_{L} \left(  r(L)\right) \\
%\varphi_{h_{L} \left(  x\right)  }\left(  t\right)  dt & \mathrm{if\ }%
%d_{L}\left(  x\right)  >r\left(  L\right)
%\end{array}
%\right.  .
%\]
%($r(L)=L/\left(  \log L\right)  ^{10}$). 
Remark that the coarse graining scale at $x$, $h_{L}\left(
x\right)$, equals $\gamma s(L)$
%=\gamma L/\left(  \log L\right)  ^{3}$ 
for $d_{L}\left(
x\right)  \geq2s(L),$ and 
 $h_{L}\left(  x\right)  \leq\left(  \gamma/2\right)
s(L)$ for $x\in{\operatorname*{Shell}\nolimits}_{L} \left(  r,2s(L)\right)  .$
%We then write $\hat{\pi}_{\mathcal{S}}$ for the corresponding transition
%probabilities. 
By a slight abuse of notation, we write $\hat{\pi}_{m}$ for the
transition probabilities on $\mathbb{Z}^{d}$ with the
constant in $x$ coarse-graining scheme $\Psi_m$.
%
%a smoothing scheme $\left(
%s_{x}\right)  $ which is constant in $x,$ and given by $\varphi_{m}\left(
%t\right)  dt.$ 
We also write $\hat{\pi}_{m}\left(  x\right)  $ for $\hat{\pi
}_{m}\left(  0,x\right)  .$ For $x\in V_{L-2s(L)},$ $\hat{\pi}%
\left(  x,\cdot\right)  =\hat{\pi}_{\gamma s(L)}\left(  x,\cdot\right)  $
under either $\mathcal{S}_i$, $i=1,2$.

Let $m\in\mathbb{R}^{+}.$ $\hat{\pi}_{m}$ is centered, and the covariances
satisfy%
\[
\sum_{x}x_{i}x_{j}\hat{\pi}_{m}\left(  x\right)  =\alpha\left(  m\right)
\delta_{ij},
\]
where for some $0<\alpha_{1}<\alpha_{2}$%
\[
\alpha_{1}m^{2}\leq\alpha\left(  m\right)  \leq\alpha_{2}m^{2}.
\]
(It is evident that $\sigma_m^2
\overset{\mathrm{def}}{=}
\alpha\left(  m\right)  /m^{2}$ converges as
$m\rightarrow\infty.)$ 
%We also write $\sigma_m^2=\alpha(m)$.
Using Lemma \ref{Le_Lawler_Exit} a), one sees that for $1<a<b<2,$ one has for
some $\delta$ (which may depend on $a,b)$%
\begin{equation}
\inf_{am\leq\left\vert x\right\vert \leq bm}\hat{\pi}_{m}\left(  x\right)
\geq\delta m^{-d}. \label{Est_Pi_from_below}%
\end{equation}
Furthermore, by definition, we have $\hat{\pi}_{m}\left(  x\right)  =0$ for
$\left\vert x\right\vert \geq2m.$

We will also use the following fact, proved in Lemma \ref{Le_Lipshitz_Pi}:%
\begin{equation}
	\label{eq-170603cc}
\left\vert \hat{\pi}_{m}\left(  x\right)  -\hat{\pi}_{m}\left(  y\right)
\right\vert \leq Cm^{-d}\left\vert \frac{x-y}{m}\right\vert ^{1/15}.
\end{equation}

In what follows, we write $\hat{\pi}_m^{\ast n}$ for the $n$-fold convolution
of $\hat{\pi}_m$.
\begin{proposition}
\label{Ap_Prop_LCLT}
%Uniformly in $x\in V_L$,
\[
\hat{\pi}_{m}^{\ast n}\left(  x\right)  =\frac{1}{\left(  2\pi m^{2}\sigma
_{m}^{2}n\right)  ^{d/2}}\exp\left[  -\frac{\left\vert x\right\vert ^{2}%
}{2m^{2}\sigma_{m}^{2}n}\right]  +O\left(  m^{-d}n^{-\left(  d+2\right)
/2}\left(  \log n\right)  ^{4}\right)
\]

\end{proposition}
\begin{proof}
[Proof of Proposition \ref{Ap_Prop_LCLT}]The proof is standard, but we need to
keep track of the $m$-dependence, and we are not aware of a reference for that
in the literature. Let%
\[
\chi_{m}\left(  z\right)  \overset{\mathrm{def}}{=}\sum_{x}\mathrm{e}^{iz\cdot
x/m}\hat{\pi}_{m}\left(  x\right)  ,\ z\in B_{m}\overset{\mathrm{def}}%
{=}\left[  -m\pi,m\pi\right]  ^{d}%
\]
By Fourier inversion, we have%
\[
\hat{\pi}_{m}^{\ast n}\left(  x\right)  =\left(  2\pi\right)  ^{-d}m^{-d}%
\int_{B_{m}}\mathrm{e}^{-iz\cdot x/m}\chi_{m}\left(  z\right)  ^{n}dz.
\]
We will choose $0<a<A,$ $b>0,$ and $\alpha\in\left(  0,1\right)  $ (not
depending on $n,m$) and split%
\begin{align}
\int_{B_{m}}\mathrm{e}^{-iz\cdot x/m}\chi_{m}\left(  z\right)  ^{n}dz  &
=\int\limits_{\left\vert z\right\vert \leq\frac{b\log n}{\sqrt{n}}}%
+\int\limits_{\frac{b\log n}{\sqrt{n}}<\left\vert z\right\vert \leq a}%
+\int\limits_{a<\left\vert z\right\vert \leq A}+\int\limits_{A<\left\vert
z\right\vert \leq m^{\alpha}}+\int\limits_{m^{\alpha}<\left\vert z\right\vert
,\ z\in B_{m}}\nonumber\\
&  \overset{\mathrm{def}}{=}
A_{1}+A_{2}+A_{3}+A_{4}+A_{5}.\nonumber
\end{align}

From Taylor's formula, we get%
\[
\chi_{m}\left(  z\right)  =1-\frac{\left\vert z\right\vert ^{2}\sigma_{m}^{2}%
}{2}+O\left(  \left\vert z\right\vert ^{4}\right)  ,
\]
and therefore, for $\left\vert z\right\vert \leq1/C$,%
\[
\log\chi_{m}\left(  z\right)  =-\frac{\left\vert z\right\vert ^{2}\sigma
_{m}^{2}}{2}+O\left(  \left\vert z\right\vert ^{4}\right)  .
\]
From that we get for $b$ sufficiently large
 and $n\geq C\left(  b\right)  $,%
\begin{align}
A_{1}  &  =\left(  1+O\left(  \frac{\left(  \log n\right)  ^{4}}{n}\right)
\right)  \int\limits_{\left\vert z\right\vert \leq\frac{b\log n}{\sqrt{n}}%
}\exp\left[  -i\frac{z\cdot x}{m}-\frac{n\left\vert z\right\vert ^{2}%
\sigma_{m}^{2}}{2}\right]  dz\nonumber\\
&  =\left(  1+O\left(  \frac{\left(  \log n\right)  ^{4}}{n}\right)  \right)
\int\exp\left[  -i\frac{z\cdot x}{m}-\frac{n\left\vert z\right\vert ^{2}%
\sigma_{m}^{2}}{2}\right]  dz+O\left(  n^{-d/2-1}\right) \label{Est_A1}\\
&  =\frac{\left(  2\pi\right)  ^{d/2}}{n^{d/2}\sigma_{m}^{d}}\exp\left[
-\frac{\left\vert x\right\vert ^{2}}{2m^{2}\sigma_{m}^{2}n}\right]  +O\left(
n^{-d/2-1}\left(  \log n\right)  ^{4}\right)  .\nonumber
\end{align}
In order to prove the proposition, it therefore suffices to prove that
$A_{2},\ldots,A_{5}$ are of order $O\left(  n^{-d/2-1}\right)  ,$ uniformly in
$L.$

To handle $A_{2},$ we choose $a$ such that $\log\chi_{m}\left(  z\right)
\leq-\left\vert z\right\vert ^{2}\sigma_{m}^{2}/3$ for $\left\vert
z\right\vert \leq a.$ Then%
\begin{equation}
\left\vert A_{2}\right\vert \leq\int_{\frac{b\log n}{\sqrt{n}}<\left\vert
z\right\vert }\exp\left[  -\left\vert z\right\vert ^{2}n\sigma_{m}%
^{2}/3\right]  dz=O\left(  n^{-d/2-1}\right)  .\nonumber
\end{equation}
if we choose $b$ sufficiently large.

For $A_{3},$ we use the following fact, which is an easy consequence of
(\ref{Est_Pi_from_below}): for any $a<A,$ one has%
\begin{equation}
\sup_{m,a\leq\left\vert z\right\vert \leq A}\left\vert \chi_{m}\left(
z\right)  \right\vert <1\,. \label{Est_Fourier}%
\end{equation}
Using this, we immediately get%
\begin{equation}
\left\vert A_{3}\right\vert \leq CA^{d}\left(  1-1/C\right)  ^{n}\,.
\label{Est_A3}%
\end{equation}

We come now to $A_{4}$ which is more difficult. First 
remark that since the coarse graining scheme is isotropic,
%by the
%assumed symmetry under lattice isomorphisms, 
we only have to consider
$z$-values with all components positive. Put $\left\vert z\right\vert
_{\infty}\overset{\mathrm{def}}{=}\max\left(  z_{1},\ldots,z_{d}\right)  .$
For simplicity, we assume that $z_{1}$ is the biggest component of $z,$ so
that $\left\vert z\right\vert _{\infty}=z_{1}.$ Let $M\overset{\mathrm{def}%
}{=}\left[  2\pi m/z_{1}\right]  ,$ and $K\overset{\mathrm{def}}{=}\left[
\left(  2m+1\right)  /M\right]  .$ We may assume that $M<m$ by choosing $A$
large enough. We write%
\begin{align*}
\chi_{m}\left(  z\right)   &  =\sum_{\left(  x_{2},\ldots,x_{d}\right)  }%
\exp\left[  \frac{i}{m}\sum\nolimits_{s=2}^{d}x_{s}z_{s}\right] \\
&  \times\left\{  \sum_{j=1}^{K}\sum_{x_{1}=-m+\left(  j-1\right)
M}^{-m+jM-1}\mathrm{e}^{ix_{1}z_{1}/m}\hat{\pi}_{m}\left(  x\right)
+\sum_{x_{1}=-m+KM}^{m}\mathrm{e}^{ix_{1}z_{1}/m}\hat{\pi}_{m}\left(
x\right)  \right\}  .
\end{align*}
In the first summand, inside the $x_{1}$-summation, we write for each $j$
separately, $\hat{\pi}_{m}\left(  x\right)  =\hat{\pi}_{m}\left(  x\right)
-\hat{\pi}_{m}\left(  x^{\prime}\right)  +\hat{\pi}_{m}\left(  x^{\prime
}\right)  ,$ where $x^{\prime}=\left(  -m+\left(  j-1\right)  M,x_{2}%
,\ldots,x_{d}\right)  .$ Then, by (\ref{eq-170603cc}),
%we estimate%
\[
\left\vert \hat{\pi}_{m}\left(  x\right)  -\hat{\pi}_{m}\left(  x^{\prime
}\right)  \right\vert \leq Cm^{-d}\left(  \frac{x_{1}+m-\left(  j-1\right)
M}{m}\right)  ^{1/15}.
\]
Therefore,%
\[
\left\vert \sum\nolimits_{x_{1}=-m+\left(  j-1\right)  M}^{-m+jM-1}%
\mathrm{e}^{ix_{1}z_{1}/m}\left(  \hat{\pi}_{m}\left(  x\right)  -\hat{\pi
}_{m}\left(  x^{\prime}\right)  \right)  \right\vert \leq Cm^{-d+1}\frac
{1}{z_{1}^{16/15}},
\]
and therefore,%
\[
\left\vert \sum_{j=1}^{K}\sum_{x_{1}=-m+\left(  j-1\right)  M}^{-m+jM-1}%
\mathrm{e}^{ix_{1}z_{1}/m}\left(  \hat{\pi}_{m}\left(  x\right)  -\hat{\pi
}_{m}\left(  x^{\prime}\right)  \right)  \right\vert \leq Cm^{-d+1}\left\vert
z\right\vert ^{-1/15}.
\]%
Also,
\[
\left\vert \sum_{j=1}^{K}\sum_{x_{1}=-m+\left(  j-1\right)  M}^{-m+jM-1}%
\mathrm{e}^{ix_{1}z_{1}/m}\hat{\pi}_{m}\left(  x^{\prime}\right)  \right\vert
\leq K\hat{\pi}_{m}\left(  x^{\prime}\right)  \left\vert \frac{1-\exp\left[
iz_{1}M/m\right]  }{1-\exp\left[  iz_{1}/m\right]  }\right\vert \leq
C\left\vert z\right\vert m^{-d},
\]%
\[
\left\vert \sum_{x_{1}=-m+KM}^{m}\mathrm{e}^{ix_{1}z_{1}/m}\hat{\pi}%
_{m}\left(  x\right)  \right\vert \leq m^{-d+1}\left\vert z\right\vert ^{-1}.
\]
Therefore, we get the estimate%
\[
\left\vert \chi_{m}\left(  z\right)  \right\vert \leq C_{1}\left(  \left\vert
z\right\vert ^{-1/15}+\frac{\left\vert z\right\vert }{m}\right)  .
\]
From this, we get%
\begin{equation}
\left\vert A_{4}\right\vert \leq C_{1}^{n}\int_{A\leq\left\vert z\right\vert
\leq m^{\alpha}}\left(  \left\vert z\right\vert ^{-1/15}+\frac{\left\vert
z\right\vert }{m}\right)  ^{n}dz\leq2^{-n} \label{Est_A4}%
\end{equation}
for large enough $A$ and $m.$

For $A_{5},$ we need a slight modification. Let again $z_{1}>0$ be the largest
of the $z$-components. Then we write%
\[
\hat{\pi}_{m}\left(  x\right)  =\sum_{y=-m}^{x_{1}}\left(  \hat{\pi}%
_{m}\left(  y,x_{2},\ldots,x_{d}\right)  -\hat{\pi}_{m}\left(  y-1,x_{2}%
,\ldots,x_{d}\right)  \right)  ,
\]%
\begin{align*}
\chi_{m}\left(  z\right)   &  =2i\sum_{x_{2},\ldots,x_{d}}\exp\left[  \frac
{i}{m}\sum\nolimits_{s=2}^{d}x_{s}z_{s}\right] 
\\
&  
\times\sum_{y=-m}^{m}\left(  \hat{\pi}_{m}\left(  y,x_{2},\ldots
,x_{d}\right)  -\hat{\pi}_{m}\left(  y-1,x_{2},\ldots,x_{d}\right)  \right) \\
&  \times\frac{\mathrm{e}^{i\left(  z_{1}/m\right)  \left(  y-1/2\right)
}-\mathrm{e}^{i\left(  z_{1}/m\right)  \left(  m+1/2\right)  }}{\sin\left(
z_{1}/2m\right)  }.
\end{align*}
Therefore%
\[
\left\vert \chi_{m}\left(  z\right)  \right\vert \leq Cm^{d-1}\frac{m}{z_{1}%
}\sum_{y=-m}^{m}\left\vert \hat{\pi}_{m}\left(  y,x_{2},\ldots,x_{d}\right)
-\hat{\pi}_{m}\left(  y-1,x_{2},\ldots,x_{d}\right)  \right\vert \leq
C\frac{m^{14/15}}{\left\vert z\right\vert },
\]
and if $\alpha>1-\gamma$%
\begin{align}
\left\vert A_{5}\right\vert  &  \leq m^{-d}\int_{m^{\alpha}\leq\left\vert
z\right\vert }\left\vert \chi_{m}\left(  z\right)  \right\vert ^{n}dz\leq
C^{n}m^{-d}m^{14n/15}\int_{m^{\alpha}}^{\infty}r^{d-1}r^{-n}dz\label{Est_A5}\\
&  \leq C^{n}m^{-d}m^{14n/15}m^{\alpha\left(  d-n\right)  }\leq
2^{-n},\nonumber
\end{align}
if $m$ and $n$ are large enough.
Combining (\ref{Est_A1})-(\ref{Est_A5}), Proposition
\ref{Ap_Prop_LCLT} follows.
\end{proof}

We next need a simple large deviation estimate

\begin{lemma}
\label{Ap_Le_LDBound}There exists $C>0,$ such that for $\left\vert
x\right\vert \geq2m$,%
\[
\hat{\pi}_{m}^{\ast n}\left(  x\right)  \leq Cm^{-d}\exp\left[  -\frac
{\left\vert x\right\vert ^{2}}{Cnm^{2}}\right]  \,.
\]

\end{lemma}
\begin{proof}
[Proof of Lemma \ref{Ap_Le_LDBound}]If $\left\vert x\right\vert \geq r,$ then
one of the $d$ components of $x$ satisfies $\left\vert x_{i}\right\vert \geq
r/\sqrt{d}.$ By rotational symmetry, we get%
\[
\sum_{x:\left\vert x\right\vert \geq r}\hat{\pi}_{m}^{\ast n}\left(  x\right)
=dP\left(  \left\vert \sum\nolimits_{j=1}^{n}\xi_{j}\right\vert \geq
r/\sqrt{d}\right)  ,
\]
where the $\xi_{j}$ are i.i.d. with the one-dimensional marginal of $\hat{\pi
}$ as its distribution. Then,%
\[
P\left(  \left\vert \sum\nolimits_{j=1}^{n}\xi_{j}\right\vert \geq r/\sqrt
{d}\right)  \leq2\exp\left[  -nI\left(  \frac{r}{\sqrt{d}n}\right)  \right]
\,,\]
where%
\[
I\left(  t\right)  =\sup\left\{  \lambda t-\log E\left(  \mathrm{e}^{\lambda
t}\right)  \right\}  .
\]
By symmetry $I^{\prime}\left(  0\right)  =0,$ and from our assumptions, we
have $I^{\prime\prime}\left(  0\right)  \geq1/Cm^{2}.$ Furthermore, $I\left(
t\right)  =\infty$ if $\left\vert t\right\vert >2.$ By convexity of $I,$ we
therefore have $I\left(  t\right)  \geq t^{2}/Cm^{2}.$ Combining these
estimates gives%
\[
\sum_{x:\left\vert x\right\vert \geq r}\hat{\pi}_{m}^{\ast n}\left(  x\right)
\leq C\exp\left[  \frac{r^{2}}{Cnm^{2}}\right]  .
\]
From this, we get%
\begin{align*}
\hat{\pi}^{\ast n}\left(  x\right)   &  =\sum_{y}\hat{\pi}_{m}^{\ast\left(
n-1\right)  }\left(  y\right)  \hat{\pi}_{m}\left(  x-y\right) \\
&  \leq Cm^{-d}\sum_{y:\left\vert y\right\vert \geq\left\vert x\right\vert
-2m}\hat{\pi}_{m}^{\ast\left(  n-1\right)  }\left(  y\right)  \leq Cm^{-d}%
\exp\left[  -\frac{\left(  \left\vert x\right\vert -2m\right)  ^{2}}{C\left(
n-1\right)  m^{2}}\right] \\
&  \leq Cm^{-d}\exp\left[  -\frac{\left\vert x\right\vert ^{2}}{Cnm^{2}%
}\right]  .
\end{align*}

\end{proof}

Let%
\begin{equation}
G_{m}\left(  x\right)  \overset{\mathrm{def}}{=}\sum_{n=0}^{\infty}\hat{\pi
}_{m}^{\ast n}\left(  x\right)  . \label{Green_CGRW}%
\end{equation}

\begin{corollary}
\label{Ap_Cor_Green} For $\left\vert x\right\vert \geq m,$ we have for some
constant $c\left(  d\right)  $%
\[
G_{m}\left(  x\right)  =c\left(  d\right)  \frac{1}{\alpha\left(  m\right)
}\left\vert x\right\vert ^{-d+2}+O\left(  \left\vert x\right\vert ^{-d}\left(
\log\frac{\left\vert x\right\vert }{m}\right)  ^{5d}\right)  .
\]
For $\left\vert x\right\vert \leq m,$ we have%
\[
G_{m}\left(  x\right)  =\delta_{0,x}+O\left(  m^{-d}\right)
\,.\]
\end{corollary}
\begin{proof}
[Proof of Corollary \ref{Ap_Cor_Green}]Assume $\left\vert x\right\vert \geq m$
and set%
\[
N\left(  x,m\right)  \overset{\mathrm{def}}{=}\frac{\left\vert x\right\vert
^{2}}{\alpha\left(  m\right)  }\left(  \log\frac{\left\vert x\right\vert ^{2}%
}{\alpha\left(  m\right)  }\right)  ^{-10}.
\]
Then%
\begin{align*}
\sum_{n=N}^{\infty}\hat{\pi}^{\ast n}\left(  x\right)   &  =\sum_{n=N}%
^{\infty}\frac{1}{\left(  2\pi d\alpha\left(  m\right)  n\right)  ^{d/2}}%
\exp\left[  -\frac{\left\vert x\right\vert ^{2}}{2\alpha\left(  m\right)
n}\right] \\
&  +\sum_{n=N}^{\infty}O\left(  \alpha\left(  m\right)  ^{-d/2}n^{-\left(
d+2\right)  /2}\right)  .
\end{align*}%
and we note that
\[
\sum_{n=N}^{\infty}O\left(  \alpha\left(  m\right)  ^{-d/2}n^{-\left(
d+2\right)  /2}\right)  =O\left(  \left\vert x\right\vert ^{-d}\left(
\log\frac{\left\vert x\right\vert ^{2}}{\alpha\left(  m\right)  }\right)
^{5d}\right)
\,.\]
Putting%
$t_{n}\overset{\mathrm{def}}{=}{2\alpha\left(  m\right)  n}/{\left\vert
x\right\vert ^{2}}$,
we get%
\begin{align*}
&  \sum_{n=N}^{\infty}\frac{1}{\left(  2\pi d\alpha\left(  m\right)  n\right)
^{d/2}}\exp\left[  -\frac{\left\vert x\right\vert ^{2}}{2\alpha\left(
m\right)  n}\right] \\
&  =\frac{\left\vert x\right\vert ^{-d+2}}{2\left(  \pi d\right)  ^{d/2}%
\alpha\left(  m\right)  }\sum_{n=N}^{\infty}\frac{1}{\left(  t_{n}\right)
^{d/2}}\exp\left[  -\frac{1}{t_{n}}\right]  \left(  t_{n}-t_{n-1}\right) \\
&  =\frac{\left\vert x\right\vert ^{-d+2}}{2\left(  \pi d\right)  ^{d/2}%
\alpha\left(  m\right)  }\int_{0}^{\infty}t^{-d/2}\exp\left[  -t^{-1}\right]
dt+O\left(  \left\vert x\right\vert ^{-d}\right)  .
\end{align*}
This proves the corollary 
for $\left\vert x\right\vert \geq m$ with%
\[
c\left(  d\right)  =\frac{1}{2\left(  \pi d\right)  ^{d/2}}\int_{0}^{\infty
}t^{-d/2}\exp\left[  -t^{-1}\right]  dt.
\]
For $\left\vert x\right\vert \leq m,$ the estimate is evident from Proposition
\ref{Ap_Prop_LCLT}.
\end{proof}

\noindent
\begin{proof}
[Proof of Lemma \ref{Cor_Green} a)]There exists a $\theta,$ such that for any
$y\in{\operatorname*{Shell}\nolimits}_{L}(r(L))$,
%with $\left\vert y\right\vert \geq S\left(  r\right)  ,$
there exists a unit vector $x\in\mathbb{R}^{d}$ such that $\left(
y+C_{\theta}\left(  x\right)  \right)  \cap\partial V_{3r\left(  L\right)
}\left(  y\right)  \cap V_{L}=\emptyset.$ Using this, we see from
(\ref{ConeEst}), that our coarse grained Markov chain has after every visit of
${\operatorname*{Shell}\nolimits}_{L}\left(  r(L)\right)  $ a probability of
at least $\delta\left(  \theta\right)  $ to leave $V_{L}$ in the next step.
Therefore, the expected number of visits in this shell is finite, uniformly in
the starting point.
\end{proof}

\noindent
\begin{proof}
[Proof of Lemma \ref{Cor_Green} b)]If $x\in{\operatorname*{Shell}%
\nolimits}_{L}\left(  r,2s\right)  ,$ then $\hat{\pi}\left(  x,\cdot\right)  $
is an averaging over exit distributions from (discrete) balls $V_{u}\left(
x\right)  ,$ the averaging taken over $u$'s with $u\geq\left(  \gamma
/2\right)  d_{L} \left(  x\right)  .$ Therefore, there exists a $\delta>0,$
such that $\hat{\pi}\left(  x,{\operatorname*{Shell}\nolimits}_{L}\left(
d_{L}\left(  x\right)  \left(  1-\gamma/4\right)  \right)  \right)  \geq
\delta.$ Therefore, if $x\in{\operatorname*{Shell}\nolimits}_{L}\left(
a,a+\gamma/8\right)  ,$ $r(L)\leq a\leq2s(L),$ we have\newline$\hat{\pi
}\left(  x,{\operatorname*{Shell}\nolimits}_{L} \left(  a\left(
1-\gamma/8\right)  \right)  \right)  \geq\delta.$ Therefore, a Markov chain
with transition probabilities $\hat{\pi}$ which starts in
${\operatorname*{Shell}\nolimits}_{L}\left(  a,a+\gamma s(L)/8\right)  $ has
probability at least $\delta$ to reach in one step ${\operatorname*{Shell}%
\nolimits}_{L}\left(  a\left(  1-\gamma/8\right)  \right)  .$ By Lemma
\ref{Le_Lawler_Exit} c), an nearest neighbor chain starting in
${\operatorname*{Shell}\nolimits}_{L}\left(  a\left(  1-\gamma/8\right)
\right)  $ has a probability at least $\varepsilon\left(  \gamma\right)  >0$
of exiting $V_{L}$ before reentering into ${\operatorname*{Shell}%
\nolimits}_{L} \left(  a,a+\gamma/8\right)  .$ This evidently then applies
also to our coarse grained random walk.

We conclude that for the coarse grained chain starting in $x\in
{\operatorname*{Shell}\nolimits}_{L}\left(  a,a+\gamma s/8\right)  $, there is
a positive probability $\varepsilon>0,$ not depending on $x,a,$ that the chain
exits from $V_{L}$ before reentering this shell. It therefore follows that the
expected number of visits in ${\operatorname*{Shell}\nolimits}_{L}\left(
a,a+\gamma s/8\right)  $ is bounded, uniformly in the starting point of the
chain, and $a.$ From this the conclusion follows by summing over a finite
number of such shells.
\end{proof}

As a preparation for the proof of parts c) and d) of Lemma \ref{Cor_Green}, we
prove a preliminary result about our coarse grained random walk.

%We fix $L$ and the scheme $\mathcal{S}.$ We write $\hat{g}$ for $\hat
%{g}_{\mathcal{S},L}.$

\begin{lemma}
\label{Ap_Le_GreenFromBoundary}%
\[
\sup_{x\in{\operatorname*{Shell}\nolimits}_{L} \left(  2s(L)\right)  }%
\sum_{y\in V_{L-2s(L)}}\hat{g}_{L} \left(  x,y\right)  \leq C\left(  \log
L\right)  ^{3}.
\]

\end{lemma}
\begin{proof}
[Proof of Lemma \ref{Ap_Le_GreenFromBoundary}]The expression $\sum_{y\in
V_{L-2s(L)}}\hat{g}_{L}\left(  x,y\right)  $ is the expected total time that
the random walk spends in $V_{L-2s}\subset V_{L}.$ When starting in
${\operatorname*{Shell}\nolimits}_{L}\left(  2s(L)\right)  ,$ the walk has a
probability bounded from below, say by $\varepsilon_{1}>0,$ of never entering
$V_{L-2s(L)}$ before exiting $V_{L},$ uniformly in the starting point. If the
walk enters $V_{L-2s(L)},$ it has to enter through ${\operatorname*{Shell}%
\nolimits}_{L}\left(  2s,4s\right)  .$ Therefore%
\[
\sup_{x\in{\operatorname*{Shell}\nolimits}_{L}\left(  2s(L)\right)  }%
\sum_{y\in V_{L-2s}}\hat{g}_{L}\left(  x,y\right)  \leq\varepsilon_{1}%
^{-1}\left(  1+\sup_{x\in{\operatorname*{Shell}\nolimits}_{L} \left(
2s(L),4s(L)\right)  }E_{x}\left(  T_{{\operatorname*{Shell}\nolimits}%
_{L}\left(  2s(L)\right)  }^{\mathrm{CG}}\right)  \right)  ,
\]
where $T_{A}^{\mathrm{CG}}$ stands for the first entrance time into $A$ by the
coarse grained random walk with transition kernel $\hat{\pi}_{s}$ from
$V_{L-2s(L)}.$ It therefore suffices to prove%
\[
\sup_{x\in{\operatorname*{Shell}\nolimits}_{L}\left(  2s(L),4s(L)\right)  }
E_{x}\left(  T_{{\operatorname*{Shell}\nolimits}_{L}\left(  2s(L)\right)
}^{\mathrm{CG}}\right)  \leq C\left(  \log L\right)  ^{3},
\]

Consider the shells $R_{j}\overset{\mathrm{def}}{=} {\operatorname*{Shell}%
\nolimits}_{L}\left(  js(L),\left(  j+1\right)  s(L)\right)  ,$ $j\geq2,$ and
let $T_{j}$ be the first entrance time of our (coarse grained) random walk
into $R_{j}.$One then has
\[
P_{x}\left(  T_{R_{j}}^{\mathrm{CG}}<T_{{\operatorname*{Shell}\nolimits}%
_{L}\left(  2s(L) \right)  }^{\mathrm{CG}}\right)  \leq CP_{x}\left(
T_{R_{j}}^{\mathrm{RW}}<T_{ {\operatorname*{Shell}\nolimits}_{L}\left(
2s(L)\right)  }^{\mathrm{RW}}\right)  ,
\]
and the right hand side we can estimate by Lemma \ref{Le_Lawler_Exit} c),
giving%
\[
P_{x}\left(  T_{R_{j}}^{\mathrm{RW}}<T_{{\operatorname*{Shell}\nolimits}_{L}
\left(  2s(L)\right)  }^{\mathrm{RW}}\right)  \leq\frac{C}{j},
\]
and therefore we get%
\[
P_{x}\left(  T_{R_{j}}^{\mathrm{CG}}<T_{{\operatorname*{Shell}\nolimits}%
_{L}\left(  2s(L) \right)  }^{\mathrm{CG}}\right)  \leq\frac{C}{j}.
\]
If $x\in R_{j},$ we estimate the expected number of visits in $R_{j}$ by
Corollary \ref{Ap_Cor_Green}, which gives%
\[
\sup_{x\in R_{j}}\sum_{y\in R_{j}}G_{\gamma s(L)} \left(  x,y\right)  \leq
C\left(  \log L\right)  ^{3}.
\]
Combining these estimates completes the proof of Lemma
\ref{Ap_Le_GreenFromBoundary}
\end{proof}

Let $\sigma$ be the first entrance time of $\left\{  X_{n}\right\}  $ into
${\operatorname*{Shell}\nolimits}_{L}(2s(L))$.
%boundary layer $\mathcal{L}\overset{\mathrm{def}}{=}
%\partial_{\mathrm{in,}%
%2s\left(  L\right)  }\left(  V_{L}\right)  ,$
%and let $\mathcal{L}^{\prime
%}\overset{\mathrm{def}}{=}\partial_{\mathrm{in,}4s
%\left(  L\right)  }\left(
%V_{L}\right)  \ \backslash$ $\partial_{\mathrm{in,}2s\left(  L\right)
%}\left(  V_{L}\right)  .$
Before time $\sigma$,
%first entrance into $\mathcal{L},$
the Markov process $\{X_{n}\}$ proceeds as a random walk on $\mathbb{Z}^{d}$
with jump distribution $\hat{\pi}_{m},$ where $m=\gamma s\left(  L\right)  .$

\noindent

\begin{proof}
[Proof of Lemma \ref{Cor_Green} c), d), e)]From Corollary \ref{Ap_Cor_Green},
we get%
\[
\sup_{x\in V_{L}}\sum_{y\in V_{L-2s(L)}}
%\operatorname*{Bulk}_{L} }
G_{\gamma s(L)} \left(  x,y\right)  \leq C\left(  \log L\right)  ^{6}.
\]
Evidently, from Lemma \ref{Ap_Le_GreenFromBoundary}, we get%
\[
\sup_{x\in V_{L}}\sum_{y\in V_{L-2s(L)}}
%\operatorname*{Bulk}_{L}}
\left\vert G_{\gamma s(L)}\left(  x,y\right)  -\hat{g}_{L}\left(  x,y\right)
\right\vert \leq C\left(  \log L\right)  ^{3},
\]
which implies the statement d).

e) follows by the same approximation and%
\[
\sup_{x,x^{\prime}\in V:\left\vert x-x^{\prime}\right\vert \leq s}\sum
_{y\in\operatorname*{Bulk}_{L}}\left\vert G_{\gamma s(L)}\left(  x,y\right)
-G_{\gamma s(L)}\left(  x^{\prime},y\right)  \right\vert \leq C\left(  \log
L\right)  ^{3},
\]
which follows again from Corollary \ref{Ap_Cor_Green}.
%\texttt{More details here won't hurt..}

We turn to the proof of part c). For $x=y$, the result is obvious from the
transience of simple random walk. In the sequel, we thus always take $x\neq
y$.
%Recall that $\gamma<1/10$,
Write $A_{y} \overset{\mathrm{def}}{=} \{z: |z-y|\leq s(L)\}$. We first prove
the result for $x\in A_{y} $ and $d_{L}(y)\geq5s(L)$. In that case,
%Then,
%with
%$\sigma^{\prime}$ denoting the first entrance time of the simple random
%walk into ${\operatorname*{Shell}}_{L}(2s(L))$,%
\[
\sup_{x\in A_{y}: x\neq y} \hat{g}_{L}(x,y) \leq G_{\gamma s(L)}(x,y)+
\max_{z\in{\operatorname*{Shell}}_{L}(2s(L))} P_{z}^{\mathrm{R}W}(T_{A_{y}%
}<T_{V_{L}}) \sup_{x\in A_{y}: x\neq y} \hat{g}_{L}(x,y) \,.
\]
Since
\[
\max_{z\in{\operatorname*{Shell}}_{L}(2s(L))} P_{z}^{\mathrm{R}W}(T_{A_{y}%
}<T_{V_{L}})<1
\]
uniformly in $L$ by Donsker's invariance principle, we conclude that
\[
\sup_{x\in A_{y}: x\neq y} \hat{g}_{L}(x,y) \leq C G_{\gamma s(L)}(x,y)\,.
\]
Corollary \ref{Ap_Cor_Green} then completes the proof in this case.

Consider next $x\in A_{y}$ but $s(L)\leq d_{L}(y)\leq5s(L)$, and set $B_{y}
\overset{\mathrm{def}}{=} \{z: |z-y|\leq s(L)/2\}$ and $C_{y} \overset
{\mathrm{def}}{=} \{z: |z-y|\leq5s(L)\}$. We note that
\[
\sup_{x\in A_{y}: x\neq y} \hat{g}_{L}(x,y) \leq\frac{C}{s(L)^{d}}+
\sup_{x\not \in A_{y}} \hat{g}_{L}(x,y) \leq\frac{C}{s(L)^{d}}+ \sup
_{z\not \in A_{y}} P_{z}^{\mathrm{R}W}(T_{B_{y}}<T_{V_{L}}) \sup_{x\in A_{y}:
x\neq y} \hat{g}_{L}(x,y)\,.
\]
Since $\sup_{z\not \in A_{y}} P_{z}^{\mathrm{R}W}(T_{B_{y}}<T_{V_{L}})<1$
uniformly in $L$, again by Donsker's invariance principle, we conclude that
\[
\sup_{x\in A_{y}: x\neq y} \hat{g}_{L}(x,y) \leq\frac{C}{s(L)^{d}}\,,
\]
which proves the claim in this case.

We next consider $x\not \in A_{y}$. Let $\sigma^{\prime}$ denote the first
entrance time of the simple random walk into ${\operatorname*{Shell}}%
_{L}(2s(L))$. Clearly, $\sigma^{\prime}\leq\sigma$. We then have
\begin{align}
\hat{g}_{L}(x,y) &  \leq G_{\gamma s(L)}(x,y)+C
\!\!\!\!\!
\!\!\!\!\!
\sum_{z\in
{\operatorname*{Shell}}_{L}(2s(L))}
\!\!\!\!\!
P_{x}^{\mathrm{R}W}(X_{\sigma^{\prime}%
}=z)P_{x}^{\mathrm{R}W}(T_{A_{y}}<T_{V_{L}})
\!\!\!\!\!
\sup_{w\in A_{y}:w\neq y}\hat
{g}_{L}(w,y)\nonumber\\
&  \leq\frac{C}{s(L)^{2}|x-y|^{d-2}}+\frac{Cd_{L}(x)d_{L}(y)}{s(L)^{2}}%
\sum_{z\in{\operatorname*{Shell}}_{L}(2s(L))}\frac{1}{(|x-z|\vee
1)^{d}(|y-z|\vee1)^{d}}\label{eq-270605a}\\
&  \leq\frac{C}{s(L)^{2}|x-y|^{d-2}}+\frac{C}{s(L)^{2}}\sum_{z\in
{\operatorname*{Shell}}_{L}(2s(L))}\frac{1}{(|x-z|\vee1)^{d-1}(|y-z|\vee
1)^{d-1}}\nonumber\\
&  \leq\frac{C}{s(L)^{2}|x-y|^{d-2}}\,,\nonumber
\end{align}
where the second inequality uses Corollary \ref{Ap_Cor_Green}, the estimate on
$\hat{g}_{L}(x,y)$ for $x\in A_{y}$ that was already proved, and Lemma
\ref{Le_MainExit}. This completes the proof.
\end{proof}


\begin{thebibliography}{99}                                                                                               %


\bibitem {BSZ}E.~Bolthausen, E., Sznitman, A.S., and Zeitouni, O., Cut points
and diffusive random walks in random environment.
\newblock {\em Ann. I.H. Poincar\'e}, \textbf{PR 39} (2003), pp. 527--555.

\bibitem {BK}Bricmont, J. and Kupiainen, A., \newblock Random walks in
asymmetric random environments. \newblock {\em Comm. Math. Phys.},
\textbf{142} (1991), pp. 345--420.

\bibitem {gilbarg}Gilbarg, D. and Trudinger, N. S., \textit{Elliptic Partial
Differential Equations of Second Order}, 2nd Ed., Springer (2001).

\bibitem {krylov}Krylov, N. V., \textit{Lectures on Elliptic and Parabolic
Equations in H\"{o}lder Spaces}, Graduate Studies in mathematics, Vol. 12,
American Math. Soc. (1996).

\bibitem {Lawler}Lawler, G. F., \textit{Intersections of Random Walks },
Birkh\"{a}user 1991.

\bibitem {Ledoux}Ledoux, M.: \textit{The Concentration of Measure Phenomenon.
}American Math. Soc. (2001).
%\bibitem{sznitman1} Sznitman, A.S.,
%On a class of transient random walks in random environment,
%{\em Ann. Probab.} {\bf  29}  (2001),  pp.  724--765


\bibitem {sznitmanLN}Sznitman, A.S., \newblock Topics in random walk in random
environment. \newblock {\em \rm Notes of course at School and Conference
on Probability
Theory, May 2002, ICTP Lecture Series, Trieste, 203--266}, 2004.

\bibitem {sznitman1}Sznitman, A.S., An effective criterion for ballistic
behavior of random walks in random environment, \emph{Probab. Theory Related
Fields} \textbf{122} (2002), pp. 509--544.

\bibitem {sznitman2}Sznitman, A.S., On new examples of ballistic random walks
in random environment, \emph{Ann. Probab.} \textbf{31} (2003), pp. 285--322.

\bibitem {SZ}Sznitman, A.S. and Zeitouni, O., An invariance principle for
isotropic diffusions in random environment, to appear,
\textit{Invent. Math.} (2006).

\bibitem {zaitsev}Zaitsev, A. Yu., Estimates for the strong approximations in
multidimensional central limit theorem. Proceedings of ICM 2002,
\textit{Documenta Mathematica}, vol III (2002), pp. 107--116. Higher Ed.
Press, Beijing.

\bibitem {zeitouniLN}Zeitouni, O., \newblock Random walks in random
environment.
\newblock {\em Lecture Notes in Mathematics, {\rm Springer, Berlin}},
\textbf{1837} (2004), pp. 190--312.


\end{thebibliography}
\end{document}